\documentclass{elsarticle}
  
\usepackage[margin=1in]{geometry}

\usepackage{bbm}
\usepackage{graphicx}
\usepackage{pstool}
\usepackage{wrapfig}
\usepackage{caption}
\usepackage{subcaption}
\usepackage[section]{placeins} 

\usepackage{fancyhdr} 
\usepackage{lastpage} 
\usepackage{extramarks} 

\usepackage{xcolor} 
\usepackage{enumerate} 
\usepackage{paralist} 
\usepackage{amsmath, amsthm, amssymb, mathtools}
\usepackage{mathabx, pifont, stmaryrd} 
\usepackage[explicit]{titlesec} 
\usepackage{etoolbox} 
\usepackage{bibentry} 
\makeatletter\let\saved@bibitem\@bibitem\makeatother 
\usepackage[colorlinks, bookmarksopen, bookmarksnumbered,
            citecolor=red,urlcolor=red]{hyperref} 
\makeatletter\let\@bibitem\saved@bibitem\makeatother 

\usepackage{algorithm}
\usepackage{algorithmic}

\newtheorem{prop}{Proposition}
\newtheorem{remark}{Remark}

\theoremstyle{definition}

\usepackage{bm} 
\usepackage{amsfonts} 

\DeclareMathOperator*{\argmin}{arg\,min}

\newcommand{\func}[3]{\ensuremath{#1 : #2 \rightarrow #3}}
\newcommand{\spanV}{\ensuremath{\mathrm{span\hspace{0.05cm}}}}
\newcommand{\norm}[1]{\ensuremath{\left\| #1 \right\|}}

\newcommand{\suchthat}{\mathrel{}\middle|\mathrel{}}

\newcommand{\optunc}[2]{\underset{#1}{\text{minimize}} ~~ #2}

\newcommand{\pder}[2]{\ensuremath{\frac{\partial #1}{\partial #2}}}

\newcommand{\Dcal}{\ensuremath{\mathcal{D}}}
\newcommand{\Ecal}{\ensuremath{\mathcal{E}}}

\newcommand{\Gcal}{\ensuremath{\mathcal{G}}}
\newcommand{\Hcal}{\ensuremath{\mathcal{H}}}
\newcommand{\Ical}{\ensuremath{\mathcal{I}}}
\newcommand{\Jcal}{\ensuremath{\mathcal{J}}}

\newcommand{\Lcal}{\ensuremath{\mathcal{L}}}
\newcommand{\Mcal}{\ensuremath{\mathcal{M}}}

\newcommand{\Ocal}{\ensuremath{\mathcal{O}}}

\newcommand{\Tcal}{\ensuremath{\mathcal{T}}}

\newcommand{\Vcal}{\ensuremath{\mathcal{V}}}
\newcommand{\Wcal}{\ensuremath{\mathcal{W}}}

\newcommand{\Vboldcal}{\ensuremath{\boldsymbol{\mathcal{V}}}}

\newcommand{\Gbb}{\ensuremath{\mathbb{G}}}

\newcommand{\Nbb}{\ensuremath{\mathbb{N} }}
\newcommand{\Pbb}{\ensuremath{\mathbb{P} }}

\newcommand{\Rbb}{\ensuremath{\mathbb{R} }}

\newcommand\Fbm{{\ensuremath{\bm{F}}}}

\newcommand\Ibm{{\ensuremath{\bm{I}}}}
\newcommand\Jbm{{\ensuremath{\bm{J}}}}

\newcommand\Rbm{{\ensuremath{\bm{R}}}}

\newcommand\Ubm{{\ensuremath{\bm{U}}}}
\newcommand\Vbm{{\ensuremath{\bm{V}}}}
\newcommand\Wbm{{\ensuremath{\bm{W}}}}

\newcommand\cbm{{\ensuremath{\bm{c}}}}

\newcommand\vbm{{\ensuremath{\bm{v}}}}
\newcommand\wbm{{\ensuremath{\bm{w}}}}
\newcommand\xbm{{\ensuremath{\bm{x}}}}
\newcommand\ybm{{\ensuremath{\bm{y}}}}
\newcommand\zbm{{\ensuremath{\bm{z}}}}

\newcommand\chibold{{\ensuremath{\boldsymbol{\chi}}}}

\newcommand\etabold{{\ensuremath{\boldsymbol{\eta}}}}

\newcommand\Phibold{{\ensuremath{\boldsymbol{\Phi}}}}

\newcommand\Psibold{{\ensuremath{\boldsymbol{\Psi}}}}

\newcommand\zerobold{\ensuremath{\mathbf{0}}}

\usepackage{tikz}
\usepackage{pgfplots}
\usepackage{pgfplotstable, filecontents, booktabs}
\pgfplotsset{compat=1.9}

\usetikzlibrary{pgfplots.groupplots}
\usepgfplotslibrary{fillbetween}
\usetikzlibrary{calc,fit,matrix,arrows,automata,positioning,shapes}
\usetikzlibrary{arrows.meta}

\pgfplotsset{select coords between index/.style 2 args={
    x filter/.code={
        \ifnum\coordindex<#1\fi
        \ifnum\coordindex>#2\fi
    }
}}

\tikzset{
 invisible/.style={opacity=0},
 visible on/.style={alt={#1{}{invisible}}},
 alt/.code args={<#1>#2#3}{%
   \alt<#1>{\pgfkeysalso{#2}}{\pgfkeysalso{#3}}
 },
}

\newcommand{\colorbarMatlabParula}[5]{
\begin{tikzpicture}
\begin{axis}[
   hide axis, scale only axis,
   height=0pt, width=0pt,
   colormap={parula}{rgb255=(62,38,168) rgb255=(62,39,172) rgb255=(63,40,175) rgb255=(63,41,178) rgb255=(64,42,180) rgb255=(64,43,183) rgb255=(65,44,186) rgb255=(65,45,189) rgb255=(66,46,191) rgb255=(66,47,194) rgb255=(67,48,197) rgb255=(67,49,200) rgb255=(67,50,202) rgb255=(68,51,205) rgb255=(68,52,208) rgb255=(69,53,210) rgb255=(69,55,213) rgb255=(69,56,215) rgb255=(70,57,217) rgb255=(70,58,220) rgb255=(70,59,222) rgb255=(70,61,224) rgb255=(71,62,225) rgb255=(71,63,227) rgb255=(71,65,229) rgb255=(71,66,230) rgb255=(71,68,232) rgb255=(71,69,233) rgb255=(71,70,235) rgb255=(72,72,236) rgb255=(72,73,237) rgb255=(72,75,238) rgb255=(72,76,240) rgb255=(72,78,241) rgb255=(72,79,242) rgb255=(72,80,243) rgb255=(72,82,244) rgb255=(72,83,245) rgb255=(72,84,246) rgb255=(71,86,247) rgb255=(71,87,247) rgb255=(71,89,248) rgb255=(71,90,249) rgb255=(71,91,250) rgb255=(71,93,250) rgb255=(70,94,251) rgb255=(70,96,251) rgb255=(70,97,252) rgb255=(69,98,252) rgb255=(69,100,253) rgb255=(68,101,253) rgb255=(67,103,253) rgb255=(67,104,254) rgb255=(66,106,254) rgb255=(65,107,254) rgb255=(64,109,254) rgb255=(63,110,255) rgb255=(62,112,255) rgb255=(60,113,255) rgb255=(59,115,255) rgb255=(57,116,255) rgb255=(56,118,254) rgb255=(54,119,254) rgb255=(53,121,253) rgb255=(51,122,253) rgb255=(50,124,252) rgb255=(49,125,252) rgb255=(48,127,251) rgb255=(47,128,250) rgb255=(47,130,250) rgb255=(46,131,249) rgb255=(46,132,248) rgb255=(46,134,248) rgb255=(46,135,247) rgb255=(45,136,246) rgb255=(45,138,245) rgb255=(45,139,244) rgb255=(45,140,243) rgb255=(45,142,242) rgb255=(44,143,241) rgb255=(44,144,240) rgb255=(43,145,239) rgb255=(42,147,238) rgb255=(41,148,237) rgb255=(40,149,236) rgb255=(39,151,235) rgb255=(39,152,234) rgb255=(38,153,233) rgb255=(38,154,232) rgb255=(37,155,232) rgb255=(37,156,231) rgb255=(36,158,230) rgb255=(36,159,229) rgb255=(35,160,229) rgb255=(35,161,228) rgb255=(34,162,228) rgb255=(33,163,227) rgb255=(32,165,227) rgb255=(31,166,226) rgb255=(30,167,225) rgb255=(29,168,225) rgb255=(29,169,224) rgb255=(28,170,223) rgb255=(27,171,222) rgb255=(26,172,221) rgb255=(25,173,220) rgb255=(23,174,218) rgb255=(22,175,217) rgb255=(20,176,216) rgb255=(18,177,214) rgb255=(16,178,213) rgb255=(14,179,212) rgb255=(11,179,210) rgb255=(8,180,209) rgb255=(6,181,207) rgb255=(4,182,206) rgb255=(2,183,204) rgb255=(1,183,202) rgb255=(0,184,201) rgb255=(0,185,199) rgb255=(0,186,198) rgb255=(1,186,196) rgb255=(2,187,194) rgb255=(4,187,193) rgb255=(6,188,191) rgb255=(9,189,189) rgb255=(13,189,188) rgb255=(16,190,186) rgb255=(20,190,184) rgb255=(23,191,182) rgb255=(26,192,181) rgb255=(29,192,179) rgb255=(32,193,177) rgb255=(35,193,175) rgb255=(37,194,174) rgb255=(39,194,172) rgb255=(41,195,170) rgb255=(43,195,168) rgb255=(44,196,166) rgb255=(46,196,165) rgb255=(47,197,163) rgb255=(49,197,161) rgb255=(50,198,159) rgb255=(51,199,157) rgb255=(53,199,155) rgb255=(54,200,153) rgb255=(56,200,150) rgb255=(57,201,148) rgb255=(59,201,146) rgb255=(61,202,144) rgb255=(64,202,141) rgb255=(66,202,139) rgb255=(69,203,137) rgb255=(72,203,134) rgb255=(75,203,132) rgb255=(78,204,129) rgb255=(81,204,127) rgb255=(84,204,124) rgb255=(87,204,122) rgb255=(90,204,119) rgb255=(94,205,116) rgb255=(97,205,114) rgb255=(100,205,111) rgb255=(103,205,108) rgb255=(107,205,105) rgb255=(110,205,102) rgb255=(114,205,100) rgb255=(118,204,97) rgb255=(121,204,94) rgb255=(125,204,91) rgb255=(129,204,89) rgb255=(132,204,86) rgb255=(136,203,83) rgb255=(139,203,81) rgb255=(143,203,78) rgb255=(147,202,75) rgb255=(150,202,72) rgb255=(154,201,70) rgb255=(157,201,67) rgb255=(161,200,64) rgb255=(164,200,62) rgb255=(167,199,59) rgb255=(171,199,57) rgb255=(174,198,55) rgb255=(178,198,53) rgb255=(181,197,51) rgb255=(184,196,49) rgb255=(187,196,47) rgb255=(190,195,45) rgb255=(194,195,44) rgb255=(197,194,42) rgb255=(200,193,41) rgb255=(203,193,40) rgb255=(206,192,39) rgb255=(208,191,39) rgb255=(211,191,39) rgb255=(214,190,39) rgb255=(217,190,40) rgb255=(219,189,40) rgb255=(222,188,41) rgb255=(225,188,42) rgb255=(227,188,43) rgb255=(230,187,45) rgb255=(232,187,46) rgb255=(234,186,48) rgb255=(236,186,50) rgb255=(239,186,53) rgb255=(241,186,55) rgb255=(243,186,57) rgb255=(245,186,59) rgb255=(247,186,61) rgb255=(249,186,62) rgb255=(251,187,62) rgb255=(252,188,62) rgb255=(254,189,61) rgb255=(254,190,60) rgb255=(254,192,59) rgb255=(254,193,58) rgb255=(254,194,57) rgb255=(254,196,56) rgb255=(254,197,55) rgb255=(254,199,53) rgb255=(254,200,52) rgb255=(254,202,51) rgb255=(253,203,50) rgb255=(253,205,49) rgb255=(253,206,49) rgb255=(252,208,48) rgb255=(251,210,47) rgb255=(251,211,46) rgb255=(250,213,46) rgb255=(249,214,45) rgb255=(249,216,44) rgb255=(248,217,43) rgb255=(247,219,42) rgb255=(247,221,42) rgb255=(246,222,41) rgb255=(246,224,40) rgb255=(245,225,40) rgb255=(245,227,39) rgb255=(245,229,38) rgb255=(245,230,38) rgb255=(245,232,37) rgb255=(245,233,36) rgb255=(245,235,35) rgb255=(245,236,34) rgb255=(245,238,33) rgb255=(246,239,32) rgb255=(246,241,31) rgb255=(246,242,30) rgb255=(247,244,28) rgb255=(247,245,27) rgb255=(248,247,26) rgb255=(248,248,24) rgb255=(249,249,22) rgb255=(249,251,21) },
   colorbar horizontal,
   point meta min=#1, point meta max=#5,
   colorbar style={width=10cm, xtick={#1,#2,#3,#4,#5}}
]
\addplot [draw=none] coordinates {(0,0)};
\end{axis}
\end{tikzpicture}
}

\newcommand{\colorbarMatlabJet}[5]{
\begin{tikzpicture}
\begin{axis}[
   hide axis, scale only axis,
   height=0pt, width=0pt,
   colormap={jet}{rgb255=(0,0,131) rgb255=(0,0,135) rgb255=(0,0,139) rgb255=(0,0,143) rgb255=(0,0,147) rgb255=(0,0,151) rgb255=(0,0,155) rgb255=(0,0,159) rgb255=(0,0,163) rgb255=(0,0,167) rgb255=(0,0,171) rgb255=(0,0,175) rgb255=(0,0,179) rgb255=(0,0,183) rgb255=(0,0,187) rgb255=(0,0,191) rgb255=(0,0,195) rgb255=(0,0,199) rgb255=(0,0,203) rgb255=(0,0,207) rgb255=(0,0,211) rgb255=(0,0,215) rgb255=(0,0,219) rgb255=(0,0,223) rgb255=(0,0,227) rgb255=(0,0,231) rgb255=(0,0,235) rgb255=(0,0,239) rgb255=(0,0,243) rgb255=(0,0,247) rgb255=(0,0,251) rgb255=(0,0,255) rgb255=(0,4,255) rgb255=(0,8,255) rgb255=(0,12,255) rgb255=(0,16,255) rgb255=(0,20,255) rgb255=(0,24,255) rgb255=(0,28,255) rgb255=(0,32,255) rgb255=(0,36,255) rgb255=(0,40,255) rgb255=(0,44,255) rgb255=(0,48,255) rgb255=(0,52,255) rgb255=(0,56,255) rgb255=(0,60,255) rgb255=(0,64,255) rgb255=(0,68,255) rgb255=(0,72,255) rgb255=(0,76,255) rgb255=(0,80,255) rgb255=(0,84,255) rgb255=(0,88,255) rgb255=(0,92,255) rgb255=(0,96,255) rgb255=(0,100,255) rgb255=(0,104,255) rgb255=(0,108,255) rgb255=(0,112,255) rgb255=(0,116,255) rgb255=(0,120,255) rgb255=(0,124,255) rgb255=(0,128,255) rgb255=(0,131,255) rgb255=(0,135,255) rgb255=(0,139,255) rgb255=(0,143,255) rgb255=(0,147,255) rgb255=(0,151,255) rgb255=(0,155,255) rgb255=(0,159,255) rgb255=(0,163,255) rgb255=(0,167,255) rgb255=(0,171,255) rgb255=(0,175,255) rgb255=(0,179,255) rgb255=(0,183,255) rgb255=(0,187,255) rgb255=(0,191,255) rgb255=(0,195,255) rgb255=(0,199,255) rgb255=(0,203,255) rgb255=(0,207,255) rgb255=(0,211,255) rgb255=(0,215,255) rgb255=(0,219,255) rgb255=(0,223,255) rgb255=(0,227,255) rgb255=(0,231,255) rgb255=(0,235,255) rgb255=(0,239,255) rgb255=(0,243,255) rgb255=(0,247,255) rgb255=(0,251,255) rgb255=(0,255,255) rgb255=(4,255,251) rgb255=(8,255,247) rgb255=(12,255,243) rgb255=(16,255,239) rgb255=(20,255,235) rgb255=(24,255,231) rgb255=(28,255,227) rgb255=(32,255,223) rgb255=(36,255,219) rgb255=(40,255,215) rgb255=(44,255,211) rgb255=(48,255,207) rgb255=(52,255,203) rgb255=(56,255,199) rgb255=(60,255,195) rgb255=(64,255,191) rgb255=(68,255,187) rgb255=(72,255,183) rgb255=(76,255,179) rgb255=(80,255,175) rgb255=(84,255,171) rgb255=(88,255,167) rgb255=(92,255,163) rgb255=(96,255,159) rgb255=(100,255,155) rgb255=(104,255,151) rgb255=(108,255,147) rgb255=(112,255,143) rgb255=(116,255,139) rgb255=(120,255,135) rgb255=(124,255,131) rgb255=(128,255,128) rgb255=(131,255,124) rgb255=(135,255,120) rgb255=(139,255,116) rgb255=(143,255,112) rgb255=(147,255,108) rgb255=(151,255,104) rgb255=(155,255,100) rgb255=(159,255,96) rgb255=(163,255,92) rgb255=(167,255,88) rgb255=(171,255,84) rgb255=(175,255,80) rgb255=(179,255,76) rgb255=(183,255,72) rgb255=(187,255,68) rgb255=(191,255,64) rgb255=(195,255,60) rgb255=(199,255,56) rgb255=(203,255,52) rgb255=(207,255,48) rgb255=(211,255,44) rgb255=(215,255,40) rgb255=(219,255,36) rgb255=(223,255,32) rgb255=(227,255,28) rgb255=(231,255,24) rgb255=(235,255,20) rgb255=(239,255,16) rgb255=(243,255,12) rgb255=(247,255,8) rgb255=(251,255,4) rgb255=(255,255,0) rgb255=(255,251,0) rgb255=(255,247,0) rgb255=(255,243,0) rgb255=(255,239,0) rgb255=(255,235,0) rgb255=(255,231,0) rgb255=(255,227,0) rgb255=(255,223,0) rgb255=(255,219,0) rgb255=(255,215,0) rgb255=(255,211,0) rgb255=(255,207,0) rgb255=(255,203,0) rgb255=(255,199,0) rgb255=(255,195,0) rgb255=(255,191,0) rgb255=(255,187,0) rgb255=(255,183,0) rgb255=(255,179,0) rgb255=(255,175,0) rgb255=(255,171,0) rgb255=(255,167,0) rgb255=(255,163,0) rgb255=(255,159,0) rgb255=(255,155,0) rgb255=(255,151,0) rgb255=(255,147,0) rgb255=(255,143,0) rgb255=(255,139,0) rgb255=(255,135,0) rgb255=(255,131,0) rgb255=(255,128,0) rgb255=(255,124,0) rgb255=(255,120,0) rgb255=(255,116,0) rgb255=(255,112,0) rgb255=(255,108,0) rgb255=(255,104,0) rgb255=(255,100,0) rgb255=(255,96,0) rgb255=(255,92,0) rgb255=(255,88,0) rgb255=(255,84,0) rgb255=(255,80,0) rgb255=(255,76,0) rgb255=(255,72,0) rgb255=(255,68,0) rgb255=(255,64,0) rgb255=(255,60,0) rgb255=(255,56,0) rgb255=(255,52,0) rgb255=(255,48,0) rgb255=(255,44,0) rgb255=(255,40,0) rgb255=(255,36,0) rgb255=(255,32,0) rgb255=(255,28,0) rgb255=(255,24,0) rgb255=(255,20,0) rgb255=(255,16,0) rgb255=(255,12,0) rgb255=(255,8,0) rgb255=(255,4,0) rgb255=(255,0,0) rgb255=(251,0,0) rgb255=(247,0,0) rgb255=(243,0,0) rgb255=(239,0,0) rgb255=(235,0,0) rgb255=(231,0,0) rgb255=(227,0,0) rgb255=(223,0,0) rgb255=(219,0,0) rgb255=(215,0,0) rgb255=(211,0,0) rgb255=(207,0,0) rgb255=(203,0,0) rgb255=(199,0,0) rgb255=(195,0,0) rgb255=(191,0,0) rgb255=(187,0,0) rgb255=(183,0,0) rgb255=(179,0,0) rgb255=(175,0,0) rgb255=(171,0,0) rgb255=(167,0,0) rgb255=(163,0,0) rgb255=(159,0,0) rgb255=(155,0,0) rgb255=(151,0,0) rgb255=(147,0,0) rgb255=(143,0,0) rgb255=(139,0,0) rgb255=(135,0,0) rgb255=(131,0,0) rgb255=(128,0,0) },
   colorbar horizontal,
   point meta min=#1, point meta max=#5,
   colorbar style={width=10cm, xtick={#1, #2, #3, #4, #5}}
]
\addplot [draw=none] coordinates {(0,0)};
\end{axis}
\end{tikzpicture}
}

\usepackage{setspace}

\newbool{fastcompile}
\setbool{fastcompile}{false}

\newcommand{\pgrad}{\nabla}
\newcommand{\pflux}{f}
\newcommand{\pstvc}{u}
\newcommand{\pdom}{\Omega}
\newcommand{\pnrml}{n}
\newcommand{\rgrad}{\nabla_0}
\newcommand{\rflux}{F}
\newcommand{\rstvc}{U}
\newcommand{\rdom}{\Omega_0}
\newcommand{\rnrml}{N}
\newcommand{\param}{\mu}
\newcommand{\paramsp}{\Dcal}
\newcommand{\dommap}{\Gcal}
\newcommand{\dommapsp}{\Gbb}
\newcommand{\domdet}{g}
\newcommand{\domjac}{G}
\newcommand{\pcoord}{x}
\newcommand{\rcoord}{X}

\newcommand{\dpcoord}{\xbm}

\newcommand{\dstvc}{\Ubm}
\newcommand{\dres}{\Rbm}
\newcommand{\msh}{\Ecal}
\newcommand{\mshsz}{h}
\newcommand{\mshord}{q}
\newcommand{\rnode}{\hat{X}}
\newcommand{\pnode}{\hat{x}}

\begin{document}
\title{Model reduction of convection-dominated partial differential
       equations via optimization-based implicit feature tracking}

\author[rvt1]{Marzieh Alireza Mirhoseini\fnref{fn1}}
\ead{malireza@nd.edu}

\author[rvt1]{Matthew J. Zahr\fnref{fn2}\corref{cor1}}
\ead{mzahr@nd.edu}

\address[rvt1]{Department of Aerospace and Mechanical Engineering, University
               of Notre Dame, Notre Dame, IN 46556, United States}
\cortext[cor1]{Corresponding author}

\fntext[fn1]{Graduate Student, Department of Aerospace and Mechanical
             Engineering, University of Notre Dame}
\fntext[fn2]{Assistant Professor, Department of Aerospace and Mechanical
             Engineering, University of Notre Dame}


\begin{abstract}
This work introduces a new approach to reduce the computational cost of solving
partial differential equations (PDEs) with convection-dominated solutions:
\textit{model reduction with implicit feature tracking}. Traditional model reduction
techniques use an affine subspace to reduce the dimensionality of the solution manifold
and, as a result, yield limited reduction and require extensive training due to the
slowly decaying Kolmogorov $n$-width of convection-dominated problems.
The proposed approach circumvents the slowly decaying $n$-width limitation
by using a nonlinear approximation manifold systematically defined by
composing a low-dimensional affine space with a space of bijections
of the underlying domain.
Central to the implicit feature tracking approach is a residual
minimization problem over the reduced nonlinear manifold
that simultaneously determines the reduced coordinates
in the affine space and the domain mapping that minimize
the residual of the unreduced PDE discretization.
This is analogous to standard minimum-residual
reduced-order models, except instead of only minimizing
the residual over the affine subspace of PDE states, our
method enriches the optimization space to also include
admissible domain mappings. The nonlinear trial manifold is
constructed by using the proposed residual minimization formulation to
determine domain mappings that cause parametrized features
to align in a reference domain for a set of training parameters.
Because the feature is stationary in the reference domain, i.e., the convective
nature of solution removed, the snapshots are effectively compressed to define
an affine subspace. The space of domain mappings, originally constructed using
high-order finite elements, are also compressed in a way that ensures
the boundaries of the original domain are maintained.
Several numerical experiments are provided, including
transonic and supersonic, inviscid, compressible flows,
to demonstrate the potential of the method to yield
accurate approximations to convection-dominated problems
with limited training.

\end{abstract}

\maketitle

\section{Introduction}
\label{sec:intro}
Partial differential equations (PDEs) that model convection-dominated
phenomena often arise in engineering practice and scientific
applications, ranging from the study of high-speed, turbulent
flow over vehicles to wave propagation through solid media.
The solutions of these equations are characterized by local
features or disturbances, e.g., shock or contact discontinuities,
rarefaction waves, steep gradients, and vortical structures, that
propagate throughout the domain as time evolves or a system parameter
varies. Numerical methods to approximate these solutions
require stabilization and fine, usually adaptive, grids to
adequately resolve the local features, which lead to discretizations
with a large number of degrees of freedom. This makes many-query
analyses, e.g., optimization, uncertainty quantification, and
``what-if'' scenarios, computationally demanding and often infeasible.
Projection-based model reduction \cite{benner2015survey}
is a promising approach to substantially reduce the
computational expense associated with PDE simulations
and enable many-query analyses.


\subsection{Subspace approximation and Kolmogorov $n$-width}
\label{sec:intro:sub}
Projection-based model reduction methods search for an approximate
solution to a PDE in a low-dimensional (affine) subspace, usually
defined by the span of small number of data-driven modes with global
support constructed via the reduced basis method \cite{prud2002reliable}
or proper orthogonal decomposition (POD) combined with the
method of snapshots \cite{holmes1996turbulence,sirovich1987turbulence}.
The reduced subspace
is constructed in a time-intensive offline phase, where the
(expensive) PDE solution is queried at a number of sampled
parameter configurations and compressed (or orthogonalized)
to form an orthonormal basis. From the reduced solution (trial) space,
the reduced-order model (ROM) is defined via a projection of the
governing equations or residual minimization and serves
as a surrogate for the expensive PDE in the online stage. This
approach has lead to ROMs that are many orders of magnitude
faster to query than the unreduced PDE simulation, particularly
for elliptic and parabolic problems
\cite{veroy2003posteriori,grepl2005posteriori}. However, such an
approach is ``bound to fail'' \cite{ohlberger_reduced_2016} for
convection-dominated problems, although several authors have
applied these techniques to complex convection-dominated problems
\cite{carlberg_efficient_2011, washabaugh_use_2016, washabaugh_fast_2016,blonigan2021model},
usually in the non-parametric setting \cite{carlberg_efficient_2011} or
require extensive training \cite{washabaugh_use_2016, washabaugh_fast_2016}.

The Kolmogorov $n$-width associated with a subset $\Mcal\subset B$
of a Banach space $B$, denoted $d_n(\Mcal)$, defines the worst-case error
arising from the projection of points $u\in\Mcal$ onto the
best-possible linear space $V_n\subset B$ of a given dimension
$n\in\Nbb$
\begin{equation}
 d_n(\Mcal)\coloneqq
 \inf_{\substack{V_n\subset B\\\dim V_n=n}}\,
 \sup_{u \in \Mcal}\,
 \inf_{v_n\in V_n} \norm{u-v_n}_B.
\end{equation}
The decay of the $d_n(\Mcal)$ as $n\rightarrow\infty$ defines a fundamental
limitation on model reduction owing to the subspace approximation. A slowly
decaying $n$-width implies, even if the best-possible subspace is used, a
large space is needed to obtain small errors. This is a significant barrier to
effective model reduction, which relies on the dimension of the reduced
space being moderate to realize computational efficiency.
For elliptic and parabolic problems, the decay of the $n$-width is exponential,
i.e., $d_n(\Mcal) \leq C e^{-\beta n}$ where $C<\infty$ and $\beta > 0$ are
some constants \cite{buffa_priori_2012,ohlberger_reduced_2016}, which accounts
for the success of model reduction methods for this class of problems.
However, for hyperbolic equations, the $n$-width decay is at most $n^{-1/2}$
\cite{ohlberger_reduced_2016, greif_decay_2019}, which indicates that
convection-dominated features are not amenable to approximation via
low-dimensional linear spaces.

\begin{figure}
\centering
\input{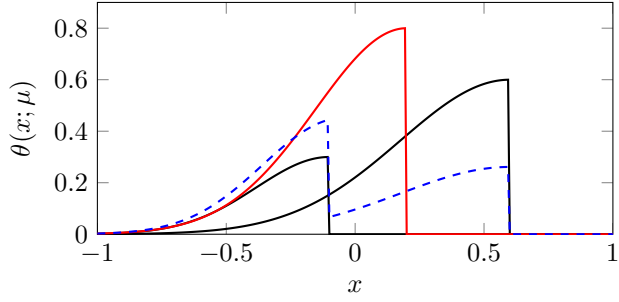}
\caption{Two instances of the cutoff Gaussian
         $\{\theta(\,\cdot\,;\param_i)\}_{i=1}^2$
         where $\param_1=(0.3,0.4,-0.1)$ and $\param_2=(0.6, 0.6, 0.6)$
         (\ref{line:demo1d0_snap}) are used to construct a two-dimensional
         linear space
         $\Vcal_2^\theta\coloneqq\spanV\{\theta(\,\cdot\,;\param_i)\}_{i=1}^2$,
         and the optimal approximation ($L^2$ sense)
         in $\Vcal_2^\theta$ (\ref{line:demo1d0_l2proj}) to a third
         instance of the cutoff Gaussian $\theta(\,\cdot\,;\param_3)$
         (\ref{line:demo1d0_test}) where $\param_3=(0.8,0.5,0.2)$.}
\label{fig:demo1d0_linapprox}
\end{figure}

To illustrate the poor approximation and slowly decaying $n$-width associated
with convection-dominated features, we consider a contrived example similar to
that in \cite{welper_interpolation_2017}.
Let $\func{\theta}{[-1,1]\times\Rbb^3}{\Rbb_{>0}}$ be the parametrized
cutoff Gaussian function defined as
\begin{equation} \label{eqn:gausscut}
 (\pcoord,(a,b,c)) \mapsto
 \theta(\pcoord;(a,b,c)) \coloneqq
 \begin{cases}
  a~\text{exp}\left[-\left(\frac{\pcoord-c}{b}\right)^2\right] & \pcoord < c, \\
  0 & \pcoord > c,
 \end{cases}
\end{equation}
where $a,b,c\in\Rbb$ define the amplitude, width, and center (and location
of discontinuity) of the Gaussian, respectively. We consider a typical
situation in model reduction whereby a two-dimensional linear space
$\Vcal_2^\theta$ is constructed as the span of two instances of the
function (snapshots), i.e.,
$\Vcal_2^\theta\coloneqq\spanV\{\theta(\,\cdot\,;\mu_i)\}_{i=1}^2$.
We arbitrarily choose $\mu_1=(0.3,0.4,-0.1)$ and $\mu_2=(0.6,0.6,0.6)$,
which corresponds to the discontinuity at $x=-0.1$ and $x=0.6$, respectively.
Then, for any element $u\in\Mcal=\{\theta(\,\cdot\,;\mu)\mid\mu\in\Rbb^3\}$,
the best approximation in $\Vcal_2^\theta$ ($L^2$ sense) will correspond to
a ``staircase'' approximation with discontinuities at $x=-0.1$ and $x=0.6$
(Figure~\ref{fig:demo1d0_linapprox}).
The staircase approximation will correspond to a large error (unless the
discontinuity in $u$ lies at or near $x=-0.1$ or $x=0.6$), which will only
slowly decrease as $\Vcal_2^\theta$ is expanded with additional snapshots.
This simple example is indicative of the poor approximations generally
obtained when applying model reduction to convection-dominated problems
and illustrates the slowly decaying Kolmogorov $n$-width.

\subsection{Proposed approach: nonlinear approximation via deforming modes,
            implicit feature tracking}
\label{sec:intro:nonl}
To avoid the fundamental linear reducibility limitation associated with
convection-dominated problems, we construct a nonlinear approximation by
composing a low-dimensional affine space with a parametrized domain mapping,
a concept first introduced in \cite{welper_interpolation_2017}.
The affine space is constructed using the method of snapshots and POD;
prior to compression, each snapshot is composed with a mapping that
causes its local features to align (same spatial location) with the
corresponding features in all other snapshots. The parametrized
domain mapping (defining the nonlinear approximation) is chosen
such that the local features present in the affine space deform
to the corresponding features in the solution being approximated,
a concept discussed in \cite{zahr_optimization-based_2018} in
the context of high-order methods for shock-dominated flows.
To demonstrate this concept,
we revisit the cutoff Gaussian example from the previous section.
Define the following parametrized family of bijections
$\func{\Hcal_\tau}{[-1,1]}{[-1,1]}$ for $\tau\in[-1,1]$ as
\begin{equation} \label{eqn:domdef_quad}
 \rcoord \mapsto \Hcal_\tau(\rcoord) \coloneqq \rcoord + \tau(1-X^2).
\end{equation}
and the corresponding mapped cutoff Gaussian function
$\func{\Theta_\tau}{[-1,1]\times\Rbb^3}{\Rbb_{>0}}$ as
\begin{equation}
 (\rcoord,\mu) \mapsto \Theta_\tau(\rcoord,\mu) \coloneqq
 \theta(\Hcal_\tau(X);\mu).
\end{equation}
For any $a,b\in\Rbb$ and $c\in[-1,1]$, if we choose the domain mapping
parameter as $\tau = c$, then all instances of $\Theta_c(\,\cdot\,;(a,b,c))$
have the discontinuity at $X=0$. In this setting, the affine space constructed
from two snapshots
$\Vcal_2^\Theta\coloneqq\spanV\{\Theta_{\tau_i}(\,\cdot\,;\mu_i)\}_{i=1}^2$
(same $\mu_1$ and $\mu_2$ from previous section and $\tau_i = (\mu_i)_3$)
effectively approximates the solution manifold
$\{\Theta_c(\,\cdot\,;(a,b,c))\mid a,b\in\Rbb,c\in[-1,1]\}$
because the convection-dominated nature of the solution has been removed
(Figure~\ref{fig:demo1d0_nlapprox}, left). A highly accurate approximation
for $\theta(x)= \Theta_\tau(\Hcal_\tau^{-1}(x);\mu)$ is obtained by
composing with the inverse of the domain mapping
(Figure~\ref{fig:demo1d0_nlapprox}, right). This concept of transformed
snapshots was first introduced in \cite{welper_interpolation_2017}.
\begin{figure}
\centering
\input{_py/demo1d0_nlapprox.tikz}
\caption{\textit{Left}: Two instances of the mapped cutoff Gaussian
  $\{\Theta_{\tau_i}(\,\cdot\,;\param_i)\}_{i=1}^2$ where
  $\param_1=(0.3,0.4,-0.1)$, $\tau_1=-0.1$,
  $\param_2=(0.6, 0.6, 0.6)$, $\tau_2=0.6$
  (\ref{line:demo1d0_snap}) used to construct a two-dimensional linear space
  $\Vcal_2^\Theta\coloneqq\spanV\{\Theta_{\tau_i}(\,\cdot\,;\param_i)\}_{i=1}^2$,
 and the optimal approximation ($L^2$ sense) in $\Vcal_2^\Theta$
  (\ref{line:demo1d0_l2proj}) to a third instance of mapped
  cutoff Gaussian $\Theta_{\tau_3}(\,\cdot\,;\param_3)$
  (\ref{line:demo1d0_test}) where $\param_3=(0.8,0.5,0.2)$, $\tau_3=0.2$.
  \textit{Right}: Push forward of the mapped functions by
  composing with the corresponding inverse mapping $\Hcal_{\tau}^{-1}$.}
\label{fig:demo1d0_nlapprox}
\end{figure}

In practice, an analytical expression for the domain mapping that causes
feature alignment will not be available; instead, we propose a formulation
that solves for the optimal domain mapping. To begin, we compose the PDE
solution with a generic domain mapping, reformulate the resulting equations
on a fixed domain with the mapping embedded into the PDE operators, and
apply a standard discretization. Then, we define a reduced-order model
as the solution in the nonlinear manifold, defined as the composition
of an affine space with a domain mapping, that minimizes the discretization
residual. This is a residual minimization approach; however, instead of
optimizing over an affine space (fixed domain), we optimize over the
nonlinear manifold induced by the domain mapping. In the remainder, we
refer to this as a reduced-order model with \textit{implicit feature tracking}
(ROM-IFT) because no attempt is made to explicitly quantify or identify local
features; they are tracked (in the sense of Figure~\ref{fig:demo1d0_nlapprox})
implicitly due to residual minimization over the nonlinear manifold,
which effectively removes the convection-dominated nature of the solution.
This is similar to the approach for implicitly tracking shocks with a
computational grid in
\cite{zahr_optimization-based_2018, corrigan_moving_2019, zahr_implicit_2020}.
To construct the affine space that generates the nonlinear manifold, we
use the implicit feature tracking framework, or any related approach
\cite{nair_transported_2019,taddei_registration_2020},
 to ensure features in the snapshots align (in the sense of
Figure~\ref{fig:demo1d0_nlapprox}) prior to compression via POD.
A low-dimensional family of domain mappings that preserve the boundary
of the domain, even for complex domains, are constructed using high-order
finite elements and POD compression.

\subsection{Connection to previous work}
\label{sec:intro:litrev}
Existing model reduction methods for convection-dominated flows that
aim to overcome the slowly decaying $n$-width can be broken into three
categories: localized affine approximations, online adaptation, nonlinear
approximations.

{\sc Localized linear (affine) approximation.}
The earliest generalizations of the affine approximation used in model
reduction define a collection of affine approximations, each one associated
with a particular region of the parameter domain \cite{haasdonk_training_2011},
time interval \cite{dihlmann_model_2011}, or state space
\cite{amsallem_nonlinear_2012, washabaugh_use_2016}. These methods benefit
from the localized nature of the approximation because the dimension of
each local affine space can be much smaller than the dimension of a
single global affine space, which can lead to significant savings
\cite{amsallem_nonlinear_2012, washabaugh_use_2016}.
However, as noted in \cite{ohlberger_reduced_2016}, the approximations
generated by these methods are limited by the affine space defined as
the union of all localized affine spaces, whose dimension is controlled
by the Kolmogorov $n$-width.

{\sc Online adaptation.} Another class of methods define a
low-dimensional affine approximation in the offline phase and
adapt it online based on error indicators to account
for e.g., insufficient training or a slowly decaying Kolmogorov $n$-width.
Such methods usually require sampling the HDM solution when the error
indicator exceeds a threshold, which breaks the offline-online decomposition.
They have been successfully used to accelerate optimization problems
\cite{arian_trust-region_2000, agarwal_trust-region_2013, yue_accelerating_2013, zahr_progressive_2015, zahr_efficient_2019, yano_globally_2020};
however, the efficiency of model reduction methods in more general
settings comes from an offline-online decomposition, where all expensive
operations are confined to the offline stage. To maintain online efficiency,
Constantine and Iaccarino \cite{constantine_reduced_2012} use the HDM to
correct the ROM in regions of the domain where discontinuities are likely to
be present, which only requires the HDM solution on a small region of the
domain. Similarly, Peherstorfer \cite{peherstorfer_model_2020} adapts the
affine approximation space with the HDM at selected locations in the spatial
domain. On the other hand, Lucia et. al. \cite{lucia_reduced_2003} use the
HDM in regions of the domain where moving discontinuities are present and
the ROM elsewhere. Taddei et. al. \cite{taddei_reduced_2015} decomposes
the space-time solution domain (scalar conservation laws) into regions where
the solution is smooth and approximates the space-time solution in each region
using an interpolation strategy.
Finally, Carlberg \cite{carlberg_adaptive_2015} avoids
directly re-sampling the HDM solution by splitting the reduced basis vectors
into multiple vectors with disjoint support, a form of $h$-refinement in the
context of model reduction.

{\sc Nonlinear trial manifolds.} Several model reduction methods
have been proposed that leverage nonlinear trial manifolds to overcome
the reduction limitation set by the Kolmogorov $n$-width. One popular
approach to construct such manifolds, and the approach adopted in this
work, use a nonlinear parametrization of an affine space, usually
a transformation of the underlying spatial domain
\cite{ohlberger_nonlinear_2013, welper_interpolation_2017, reiss_shifted_2018, rim_transport_2018, rim_displacement_2018, cagniart_model_2019, nair_transported_2019, black_projection-based_2020, taddei_registration_2020, welper_transformed_2020, torlo_model_2020, mojgani_physics-aware_2020, rim_manifold_2020, bansal2021model, ferrero2021registration, taddei2021registration},
as described in Section~\ref{sec:intro:nonl}. Many of these methods
limit the trial manifold by \textit{decoupling} the search for the
transformation and subspace approximation, whereas the proposed
approach simultaneously determines the subspace approximation
and mapping that minimize the HDM residual (Section~\ref{sec:rom:cmpr}
highlights the significance of this distinction). 
Furthermore, the domain mappings are often constructed in such a way
that they can only be applied to tensor-product domains.

The method of freezing \cite{ohlberger_nonlinear_2013}
decomposes the PDE solution as a shape and (Lie) group component and
introduces an online-efficient reduced basis approximation to solve
the resulting equations that freezes the mapping at a time instance
to evolve the shape component and subsequently updates the mapping.
Cagniart et. al. \cite{cagniart_model_2019} introduce a translational
mapping, choose the approximation to be orthogonal to the translating
space, and then choose the translation that fits the chosen subspace
approximation via residual minimization.
Transported snapshot model order reduction \cite{nair_transported_2019}
defines the transformation via a smooth transport field constructed
via low-order polynomial approximation that is trained via a
least-squares procedure; in the online phase, the transport
field is determined explicitly from the parameter configuration and
the subspace approximation is determined via residual minimization.
A registration method was proposed in \cite{taddei_registration_2020}
that aligns features in snapshots (in the sense of
Section~\ref{sec:intro:nonl}) via an optimization procedure that
minimizes the difference between the mapped solution and a reference
state; the procedure is generalized to the online setting via a regression
procedure and the subspace approximation is determined using standard reduced
basis techniques. In subsequent work \cite{taddei_space-time_2020},
the approach was extended to include hyperreduction based on empirical
quadrature \cite{yano_lp_2019}, a greedy procedure for determining sample
points, and a minimum-residual procedure to determine the subspace
approximation. Recently, Torlo \cite{torlo_model_2020} used an arbitrary
Lagrangian-Eulerian (ALE) formulation of the governing equations
and constructs mappings to explicitly track discontinuities;
regression and machine learning techniques are used to define
online-efficient procedure to generalize the mapping and the
subspace approximation is subsequently determined using standard
reduced basis techniques. To the authors' knowledge, the only
approach to date that simultaneously solves for the domain mapping
and subspace approximation is proposed in \cite{black_projection-based_2020}
that generalized the moving finite element method (MFEM) to the context of
model reduction  by minimizing the time-continuous MFEM residual and combining with
phase conditions to yield a coupled system of differential algebraic
equations for the domain mapping and subspace approximation.

Another method to construct nonlinear trial manifolds abandons the
subspace approximation and directly defines a nonlinear approximation
to the PDE state using  convolutional autoencoders
\cite{kashima_nonlinear_2016, hartman_deep_2017, omata_novel_2019,
      lee_model_2020, xu_multi-level_2020, maulik_reduced-order_2020,
      kim2020efficient}.
Such approaches have shown improvement compared to linear approximations
in a number of fluid flow applications with regard to dimensionality
reduction; however, they tend to incur high training costs and realization
of online efficiency remains an open problem.

{\sc Other approaches.} Other approaches that have been developed to
circumvent the slowly decaying Kolmogorov $n$-width of convection-dominated
problems that do not fit cleanly into these categories include:
a Lagrangian formulation of the governing equations
\cite{mojgani_lagrangian_2017},
definition of a time-dependent trial subspace using the concept
of approximate Lax pairs \cite{gerbeau_approximated_2014}, and 
combining shock-fitting methods with POD in the context of
shape optimization \cite{brooks_karhunenloeve_2004}.

\subsection{Contributions}
The primary contributions of the manuscript are:
\begin{itemize}
 \item We introduce and analyze a novel model reduction method
  based on implicit feature tracking that:
  \begin{inparaenum}[i)]
   \item uses a nonlinear manifold constructed as the composition of an
    affine space with a domain mapping and
   \item simultaneously determines the subspace approximation and
    mapping by minimizing the HDM residual, distinguishing it from
    other related approaches and generalizing the concept of
    minimum-residual reduced-order models (usually posed over
    an affine trial space) \cite{maday_blackbox_2002}.
  \end{inparaenum}
  We demonstrate the significance of this generalization via an
  illustrative example (Section~\ref{sec:rom:cmpr}) that compares
  the ROM-IFT method to (standard) residual minimization over an
  affine trial space. 
 \item We define an offline procedure for snapshot alignment
  (in the sense of Section~\ref{sec:intro:nonl}) that bootstraps
  the ROM-IFT method, unifying the offline and online algorithms.
  By solely relying on residual minimization, there is no need to
  explicitly identify or track discontinuous features as in
  \cite{brooks_karhunenloeve_2004, torlo_model_2020}
  or compute integrals over discontinuous functions as in
  \cite{welper_interpolation_2017, taddei_registration_2020}.
 \item We introduce a general procedure to parametrize
  mappings of complex domains using high-order finite elements;
  many existing transformation-based
  methods are limited to tensor-product-like domains (required to
  preserve the boundary of the domain as discussed in
  \cite{zahr_implicit_2020, taddei_registration_2020} and
  Section~\ref{sec:govern:dommap:bnd}) due to the parametrization of the
  domain mapping.
 \item We detail a globally convergent Levenberg-Marquardt solver for
  the ROM-IFT residual minimization problem.
 \item We provide numerous numerical
  experiments to thoroughly investigate the performance
  of the proposed method on one- and two-dimensional benchmark
  problems from computational fluid dynamics including an
  inviscid, compressible (supersonic) flow over a cylinder.
\end{itemize}

\subsection{Paper outline}
The remainder of the paper is organized as follows.
Section~\ref{sec:govern} introduces the governing system of
parametrized partial differential equations, its transformation
to a fixed reference domain with the domain mapping embedded in
the definition of the transformed PDE operators (leading to a
transformation-dependent solution), and the corresponding discretization.
Section~\ref{sec:govern:dommap} introduces a procedure to discretize and
parametrize domain mappings using high-order meshes in such a way
that ensures the boundary of the domain is preserved (to high-order accuracy). 
Section~\ref{sec:rom} introduces the proposed nonlinear approximation
manifold, the ROM-IFT method that approximates the PDE solution with the
function in this manifold that minimizes the discrete residual, and
some theoretical results. A direct comparison between the ROM-IFT
and traditional residual minimization is included to provide intuition
for the ROM-IFT method and support the theoretical results. An offline
procedure for constructing feature-aligned snapshots to define an affine
space that generates the nonlinear manifold is introduced.
Section~\ref{sec:solver}
introduces a Levenberg-Marquardt method for solving the ROM-IFT problem
and Section~\ref{sec:numexp} demonstrates the merits of the
proposed model reduction framework for a number of parametrized,
convection-dominated PDEs. Finally, Section~\ref{sec:concl}
concludes the paper and identifies avenues of future research.


\section{Governing equations and discretization}
\label{sec:govern}
In this section, we introduce the governing partial differential equations
(Section~\ref{sec:govern:claw}), its transformation to a reference domain
so that a domain mapping appears explicitly in the governing equations
(Section~\ref{sec:govern:tclaw}), a discretization of the transformed
conservation law to yield the fully discrete governing equations
(Section~\ref{sec:govern:disc}), and a general family of finite-dimensional
domain mappings (Section~\ref{sec:govern:dommap}).

\subsection{Parametrized system of conservation laws}
\label{sec:govern:claw}
Consider a $\param$-parametrized system of $m$ conservation laws in
$d$-dimensions of the form
\begin{equation} \label{eqn:claw-phys}
 \pgrad\cdot\pflux(\pstvc,\pgrad\pstvc;\param) = s(\pstvc,\pgrad\pstvc;\param)
 \quad\text{in}\quad\pdom,
\end{equation}
where $\param\in\paramsp$ is the parameter and $\Dcal$ the parameter domain,
$\pstvc(\pcoord;\param)\in\Rbb^m$ is the solution (assumed unique) at a point
$\pcoord\in\pdom$,
$\pgrad\coloneqq[\partial_{\pcoord_1},\cdots,\partial_{\pcoord_d}]$
is the gradient operator on the domain $\pdom\subset\Rbb^d$ such that
$\pgrad z\coloneqq [\partial_{\pcoord_1}z \cdots \partial_{\pcoord_d}z]$,
$\func{\pflux}{\Rbb^m\times\Rbb^{m\times s}\times\Dcal}{\Rbb^{m\times d}}$
is the flux function, and
$\func{s}{\Rbb^m\times\Rbb^{m\times s}\times\Dcal}{\Rbb^m}$
is the source term. The boundary of the domain is $\partial\pdom$
with outward unit normal $\func{\pnrml}{\partial\pdom}{\Rbb^d}$. The
formulation of the conservation law in (\ref{eqn:claw-phys}) is
sufficiently general to encapsulate steady second-order partial differential
equations (PDEs) in  a $d$-dimensional spatial domain or time-dependent
PDEs in a $(d-1)$-dimensional domain, i.e., a $d$-dimensional space-time
domain.

\subsection{Transformed system of conservation laws on a fixed reference domain}
\label{sec:govern:tclaw}
\begin{figure}
 \centering
 \includegraphics[width=0.32\textwidth]{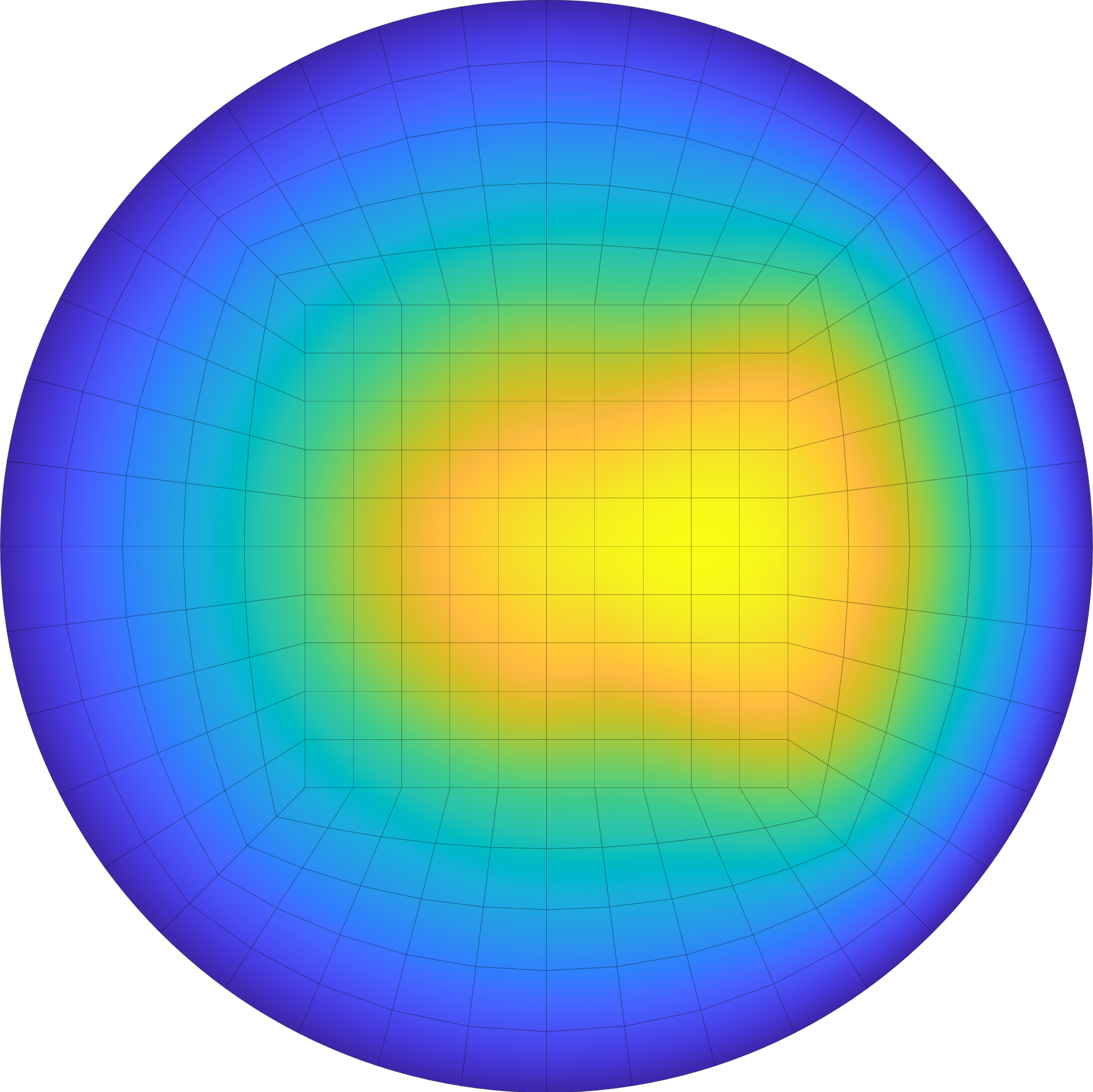} \qquad
 \includegraphics[width=0.32\textwidth]{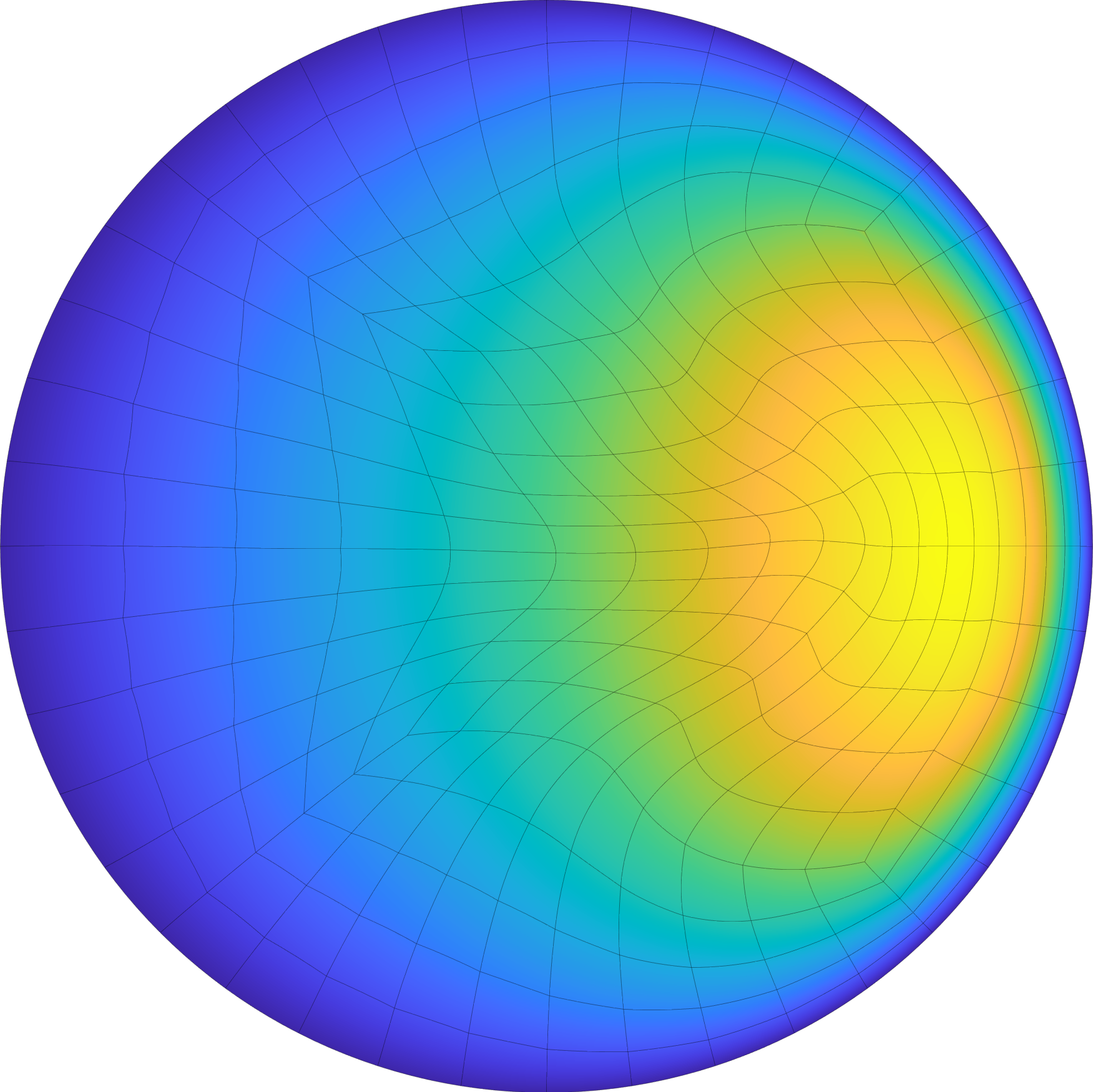}
\caption{A non-trivial domain mapping $\dommap\in\dommapsp$ (visualized
         by the edges of the mesh) in the case where $\rdom=\pdom$ (disk)
         and the corresponding solution on the reference domain
         $\rstvc(\rcoord;\dommap,\param)$ (\textit{left}) and physical
         domain $\pstvc(\pcoord;\param)$ (\textit{right}) for the
         advection-diffusion equation with homogeneous essential
         boundary conditions (Section~\ref{sec:rom:cmpr}).
         The domain mapping causes the feature
         in the solution to move in the reference domain relative to
         its actual position in the physical domain.}
\label{fig:refmap_ex}
\end{figure}

As the proposed model reduction method is fundamentally based on
deforming the domain, it is convenient to recast the PDE
(\ref{eqn:claw-phys}) on a fixed \textit{reference} domain,
$\rdom\coloneqq\bar\dommap^{-1}(\pdom)\subset\Rbb^d$,
following the approach in \cite{persson2009discontinuous}, where
$\func{\bar\dommap}{\Rbb^d}{\Rbb^d}$ is a smooth invertible mapping.
Let $\dommapsp$ be any collection of bijections from the reference
domain $\rdom$ to the physical domain $\pdom$ (Figure~\ref{fig:mapping}).
Then, for any $\dommap\in\dommapsp$, (\ref{eqn:claw-phys}) can be
written as a PDE on the reference domain as
\begin{equation} \label{eqn:claw-ref}
 \rgrad\cdot\rflux(\rstvc,\rgrad\rstvc;\dommap,\param)=S(\rstvc,\rgrad\rstvc;\dommap,\param)
 \quad\text{in}\quad\rdom,
\end{equation}
where $\rstvc(\rcoord;\dommap,\param)$ is the solution of the transformed PDE
at a point $\rcoord\in\rdom$,
$\rgrad\coloneqq[\partial_{\rcoord_1},\cdots,\partial_{\rcoord_d}]$
is the gradient operator on the domain $\rdom\subset\Rbb^d$ such that
$\rgrad z\coloneqq [\partial_{\rcoord_1}z \cdots \partial_{\rcoord_d}z]$,
$\func{\rflux}{\Rbb^m\times\Rbb^{m\times d}\times\Gbb\times\Dcal}{\Rbb^{m \times d}}$ is the transformed flux function, and
$\func{S}{\Rbb^m\times\Rbb^{m\times d}\times\Gbb\times\Dcal}{\Rbb^m}$
is the transformed source term. The boundary of the domain is $\partial\rdom$
with outward unit normal $\func{\rnrml}{\partial\rdom}{\Rbb^d}$.
\begin{figure}
\centering
\input{_py/mapping0.tikz}
\caption{Schematic of domain mapping for $\dommap\in\dommapsp$}
\label{fig:mapping}
\end{figure}
The reference and physical solutions and their gradients are related as
\begin{equation}\label{eqn:p2r-stvc}
 \rstvc(\rcoord;\dommap,\param) =
 \pstvc(\dommap(\rcoord);\param), \qquad
 \rgrad\rstvc(\rcoord;\dommap,\param) =
 \pgrad\pstvc(\dommap(\rcoord);\param)\domjac_\dommap(\rcoord),
\end{equation}
where $\domjac_\dommap\coloneqq\rgrad\dommap$ is the mapping
Jacobian and $\domdet_\dommap\coloneqq\det\domjac_\dommap$ is its
determinant. The reference and physical flux functions and source
terms are related as
\begin{equation}\label{eqn:p2r-flux}
 \rflux(w,\rgrad w;\dommap,\param) =
 \domdet_\dommap
 \pflux(w, \rgrad w \cdot\domjac_\dommap^{-1};\param)
 \domjac_\dommap^{-T}, \qquad
 S(w,\rgrad w;\dommap,\param) =  \domdet_\dommap s(w, \rgrad w \cdot\domjac_\dommap^{-1};\param)
\end{equation}
where $\func{w}{\rdom}{\Rbb^m}$ is any $m$-valued function over the reference
domain; see \cite{persson2009discontinuous} for details of derivation.
The conservation laws in
(\ref{eqn:claw-phys}) and (\ref{eqn:claw-ref}) are equivalent, that is,
if $\pstvc$ is the solution of (\ref{eqn:claw-phys}) and $\rstvc$
satisfies (\ref{eqn:p2r-stvc}), then $\rstvc$ is the solution of
(\ref{eqn:claw-ref}), and vice versa. Owing to this equivalence,
the PDE solution $\pstvc(\,\cdot\,;\param)$ is independent of
the domain mapping $\dommap\in\dommapsp$. Figure~\ref{fig:refmap_ex}
provides an example of a non-trivial domain mapping and the difference
between the solution on the reference and physical domains.

\subsection{Discretization of the transformed conservation law}
\label{sec:govern:disc}
The transformed PDE (\ref{eqn:claw-ref}) is discretized using
a suitable method depending on the properties of the flux
function $\pflux$ to yield a nonlinear system of algebraic equations
\begin{equation} \label{eqn:hdm}
 \dres(\dstvc;\dommap,\param) = \zerobold,
\end{equation}
where $\dstvc(\dommap,\param)\in\Rbb^N$ is the (assumed unique) solution
of the discretized PDE for $\dommap\in\dommapsp$, $\param\in\paramsp$,
and $\func{\dres}{\Rbb^N\times\Gbb\times\Dcal}{\Rbb^N}$ is the
nonlinear residual function defined by the spatial discretization of
(\ref{eqn:claw-ref}). Throughout the document, (\ref{eqn:hdm}) will
be referred to as the high-dimensional model (HDM).

To maintain a connection between the algebraic
solution and the spatial approximation to the PDE solution, let
$\Vcal_h$ denote the $N$-dimensional trial space (consisting of
functions over $\rdom$) associated with the discretization in
(\ref{eqn:hdm}) and let the operator $\func{\Upsilon}{\Rbb^N}{\Vcal_h}$
map a vector $\Wbm\in\Rbb^N$ to its representation as
a function over the reference domain $\rdom$,
$W_h=\Upsilon(\Wbm)\in\Vcal_h$, i.e.,
$\Wbm\in\Rbb^N$ is the vector encoding of $W_h$, a function over
$\rdom$. Thus, the $N$-vector $\dstvc(\dommap,\param)$ that satisfies
(\ref{eqn:hdm}) maps to the corresponding function $\rstvc_h\in\Vcal_h$
over the reference domain $\rdom$ as
\begin{equation} \label{eqn:hdm1}
 \rstvc_h(\,\cdot\,;\dommap,\param) \coloneqq \Upsilon(\dstvc(\dommap,\param)),
\end{equation}
for any $\dommap\in\dommapsp$ and $\param\in\paramsp$,
which in turn maps to a function over the physical domain $\pdom$
according to (\ref{eqn:p2r-stvc}), i.e.,
\begin{equation}
 \pstvc_h(\,\cdot\,;\dommap,\param) \coloneqq
 \rstvc_h(\,\cdot\,;\dommap,\param) \circ \dommap^{-1} \in
 \Vcal_h\circ\dommap^{-1},
\end{equation}
where, for any function space $\Vcal$ over $\rdom$, we define the
corresponding function space over $\pdom$ as
\begin{equation}
 \Vcal\circ\dommap^{-1}
 \coloneqq \left\{V\circ\dommap^{-1} \suchthat V\in\Vcal\right\}.
\end{equation}

\begin{remark} \label{rem:hdm-resolved}
 Recall, in the infinite-dimensional setting, the solution
 $\pstvc$ of (\ref{eqn:claw-phys}) is independent of the domain
 mapping $\dommap\in\dommapsp$. However, in the finite-dimensional
 setting, the solution in the physical domain $\pstvc_h$ will
 depend on the domain mapping. This is because the domain
 mapping deforms the solution in the reference domain causing
 its finite-dimensional approximation to change with $\dommap$;
 this is the essence of $r$-adaptivity (domain deformation
 to optimally place nodes for enhanced resolution)
 \cite{budd2009adaptivity}.
 Therefore, in Section~\ref{sec:numexp}, we test the accuracy
 of the proposed ROM-IFT method with respect to the exact
 solution (when available), rather than to the HDM solution on the same mesh.
\end{remark}

\begin{remark}
 The formulation of the HDM in this section is sufficiently general to apply
 to standard discretizations that fix the domain, i.e., $\dommap=\bar\dommap$,
 as well as $r$-adaptive methods such as those in
 \cite{zahr_optimization-based_2018, corrigan_moving_2019, zahr_implicit_2020}
 that simultaneously solve for the PDE state and mesh coordinates.
 In the $r$-adaptation case, only a slight reformulation is
 required where the HDM solution consists of the PDE state
 $\Ubm(\mu)$ and the domain mapping $\Gcal(\mu)$ for a given
 $\mu\in\Dcal$. In Section~\ref{sec:numexp}, we consider both
 types of discretizations.
\end{remark}

\subsection{Domain mapping discretization and parametrization}
\label{sec:govern:dommap}
Finally, we close this section by introducing a general family of
domain mappings using a computational mesh
(Section~\ref{sec:govern:dommap:disc}) and a parametrization that
ensures the mapping preserves the boundaries of the physical domain
(Section~\ref{sec:govern:dommap:bnd}).

\subsubsection{Discretization on mesh}
\label{sec:govern:dommap:disc}
Let $\msh_{\mshsz,\mshord}$ represent a discretization of the reference domain
$\rdom$ into non-overlapping, potentially curved, computational elements,
where $\mshsz$ is a mesh element size parameter and $\mshord\geq 1$ is
the polynomial degree associated with the curved elements, and let
$\{\rnode_I\}_{I=1}^{N_\mathrm{v}}\subset\Rbb^d$ denote the (ordered)
nodes of the mesh. Next, define the corresponding globally continuous,
piecewise polynomial function space as
\begin{equation*}
 \Wcal_{\mshsz,\mshord} =
 \left\{v \in C^0(\rdom) \suchthat
        \left.v\right|_K \in \Pbb_\mshord(K),
        ~\forall K \in \msh_{\mshsz,\mshord}\right\},
\end{equation*}
where $\Pbb_p(K)$ is the space of polynomial functions of degree at most
$p \geq 1$ on set $K\subset\rdom$. Any element of
$\dommapsp_{\mshsz,\mshord}\coloneqq[\Wcal_{\mshsz,\mshord}]^d$
is uniquely determined by its action on the nodes $\rnode_I$ of the
reference mesh. That is,
let $\{\Psi_I\}_{I=1}^{N_\mathrm{v}}$ be a nodal basis of $\Wcal_{\mshsz,\mshord}$
associated with the nodes $\{\rnode_I\}_{I=1}^{N_\mathrm{v}}$, i.e.,
$\Psi_I(\rnode_J)=\delta_{IJ}$ and define
\begin{equation}
 \dpcoord \mapsto
 \dommap_{\mshsz,\mshord}(\,\cdot\,;\dpcoord)\coloneqq\sum_{I=1}^{N_\mathrm{v}} \pnode_I \Psi_I,
\end{equation}
where the coefficients $\pnode_I\in\Rbb^d$ coincide with the action of
the mapping at the nodes $\rnode_I$ due to the choice of nodal basis,
which can be interpreted as the nodes of a mesh of $\pdom$,
and $\dpcoord\in\Rbb^{d N_\mathrm{v}}$ is the concatenation of
$\{\pnode_I\}_{I=1}^{N_\mathrm{v}}$. Figure~\ref{fig:refmap_ex} contains
the domain mapping induced by a particular choice of
$\{\pnode_I\}_{I=1}^{N_\mathrm{v}}$ for the case $\rdom=\pdom$.
The function space $\dommapsp_{\mshsz,\mshord}$
is not directly suitable for the domain mapping space $\dommapsp$
because, for a general configuration of coefficients $\dpcoord\in\Rbb^{d N_\mathrm{v}}$,
$\dommap_{\mshsz,\mshord}(\rdom;\dpcoord)\neq\pdom$.
In the next section, we will construct a subset
$\dommapsp_{\mshsz,\mshord}^\text{b}\subset\dommapsp_{\mshsz,\mshord}$
that ensures $\dommap(\rdom)=\pdom$ for any
$\dommap\in\dommapsp_{\mshsz,\mshord}^\text{b}$.

\begin{remark}
The choice to represent the domain mapping using local basis functions
associated with the mesh $\msh_{\mshsz,\mshord}$ is a
departure from similar approaches to model reduction that use
global basis functions such as Fourier modes \cite{mojgani_lagrangian_2017}
or radial basis functions \cite{taddei_registration_2020}.
This choice produces a broader class of domain
mappings, which is used in Section~\ref{sec:govern:dommap:bnd} to enforce the
constraint $\dommap(\rdom)=\pdom$ for complex domains.
\end{remark}


\subsubsection{Enforcement of physical boundaries}
\label{sec:govern:dommap:bnd}
In the finite-dimensional setting (Section~\ref{sec:govern:dommap:disc}), i.e.,
$\dommap\in\dommapsp_{\mshsz,\mshord}$, the condition
$\dommap(\rdom)=\pdom$ can only be ensured to accuracy $\Ocal(h^{q+1})$
owing to the resolution of the mesh $\Ecal_{h,q}$;
the condition will be exact for boundaries that are polynomial surfaces
of degree $\leq q$ in $\Rbb^d$, which always includes straight-sided
boundaries ($q\geq 1$). We assume that the reference and physical
domains are similar in the sense
that $\pdom$ is the image of $\rdom$ under a smooth bijection,
which implies any non-smooth feature (corners, kinks) in $\partial\pdom$
map from a corresponding feature in $\partial\rdom$. To this end,
let the boundaries $\partial\rdom$ and $\partial\pdom$ be defined
as the union of $N_\mathrm{b}$ smooth surfaces
$\{\partial\pdom_{0,i}\}_{i=1}^{N_\mathrm{b}}\subset\Rbb^d$ and
$\{\partial\pdom_i\}_{i=1}^{N_\mathrm{b}}\subset\Rbb^d$
\begin{equation}
 \partial\rdom = \bigcup_{i=1}^{N_\mathrm{b}} \partial{\pdom_{0,i}}, \qquad
 \partial\pdom = \bigcup_{i=1}^{N_\mathrm{b}} \partial{\pdom_i},
\end{equation}
which potentially intersect to form non-smooth features.
In all our numerical experiments (Section~\ref{sec:numexp}), we
take $\rdom=\pdom$ for simplicity and to guarantee this notion
of similarity between the domains.

To define a family of mappings $\dommapsp_{\mshsz,\mshord}^\text{b}$
that satisfy $\dommap(\rdom)=\pdom$ for any
$\dommap\in\dommapsp_{\mshsz,\mshord}^\text{b}$
(with accuracy $\Ocal(h^{q+1})$) in the finite-dimensional setting,
it is sufficient to ensure that
\begin{inparaenum}[1)]
 \item any node $\rnode_I\in\partial\pdom_{0,i}$
   gets mapped to $\pnode_I\in\partial\pdom_i$ for $i\in\Ical_I$,
   where $\Ical_I$ is the set of boundaries on which $\rnode_I$
   lies and
 \item the mapping is invertible.
\end{inparaenum}
Together these conditions imply that 
the mapped nodes $\{\pnode_I\}_{I=1}^{N_\mathrm{v}}$ combined with
the connectivity of $\Ecal_{\mshsz,\mshord}$ define a mesh of $\pdom$
(to order $\Ocal(\mshsz^{\mshord+1})$) and therefore
$\dommap(\rdom)=\pdom$ to the same order. The invertibility condition
will not be imposed directly on the function space; rather, it will
be enforced during the ROM-IFT optimization procedure through a
penalization term.

\begin{figure}
\centering
\includegraphics[width=0.48\textwidth]{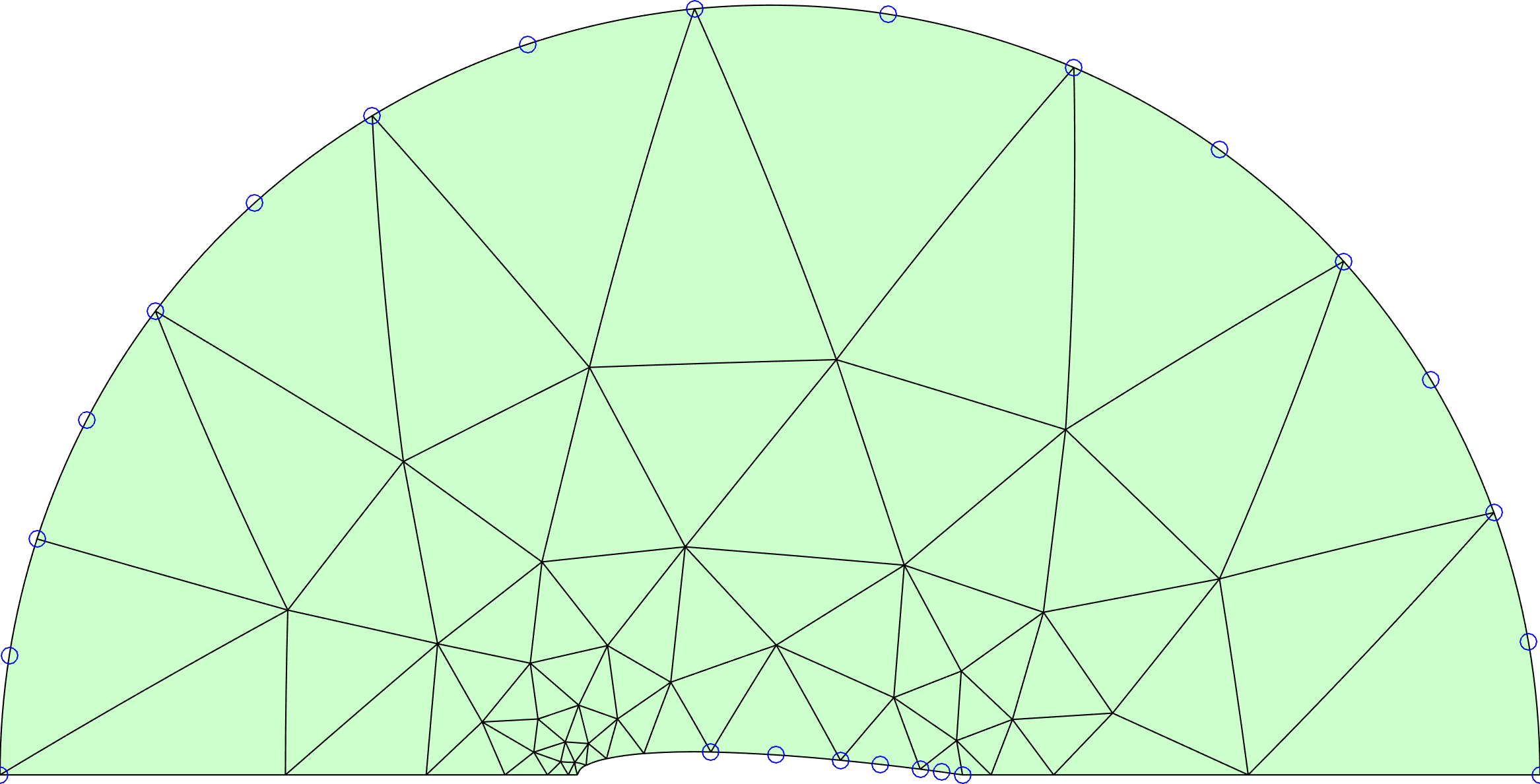} \quad
\includegraphics[width=0.48\textwidth]{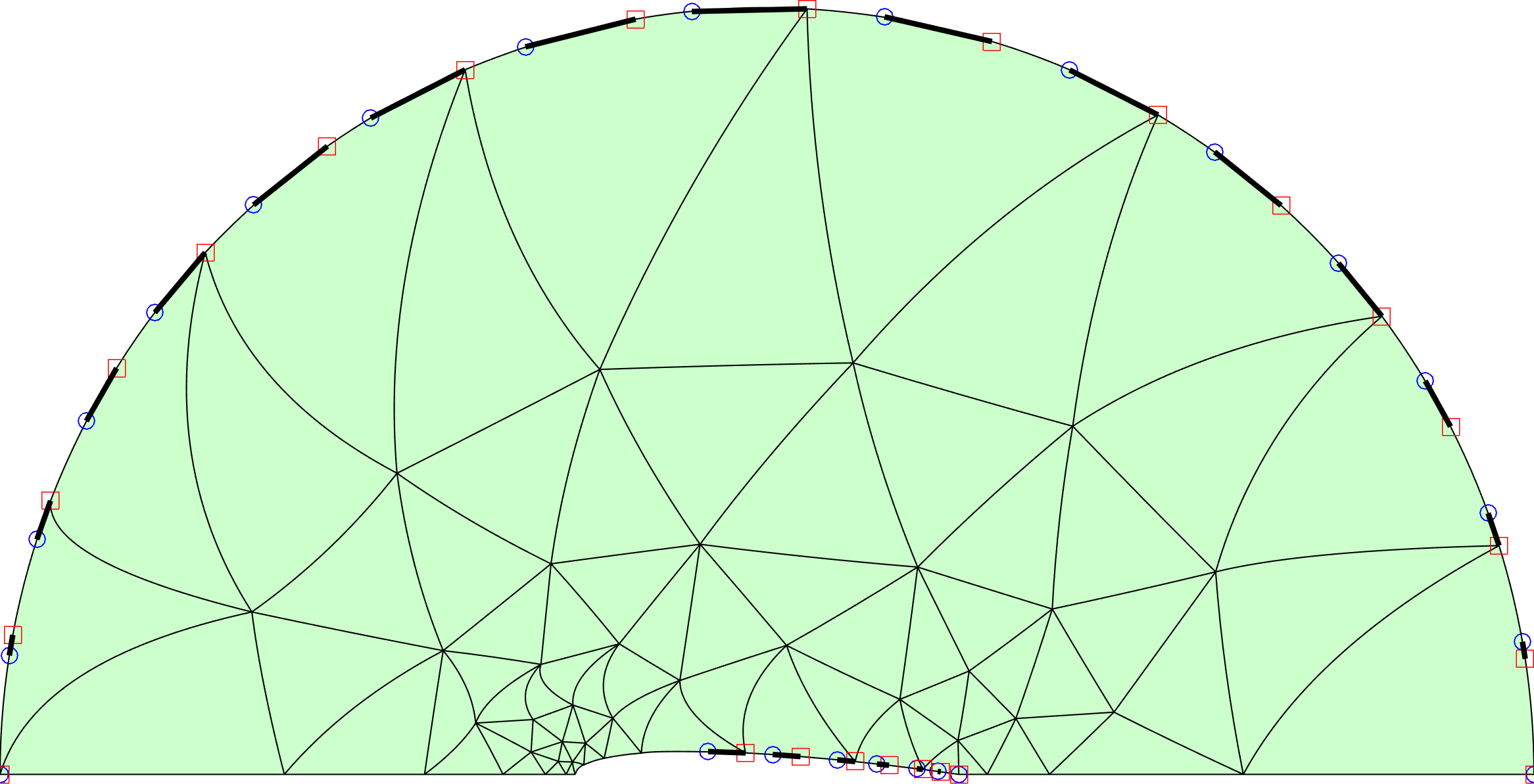}
\caption{Example of admissible mapping
         $\dommap\in\dommapsp_{\mshsz,\mshord}^\text{b}$ ($q=2$)
         in the case where $\rdom=\pdom$ (airfoil with circular farfield)
         that satisfies $\pdom=\dommap(\rdom)$. The reference domain $\rdom$
         (\textit{left}) can map to any configuration $\dommap(\rdom)$
         (\textit{right}) that corresponds to nodes sliding along the
         original boundary. For clarity, a subset of the reference nodes
         $\rnode_I$ (blue circles) and the corresponding mapped nodes
         $\pnode_I=\dommap(\rnode_I)$ (red squares) are shown and
         connected by a thick black line to highlight conformity to
         the boundary. Because high-order elements ($q=2$) are used,
         edge midnodes are included.}
\label{fig:bndslide}
\end{figure}

The first condition implies the number of boundaries on which $\rnode_I$
lies ($|\Ical_I|$) defines the number of constraints $\pnode_I$ must
satisfy, which will be embedded into the function space.
If $|\Ical_I|\geq d$, the position $\pnode_I$ is prescribed
regardless of the domain mapping (assuming all constraints independent);
otherwise, the nodes $\pnode_I$ can slide along the appropriate boundaries
\cite{huang2021robust}.
For each node $\pnode_I$ ($I=1,\dots,N_\mathrm{v}$), we identify a subset of
its entries $\Jcal_I\subset\{1,\dots,d\}$ of size
$|\Jcal_I|=N_{\text{u},I}\coloneqq d-|\Ical_I|$
as unconstrained degrees of freedom $\hat{y}_I$, i.e.,
\begin{equation}
 (\hat{y}_I)_i = (\pnode_I)_{\Jcal_I(i)}, \qquad i=1,\dots,d.
\end{equation}
In the case where a node does not lie on a boundary (interior node),
there are no constraints at that point so $\Jcal_I=\{1,\dots,d\}$ and
$\hat{y}_I = \pnode_I$ (all degrees of freedom at the node are unconstrained).
Then, the nodal position $\pnode_I$ that satisfies
the appropriate constraints is given by
\begin{equation}
 \pnode_I=\chi_I(\hat{y}_I;\rnode_I),
\end{equation}
where the mapping
$\func{\chi_I}{\Rbb^{N_{\text{u},I}}\times\Rbb^d}{\Rbb^d}$ 
is determined from either an analytical expression
\cite{zahr_implicit_2020,zahr2020r} for the surfaces $\partial\pdom_i$
for $i\in\Ical_I$ or from the CAD representation; for interior nodes,
$\chi_I$ is simply the identity map, i.e.,
$\chi_I(\hat{y}_I;\rnode_I) = \hat{y}_I$. For the special case
of straight-sided boundaries, $\chi_I$ can be determined directly
from the normal vector defining the planes $\partial\pdom_i$
\cite{huang2021robust}.
By collecting the functions $\{\chi_I\}_{I=1}^{N_\mathrm{v}}$ into
a single mapping $\func{\chibold}{\Rbb^{N_\text{u}}}{\Rbb^{dN_\mathrm{v}}}$,
we obtain a relationship between a vector of unconstrained
degrees of freedom $\ybm\in\Rbb^{N_\text{u}}$ and the physical
nodal coordinates $\xbm=\chibold(\ybm)$ that lie on the appropriate
boundaries or interior of $\pdom$, where
$N_\text{u}=\sum_{I=1}^{N_\mathrm{v}} N_{\text{u},I}$.
Finally, we define $\dommapsp_{\mshsz,\mshord}^\text{b}$
to be the following subset of $\dommapsp_{\mshsz,\mshord}$
\begin{equation}
 \dommapsp_{\mshsz,\mshord}^\text{b} \coloneqq
 \left\{
 \dommap_{\mshsz,\mshord}(\,\cdot\,;\chibold(\ybm))
 \suchthat \ybm\in\Rbb^{N_\text{u}}
 \right\}
 \subset \dommapsp_{\mshsz,\mshord}.
\end{equation}
i.e., the set of mappings of the form in $\dommapsp_{\mshsz,\mshord}$ 
that guarantees the image of the reference nodes lie on the appropriate
boundaries or interior of $\pdom$, which in turn ensures $\rdom$ maps to
$\pdom$ to $\Ocal(\mshsz^{\mshord+1})$ (provided the mapping is invertible).
Then, for any element $\dommap\in\dommapsp_{\mshsz,\mshord}^\text{b}$,
there exists $\ybm\in\Rbb^{N_\text{u}}$ such that
$\dommap=\dommap_{\mshsz,\mshord}^\text{b}(\,\cdot\,;\ybm)$,
where
\begin{equation}
 \ybm \mapsto \dommap_{\mshsz,\mshord}^\text{b}(\,\cdot\,;\ybm)
 \coloneqq \dommap_{\mshsz,\mshord}(\,\cdot\,;\chibold(\ybm)).
\end{equation}
Figure~\ref{fig:bndslide} illustrates a sample mapping from
$\dommapsp_{\mshsz,\mshord}^\text{b}$ for a non-trivial geometry.
The reader is referred to \cite{huang2021robust} for a detailed
description of this procedure including a fully automated algorithm
for constructing $\chibold$ in the case where all boundaries are
planar; for curved boundaries, see \cite{zahr2020r}.

\begin{remark}
The mesh $\Ecal_{h,q}$ is not required to be the same mesh associated
with the discretization in Section~\ref{sec:govern:disc}; however, we
make this choice in our numerical experiments for convenience.
\end{remark}

\begin{remark}
Identification of the unconstrained degrees of freedom $\Jcal_\Ical$
cannot be done arbitrarily. For example, in $d=2$ dimensions, for a
node on a boundary aligned with the $x_1$-axis, the first component
must be the unconstrained degree of freedom because perturbation in
the $x_2$-direction would cause the node to come off the boundary.
In general, the collection of surface normals at any point
$x\in\bigcap_{i\in\Ical_I} \partial\pdom_i$ can be used to
determine admissible partitions; see \cite{huang2021robust}.
\end{remark}

\begin{remark}
A simpler way to enforce that $\rnode_I$ lies on the appropriate boundary is
to fix $\pnode_I = \bar\dommap(\rnode_I)$ for any $\rnode_I\in\partial\rdom$,
i.e., fix all boundary nodes based on the nominal mapping
(Figure~\ref{fig:refmap_ex} for an example where $\bar\dommap=\mathrm{Id}$).
This is feasible for some problems; however, for problems requiring
large domain deformations to track features or those where features
intersect boundaries, e.g., transonic flow, nodes must be allowed
to slide along boundaries.
\end{remark}

\begin{remark}
While this approach is highly flexible, parametrizing the domain mapping
with all unconstrained degrees of freedom of a high-order mesh can be
burdensome. This can be avoided by using, for example, radial basis
functions to control the unconstrained degrees of freedom, which
limit the space of domain mappings, but still inherit the
boundary-preservation property.
\end{remark}

\section{Model reduction via optimization-based implicit feature tracking}
\label{sec:rom}
In this section we further motivate the limitations of the traditional
subspace approximation associated with projection-based model reduction
for convection-dominated problems (Section~\ref{sec:rom:sub}) and introduce
the proposed approach for defining a nonlinear approximation: deforming
modes (Section~\ref{sec:rom:nlman}). The nonlinear approximation is
embedded in a residual minimization framework that simultaneously
minimizes the HDM residual over the reduced coefficients of the
(deforming) modes and the deformation parameters (Section~\ref{sec:rom:ift}),
an approach to model reduction we called \textit{implicit feature tracking}
because, in the physical domain, local features in the solution tend to be
tracked (or aligned) with corresponding features in the basis without
explicitly accounting for them. Finally, the implicit feature tracking
approach is broken into an offline phase (Section~\ref{sec:rom:off})
where a collection of ``aligned'' snapshots are generated for optimal
compression into a basis, and an online phase (Section~\ref{sec:rom:on})
where the implicit feature tracking optimization formulation is leveraged
to provide accurate, low-dimensional approximations to the PDE solutions.

\subsection{Model reduction via subspace approximation}
\label{sec:rom:sub}
Central to most projection-based model reduction methods
\cite{benner2015survey}
is the following affine subspace approximation of the HDM solution
for any $\dommap\in\dommapsp$ and $\param\in\paramsp$
\begin{equation} \label{eqn:approx-ref}
 \rstvc_h(\,\cdot\,; \dommap, \param) \approx
 \rstvc_k(\,\cdot\,; \dommap, \param) \in \Vcal_k \coloneqq
 \left\{\bar{\rstvc}_h+\sum_{i=1}^k \alpha_i\Phi_{h,i} \suchthat \alpha_1,\dots,\alpha_k\in\Rbb\right\} \subseteq \Vcal_h
\end{equation}
where $\Vcal_k$ is a $k$-dimensional affine subspace of the finite-dimensional
trial space $\Vcal_h$ (function space over reference domain),
$\bar\rstvc_h\in\Vcal_h$ is the affine offset, and
$\Phi_{h,i}\in\Vcal_h$ are the reduced basis functions. When composed with
the domain mapping, this leads to the following approximation in the
physical domain
\begin{equation} \label{eqn:approx-phys}
 \pstvc(\,\cdot\,; \param) \approx
 \pstvc_h(\,\cdot\,; \dommap, \param) \approx
 \rstvc_k(\,\cdot\,; \dommap, \param)\circ\dommap^{-1} \in
 \Wcal_{k,\dommap} \coloneqq
 \Vcal_k\circ\dommap^{-1} \subseteq
 \Vcal_h\circ\dommap^{-1}.
\end{equation}
For a fixed $\dommap\in\dommapsp$, an element $w\in\Wcal_{k,\dommap}$
takes the form
\begin{equation} \label{eqn:approx-phys2}
 w = W \circ \dommap^{-1} = 
 \bar\rstvc_h\circ\dommap^{-1} +
 \sum_{i=1}^k \alpha_i \Phi_{h,i} \circ \dommap^{-1}.
\end{equation}
for some $W=\bar\rstvc_h+\sum_{i=1}^k \alpha_i \Phi_{h,i}\in\Vcal_k$.
Because $\dommap$ is fixed, $\Wcal_{k,\dommap}$ is an affine
subspace of $\Vcal_h\circ\dommap^{-1}$ with offset
$\bar\rstvc_h\circ\dommap^{-1}$ and basis functions
$\Phi_{h,i} \circ \dommap^{-1}$.

The algebraic representation of (\ref{eqn:approx-ref}) reads
\begin{equation}
 \dstvc(\dommap,\param) \approx \dstvc_k(\dommap,\param) \in \Vboldcal_k
 \coloneqq \left\{\bar\dstvc+\Phibold_k\vbm \suchthat \vbm\in\Rbb^k\right\},
\end{equation}
where $\bar\dstvc = \Upsilon(\bar{\rstvc}_h) \in \Rbb^k$ and
$\Phibold_k=[\Upsilon(\Phi_{h,1}),\dots,\Upsilon(\Phi_{h,k})]\in\Rbb^{N\times k}$
are the algebraic representation of the affine offset and reduced
basis functions, and
$\dstvc_k(\dommap,\param) = \Upsilon(\rstvc_k(\,\cdot\,; \dommap, \param))$
is the algebraic representation of the affine approximation $\rstvc_k$.
The element $\dstvc_k(\dommap, \param)\in\Vboldcal_k$
(or equivalently $\rstvc_k(\,\cdot\,; \dommap, \param)\in\Vcal_k$)
used to approximate the HDM solution is defined as the solution of
either a Galerkin or Petrov-Galerkin projection of the governing equations,
i.e., $\dstvc_k(\dommap,\param)=\bar\dstvc+\Phibold_k\vbm_k(\dommap,\param)$,
where $\vbm_k(\dommap,\param)$ satisfies either
\begin{equation} \label{eqn:gal}
 \Phibold_k^T\dres(\bar\dstvc+\Phibold_k\vbm_k(\dommap,\param);\dommap,\param) = \zerobold
\end{equation}
for a Galerkin ROM or
\begin{equation} \label{eqn:minres}
 \vbm_k(\dommap,\param) =
 \argmin_{\vbm\in\Rbb^k} \frac{1}{2}\norm{\dres(\bar\dstvc+\Phibold_k\vbm;\dommap,\param)}_2^2
\end{equation}
for a minimum-residual ROM.
Standard approaches to model reduction are formulated directly on
the physical domain, which is easily obtained by taking the
reference and physical domains to be equivalent ($\rdom=\pdom$)
and the domain mapping to be the nominal map $\dommap=\bar\dommap=\mathrm{Id}$.
In this setting, the reference and physical domain formulations
are equivalent.

\subsection{Nonlinear approximation via deforming modes}
\label{sec:rom:nlman}
\begin{figure}
\centering
\input{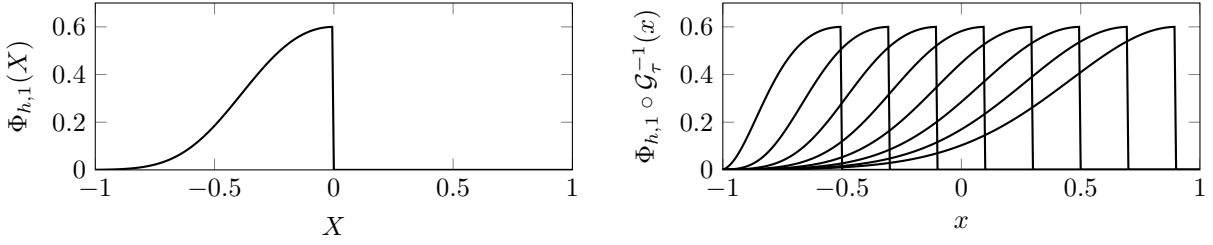}
\caption{A fixed basis function $\Phi_{h,1}$ in the reference domain
         based on a cutoff Gaussian (\textit{left}) and several instances
         instances of the deformed mode $\Phi_{h,1}\circ\Gcal_\tau^{-1}$
         in the physical domain (\textit{right}), where $\Gcal_\tau$
         is defined in (\ref{eqn:domdef_quad}) and $\tau\in[-0.5, 0.9]$.
         A linear approximation associated with the fixed mode $\Phi_{h,1}$
         would be limited to scaling this mode and unable to capture
         convection behavior, whereas the deforming mode
         $\Phi_{h,1}\circ\Gcal_\tau^{-1}$ can capture convection provided
         the mapping is chosen appropriately.}
\label{fig:demo1d0_nlman}
\end{figure}

As demonstrated in Section~\ref{sec:intro}
(Figure~\ref{fig:demo1d0_linapprox}), a linear subspace is not an
ideal low-dimensional approximation for convection-dominated problems
due to the slowly decaying Kolmogorov $n$-width.
Instead of approximating the PDE solution $\pstvc(\,\cdot\,;\param)$
with the HDM solution $\pstvc_h(\,\cdot\,;\dommap,\param)\in\Wcal_h$
for some $\dommap\in\dommapsp$, which is further approximated in the
affine subspace $\Wcal_{k,\dommap}$, we propose to directly approximate
$\pstvc(\,\cdot\,;\param)$ in the nonlinear manifold defined by
the union of all affine subspaces $\Wcal_{k,\dommap}$
\begin{equation} \label{eqn:nlman-pde2rom}
 \pstvc(\,\cdot\,; \param) \approx
 \pstvc_k(\,\cdot\,; \param) \in
 \Wcal_k \coloneqq \bigcup_{\dommap\in\dommapsp} \Wcal_{k,\dommap}.
\end{equation}
For now, we consider an arbitrary collection of domain mappings
$\dommapsp$; in Sections~\ref{sec:rom:off} and~\ref{sec:rom:on},
we will specialize it for the offline and online phase, respectively,
of the proposed method.
The nonlinear manifold $\Wcal_k$ and concept of deforming modes
is illustrated in Figure~\ref{fig:demo1d0_nlman} using the cutoff
Gaussian function.
Any element $w\in\Wcal_k$ takes the form
\begin{equation}
 w = W \circ \dommap^{-1} = 
 \bar\rstvc_h\circ\dommap^{-1} +
 \sum_{i=1}^k \alpha_i \Phi_{h,i} \circ \dommap^{-1}.
\end{equation}
for some $\dommap\in\dommapsp$ and
$W=\bar\rstvc_h+\sum_{i=1}^k \alpha_i\Phi_{h,i}\in\Vcal_k$.
Thus, the manifold approximation $\pstvc_k(\,\cdot\,;\param)$ is
uniquely defined from the domain mapping
$\dommap_k(\,\cdot\,;\param)\in\dommapsp$,
and the reference domain approximation
$\rstvc_k(\,\cdot\,;\param)\in\Vcal_k$ as
\begin{equation}
 \pstvc_k(\,\cdot\,;\param) =
 \rstvc_k(\,\cdot\,;\param) \circ \dommap_k(\,\cdot\,;\param)^{-1}.
\end{equation}
That is, a complete reduced-order model must define both
$\dommap_k$ and $\rstvc_k$, which is in contrast to the setting
in Section~\ref{sec:rom:sub} where $\rstvc_k\in\Vcal_k$ was determined
for a fixed value of the domain mapping.

\begin{remark}
For any $\dommap\in\dommapsp$, $\Wcal_{k,\dommap}$
is an affine subspace because the elements of $\Wcal_{k,\dommap}$ are
a linear superposition of the modes $\Phi_{h,i}\circ\dommap^{-1}$
and the offset $\bar\rstvc\circ\dommap^{-1}$. On the other
hand, $\Wcal_k$ is a nonlinear manifold because the modes
$\Phi_{h,i}\circ\dommap^{-1}$ deform as the domain mapping
$\dommap\in\dommapsp$ varies.
\end{remark}

\begin{remark}
Recall the HDM approximates the PDE solution
$\pstvc(\,\cdot\,;\param)\approx \pstvc_h(\,\cdot\,;\dommap,\param)=
 \rstvc_h(\,\cdot\,;\dommap,\param)\circ\dommap^{-1}$
for any $\dommap\in\dommapsp$ and the PDE solution is approximated on
the nonlinear manifold as
$\pstvc(\,\cdot\,;\param)\approx\pstvc_k(\,\cdot\,;\param)=
 \rstvc_k(\,\cdot\,;\param)\circ\dommap_k(\,\cdot\,;\param)^{-1}$.
Therefore, the manifold approximation of $\pstvc_k(\,\cdot\,;\param)$
can be viewed as an approximation to
$\pstvc_h(\,\cdot\,;\dommap_k(\,\cdot\,;\param),\param)$ in
either the physical or reference domain, i.e.,
\begin{equation} \label{eqn:nlman-hdm2rom}
 \pstvc_h(\,\cdot\,;\dommap_k(\,\cdot\,;\param),\param) \approx
 \pstvc_k(\,\cdot\,;\param)\in\Wcal_k, \qquad
 \rstvc_h(\,\cdot\,;\dommap_k(\,\cdot\,;\param),\param) \approx
 \rstvc_k(\,\cdot\,;\param)\in\Vcal_k.
\end{equation}
 Here, we assume that for reasonable $\dommap\in\dommapsp$ (well-conditioned),
 $\pstvc_h(\,\cdot\,;\dommap,\param)$ provides a
 suitable approximation to $\pstvc(\,\cdot\,;\param)$, i.e.,
 the computational grid is sufficiently refined or $r$-refinement
 is used. With this assumption, we define
 $\dommap_k(\,\cdot\,;\param)$ to ensure accuracy of the approximation
 in (\ref{eqn:nlman-hdm2rom}), from which accuracy of the
 approximation in (\ref{eqn:nlman-pde2rom}) will follow
 (Section~\ref{sec:rom:ift}).
\end{remark}

\begin{figure}
 \centering
 \includegraphics[width=0.32\textwidth]{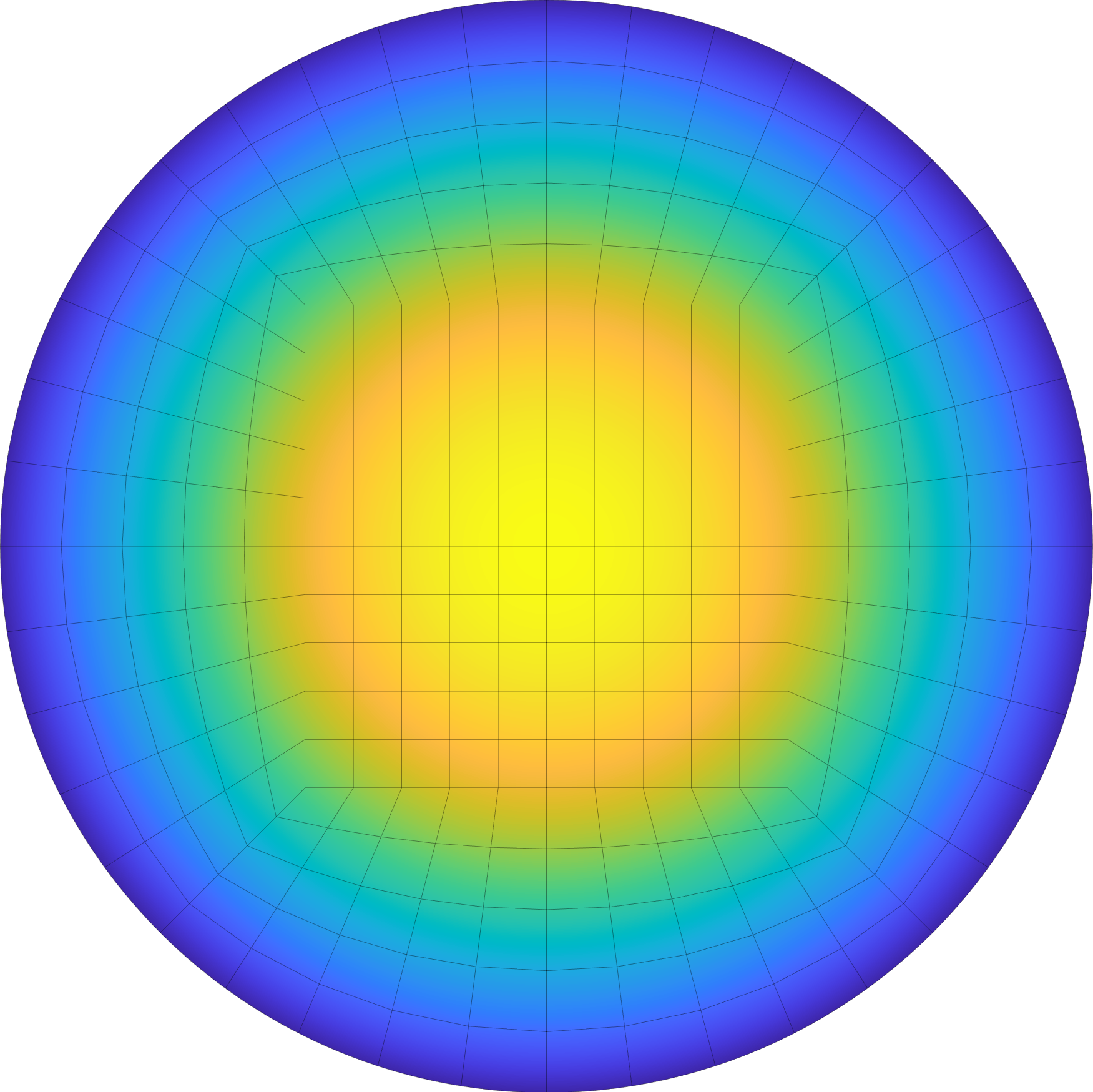} \qquad
 \includegraphics[width=0.32\textwidth]{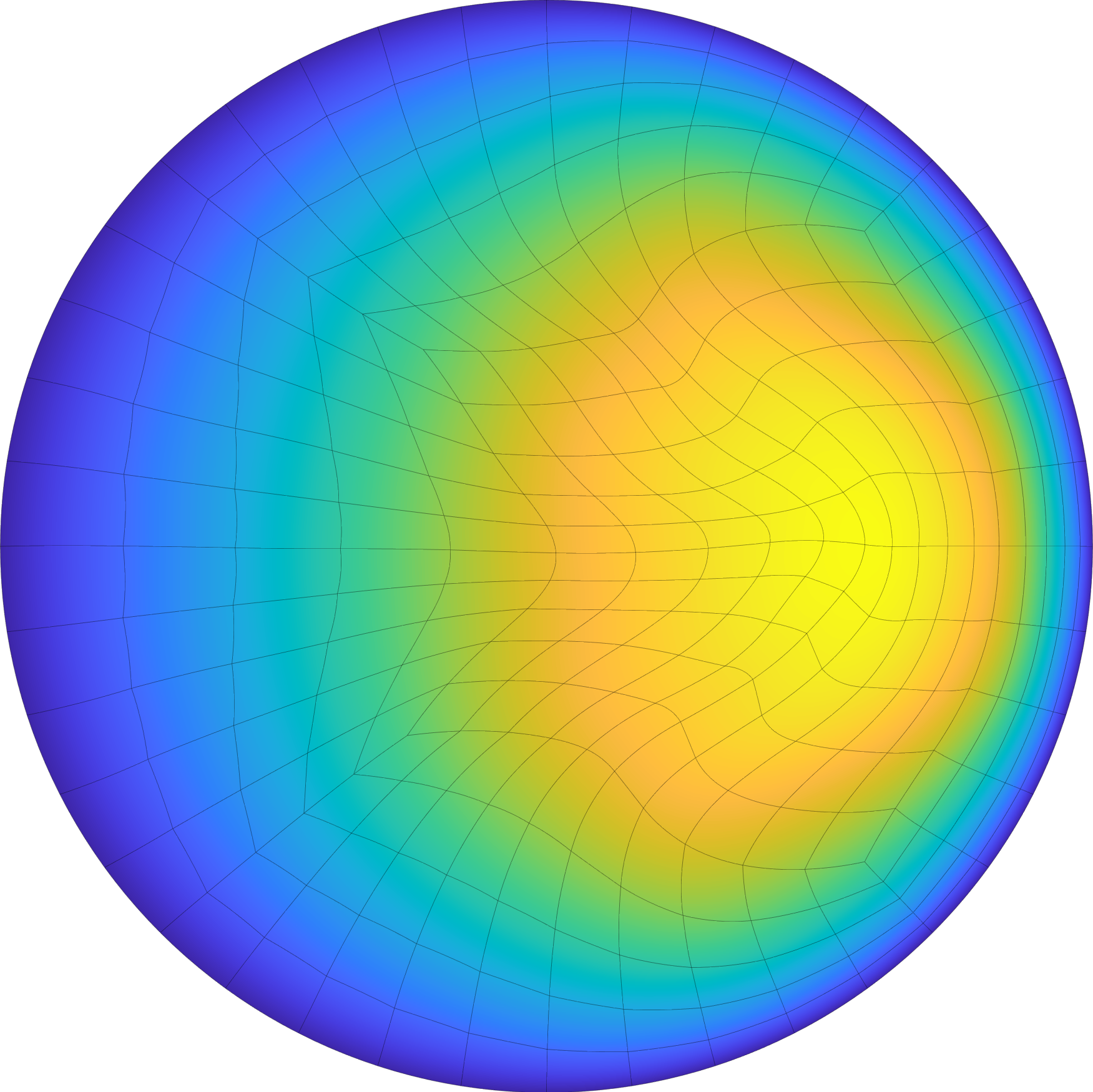}
\caption{A global basis function in reference domain $\Phi_{h,1}$
         (\textit{left}) and after composition with the domain
         mapping $\dommap\in\dommapsp$ in Figure~\ref{fig:refmap_ex}
         (\textit{right}).}
\label{fig:dual_interp}
\end{figure}

\begin{remark} \label{rem:dual}
There are two ways to interpret the approximation proposed in
(\ref{eqn:nlman-pde2rom}). In the physical domain, this reads
\begin{equation}
 \pstvc(\,\cdot\,;\param) \approx
 \bar\rstvc_h\circ\dommap^{-1}+
 \sum_{i=1}^k \alpha_i\Phi_{h,i}\circ\dommap^{-1}
\end{equation}
for $\dommap\in\dommapsp$, which shows the mapping-independent solution
$\pstvc(\,\cdot\,;\param)$ is being approximated by the deformable
modes $\Phi_{h,i}\circ\dommap^{-1}$. On the other hand,
in the reference domain, we have
\begin{equation}
 \rstvc(\,\cdot\,;\dommap,\param) \approx
 \bar\rstvc_h + \sum_{i=1}^k \alpha_i\Phi_{h,i},
\end{equation}
which shows the mapping-dependent solution $\rstvc(\,\cdot\,;\dommap,\param)$
is being approximated by the fixed modes $\Phi_{h,i}$. That is, in the
physical domain, the modes deform to the fixed solution, while in the
reference domain, the solution deforms to the fixed modes. This dual
interpretation is illustrated in Figures~\ref{fig:refmap_ex} and
\ref{fig:dual_interp}: the solution in the reference domain
(Figure~\ref{fig:refmap_ex}, \textit{left}) deforms to the fixed
mode $\Phi_{h,1}$  (Figure~\ref{fig:dual_interp}, \textit{left}), or
equivalently, the mapped mode $\Phi_{h,1}\circ\dommap$
(Figure~\ref{fig:dual_interp}, \textit{right}) deforms to
solution in the physical domain (Figure~\ref{fig:refmap_ex}, \textit{right}).
\end{remark}

\subsection{Optimization formulation for model reduction by
            implicit feature tracking}
\label{sec:rom:ift}
For convection-dominated problems, local features, e.g., shocks,
rarefactions, and vortical structures, will move in $\pdom$ as
$\param\in\paramsp$ varies, which makes fixed modes $\{\varphi_i\}_{i=1}^k$
constructed from \textit{a priori} training ill-suited to approximate
the solution $\pstvc(\,\cdot\,;\param)$ due to the slowly decaying
Kolmogorov $n$-width (or will at least require extensive training
and relatively large $k$). However, by writing the physical modes
as in (\ref{eqn:approx-phys2}), i.e., $\varphi_i=\Phi_{h,i}\circ\dommap^{-1}$,
the domain mapping $\dommap$ provides the opportunity
to deform the reference modes $\Phi_{h,i}$ in such a way that features
in $\Phi_{h,i}$ align with features in $\pstvc(\,\cdot\,;\param)$. From
the dual interpretation in Remark~\ref{rem:dual}, this can equivalently
be stated as deforming the solution in the reference domain
$\rstvc(\,\cdot\,;\dommap,\param)$ such that its local features align with
the fixed features in $\Phi_i$, effectively
removing the convection-dominated nature of the solution. 

To accomplish this feature alignment, we define a reduced-order
model as an unconstrained optimization problem over the nonlinear
manifold $\Wcal_k$
\begin{equation}\label{eqn:romift1}
\pstvc_k(\,\cdot\,;\param) \coloneqq \argmin_{w\in\Wcal_k}\hat\Jcal(w;\param).
\end{equation}
The objective function $\hat\Jcal$ is constructed to (implicitly)
promote feature alignment by penalizing an error indicator
associated with the approximation and ensure the mapping $\dommap$
is a bijection between $\rdom$ and $\pdom$.
Optimization over $\Wcal_k$ amounts to simultaneously determining the
element of the affine space $\Vcal_k$ and
the domain mapping space ($\dommapsp$), which can be explicitly written in terms
of reference domain quantities as
\begin{equation}\label{eqn:romift2}
 (\rstvc_k(\,\cdot\,;\param),\dommap_k(\,\cdot\,;\param)) \coloneqq \argmin_{(W,\dommap)\in\Vcal_k\times\dommapsp}\Jcal(W,\dommap;\param)
\end{equation}
where $\Jcal(W,\dommap;\param)\coloneqq\hat\Jcal(W\circ\dommap^{-1};\param)$.
Whether written in the physical domain or reference domain,
(\ref{eqn:romift1})-(\ref{eqn:romift2})
define the element of the nonlinear manifold used to approximate the PDE
solutions and will be referred to as the reduced-order model with implicit
feature tracking (ROM-IFT).
The algebraic representation of the ROM-IFT takes the form
\begin{equation}\label{eqn:romift3}
\rstvc_k(\,\cdot\,;\param) \coloneqq \Upsilon(\dstvc_k(\param)), \qquad
(\dstvc_k(\param),\dommap_k(\,\cdot\,;\param)) \coloneqq \argmin_{(\Vbm,\dommap)\in\Vboldcal_k\times\dommapsp} J(\Vbm,\dommap;\param)
\end{equation}
where the $\pstvc_k$ and $\rstvc_k$ are related by (\ref{eqn:approx-phys}).
Furthermore, we define $\wbm_k(\param)\in\Rbb^k$ such that
$\dstvc_k(\param) = \bar\dstvc + \Phibold_k\wbm_k(\param)$.

The objective function
$\func{J}{\Rbb^N\times\dommapsp\times\paramsp}{\Rbb}$ is composed of two
terms
\begin{equation}\label{eqn:romift4}
 J(\Vbm,\dommap; \param) \coloneqq J_\text{err}(\Vbm,\dommap; \param) + \kappa^2 J_\text{map}(\dommap),
\end{equation}
where $\func{J_\text{err}}{\Rbb^N\times\dommapsp\times\paramsp}{\Rbb}$
is an error indicator for the solution intended to promote feature
alignment, $\func{J_\text{map}}{\dommapsp}{\Rbb}$ is a function that
penalizes distortion of the domain mapping, and $\kappa\in\Rbb_{>0}$
is a weighting parameter. The error indicator is taken to be the norm
of the HDM residual
\begin{equation}\label{eqn:romift5}
J_\text{err}(\Vbm,\dommap;\param) \coloneqq \frac{1}{2}\norm{\dres(\Vbm;\dommap,\param)}_2^2,
\end{equation}
and mapping distortion term is taken as
\begin{equation}\label{eqn:romift6}
 J_\text{map}(\dommap) \coloneqq
 \frac{1}{2}\norm{\etabold(\dommap)-\etabold(\dommap_0)}_2^2,
\end{equation}
where $\func{\etabold}{\dommapsp}{\Rbb^{|\Ecal_h|}}$ is the
algebraic system corresponding to the elementwise distortion commonly used
for high-order mesh generation \cite{knupp2001algebraic}
\begin{equation}
\eta_K(\dommap) \coloneqq 
\int_K \left(\frac{\norm{\domjac_\dommap}_F^2}{\max\{\domdet_\dommap,\epsilon\}^{2/d}}\right)^2 \, dV,
\end{equation}
for any element $K$ in the computational mesh $\Ecal_{h,q}$ of $\rdom$,
$\epsilon\in\Rbb_{>0}$ is a tolerance preventing the denominator from
reaching $0$, and $\Gcal_0$ is a nominal mapping (usually
chosen to be $\bar\Gcal$, taken to be identity in this work).
The combined objective function can be written as
\begin{equation}
 J(\Vbm,\dommap;\param) = \frac{1}{2}\norm{\Fbm(\Vbm;\dommap,\param)}_2^2,
 \qquad
 (\Vbm,\dommap;\param)\mapsto\Fbm(\Vbm;\dommap,\param)\coloneqq
 \begin{bmatrix}
  \dres(\Vbm;\dommap,\param) \\
  \kappa \left(\etabold(\dommap)-\etabold(\dommap_0)\right)
 \end{bmatrix},
\end{equation}
which shows the ROM-IFT problem is a nonlinear least-squares problem;
this structure will be exploited to introduce efficient solvers in
Section~\ref{sec:solver}. In this work, we take $\epsilon=10^{-8}$
and choose the weighting parameter $\kappa$ using the algorithm
in \cite{huang2021robust}.

We conclude this section with some properties of the ROM-IFT method.
\begin{prop}
\label{prop:prop1}
For any $\param\in\paramsp$, $\hat\dstvc\in\Vboldcal_k$,
and $\hat\dommap\in\dommapsp$, we have
\begin{equation}
 J(\bar\dstvc+\Phibold_k\wbm_k(\param),\dommap_k(\,\cdot\,;\param);\param)\leq
 J(\hat\dstvc,\hat\dommap;\param).
\end{equation}
For $\kappa=0$, this reduces to
\begin{equation}
 \norm{\dres(\bar\dstvc+\Phibold_k\wbm_k(\param),\dommap_k(\,\cdot\,;\param);\param)}_2
 \leq
 \norm{\dres(\hat\dstvc, \hat\dommap; \param)}_2.
\end{equation}
\begin{proof}
Follows directly from the definition of $\wbm_k(\param)$ and
$\dommap_k(\,\cdot\,;\param)$ in (\ref{eqn:romift3}) as the arguments
that minimize $J$, which reduces to the norm of the residual function
$\dres(\,\cdot\,; \,\cdot\,, \param)$ for $\kappa=0$.
\end{proof}
\end{prop}

\begin{prop}
\label{prop:prop2}
For any $\param\in\paramsp$ and $\hat\dommap\in\dommapsp$, we have
\begin{equation}
 J(\bar\dstvc+\Phibold_k\wbm_k(\param),\dommap_k(\,\cdot\,;\param);\param) \leq
 J(\bar\dstvc+\Phibold_k\vbm_k(\hat\dommap,\param),\hat\dommap;\param),
\end{equation}
where $\vbm_k(\hat\dommap,\param)$ is the either the Galerkin
(\ref{eqn:gal}) or minimum-residual (\ref{eqn:minres}) reduced-order
model solution. For $\kappa=0$, this reduces to
\begin{equation}
 \norm{\dres(\bar\dstvc+\Phibold_k\wbm_k(\param), \dommap_k(\,\cdot\,;\param); \param)}_2
 \leq
 \norm{\dres(\bar\dstvc+\Phibold_k\vbm_k(\hat\dommap,\param), \hat\dommap; \param)}_2.
\end{equation}
\begin{proof}
Follows directly from Proposition~\ref{prop:prop1} with
$\hat\dstvc=\bar\dstvc+\Phibold_k\vbm_k(\hat\dommap,\param)$.
\end{proof}
\end{prop}

\begin{prop}
\label{prop:prop3}
Define $\Phibold_k\in\Rbb^{N\times k}$ and
$\Phibold_{k'}\in\Rbb^{N\times k'}$ such that
$\text{Ran}~\Phibold_{k'} \subseteq \text{Ran}~\Phibold_k$.
Then, for any $\param\in\paramsp$, we have
\begin{equation}
 J(\bar\dstvc+\Phibold_k\wbm_k(\param), \dommap_k(\,\cdot\,;\param); \param) \leq
 J(\bar\dstvc+\Phibold_{k'}\wbm_{k'}(\param), \dommap_{k'}(\,\cdot\,;\param); \param).
\end{equation}
For $\kappa=0$, this reduces to
\begin{equation}
 \norm{\dres(\bar\dstvc+\Phibold_k\wbm_k(\param), \dommap_k(\,\cdot\,;\param); \param)}_2
 \leq
 \norm{\dres(\bar\dstvc+\Phibold_{k'}\wbm_{k'}(\param), \dommap_{k'}(\,\cdot\,;\param); \param)}_2.
\end{equation}
\begin{proof}
Follows directly from Proposition~\ref{prop:prop1} with
$\hat\dstvc=\bar\dstvc+\Phibold_{k'}\wbm_{k'}(\param)\in\Vboldcal_{k'}\subseteq\Vboldcal_k$ and $\hat\Gcal=\Gcal_{k'}$.
\end{proof}
\end{prop}

\begin{prop}
\label{prop:prop4}
Let $\kappa=0$ and fix $\param\in\paramsp$.
Suppose there exists $\dommap\in\dommapsp$
such that $\dstvc(\dommap,\param)\in\Vboldcal_k$, where
$\dstvc(\dommap,\param)$ is the solution of
$\dres(\,\cdot\,;\dommap,\param)=\zerobold$. Then,
\begin{equation} \label{eqn:prop4a}
 \norm{\dres(\bar\dstvc+\Phibold_k\wbm_k(\param), \dommap_k(\,\cdot\,;\param); \param)} = \zerobold,
 \qquad
 \bar\dstvc+\Phibold_k\wbm_k(\param)=\dstvc(\dommap_k(\,\cdot\,;\param),\param).
\end{equation}
If the mapping $\dommap\in\dommapsp$ is the unique element such that
$\dstvc(\dommap,\param)\in\Vboldcal_k$, then
\begin{equation} \label{eqn:prop4b}
 \bar\dstvc+\Phibold_k\wbm_k(\param)=\dstvc(\dommap,\param).
\end{equation}
\begin{proof}
Equations (\ref{eqn:prop4a}) follow directly from Proposition~\ref{prop:prop1}
with $\hat\dstvc=\dstvc(\dommap,\param)$ and $\hat\dommap=\dommap$,
and from $\Ubm(\dommap_k(\,\cdot\,;\param),\param)$ as the
unique solution of $\dres(\,\cdot\,;\dommap_k(\,\cdot\,;\param),\param)=\zerobold$. The uniqueness of $\dommap$ implies $\dommap=\dommap_k(\,\cdot\,;\param)$
from which (\ref{eqn:prop4b}) follows.
\end{proof}
\end{prop}

\begin{remark}
For $\kappa = 0$, the optimization problem in
(\ref{eqn:romift3})-(\ref{eqn:romift4}) corresponds to residual
minimization over the nonlinear manifold $\Wcal_k$, which can be
directly contrasted to traditional minimum-residual approaches
\cite{legresley2006application,carlberg_efficient_2011} that freeze the domain
($\dommap=\dommap_0$) and optimize only over the affine subspace
$\Vcal_k$. Because the optimization problem is posed over a richer space,
we expect to obtain an approximation with a smaller residual than possible
with standard reduced-order models (Galerkin or minimum-residual) that use
frozen (non-deforming) modes (Proposition~\ref{prop:prop2}). As a result,
despite the focus of this work on convection-dominated PDEs, this approach
can be used for a larger class of parametrized PDEs.
\end{remark}

\begin{remark}
The proposed ROM-IFT is monotonic in the sense that enriching the
basis will only improve the solution in terms of the value of
$J$ (or the residual magnitude if $\kappa=0$)
(Proposition~\ref{prop:prop3}). Furthermore, if a unique domain
deformation exists that causes the HDM solution to lie in
the affine subspace $\Vboldcal_k$, the ROM-IFT method will
recover it ($\kappa=0$) (Proposition~\ref{prop:prop4}).
\end{remark}

\begin{remark}
The objective function of the optimization problem in
(\ref{eqn:romift3}) does not explicitly rely on identifying
or locating solution features, in contrast to those in
\cite{moretti1987computation,moretti2002thirty,taddei_reduced_2015},
which justifies the terminology ``implicit'' feature tracking.
As such, the
formulation applies generally to a wide class of problems
with various features (shocks, steep gradients, boundary
layers). 
\end{remark}


\subsection{Comparison:
            subspace residual minimization vs.
            implicit feature tracking}
\label{sec:rom:cmpr}
\begin{figure}[t]
 \centering
 \includegraphics[width=0.3\textwidth]{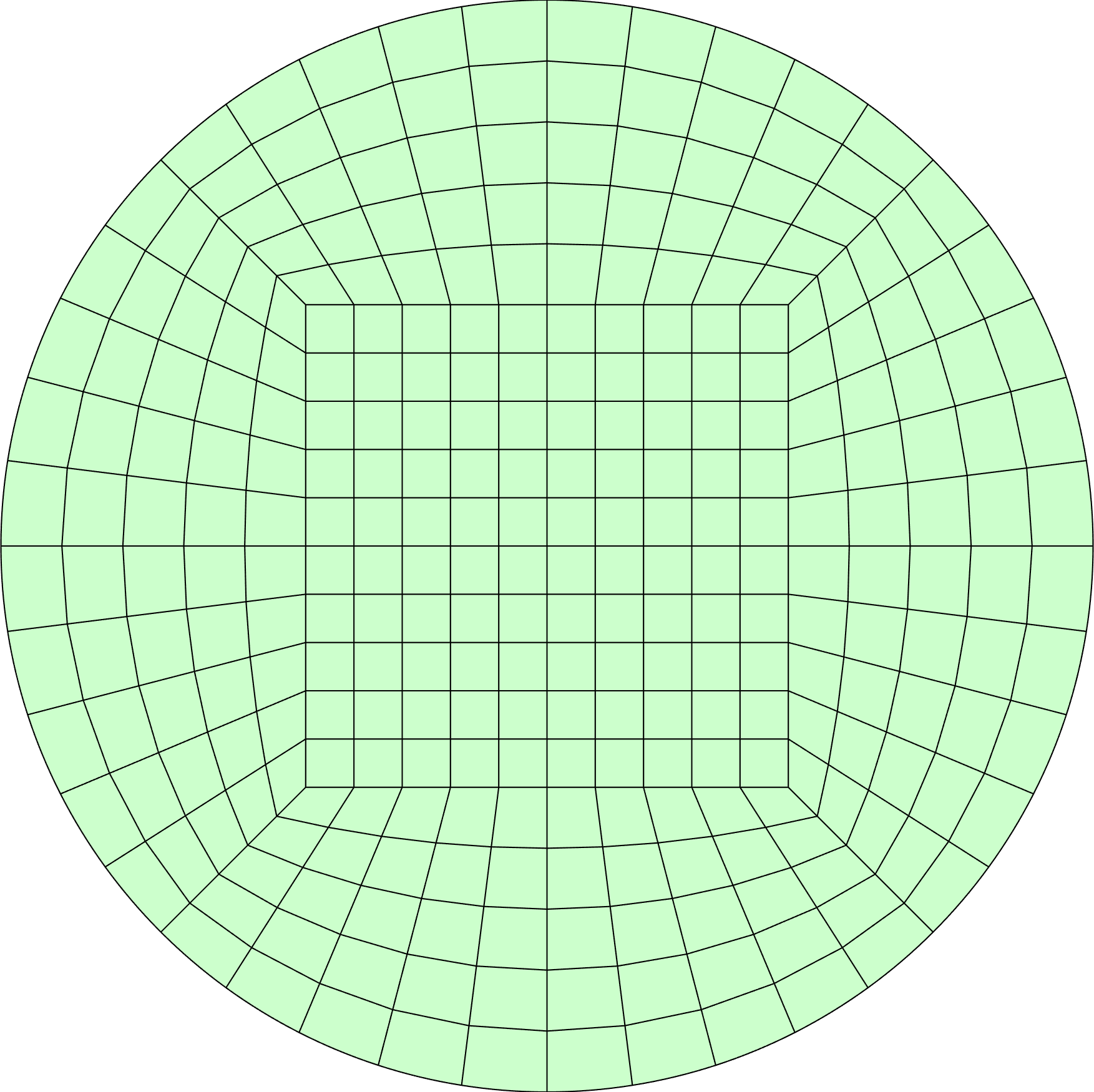} \quad
 \includegraphics[width=0.3\textwidth]{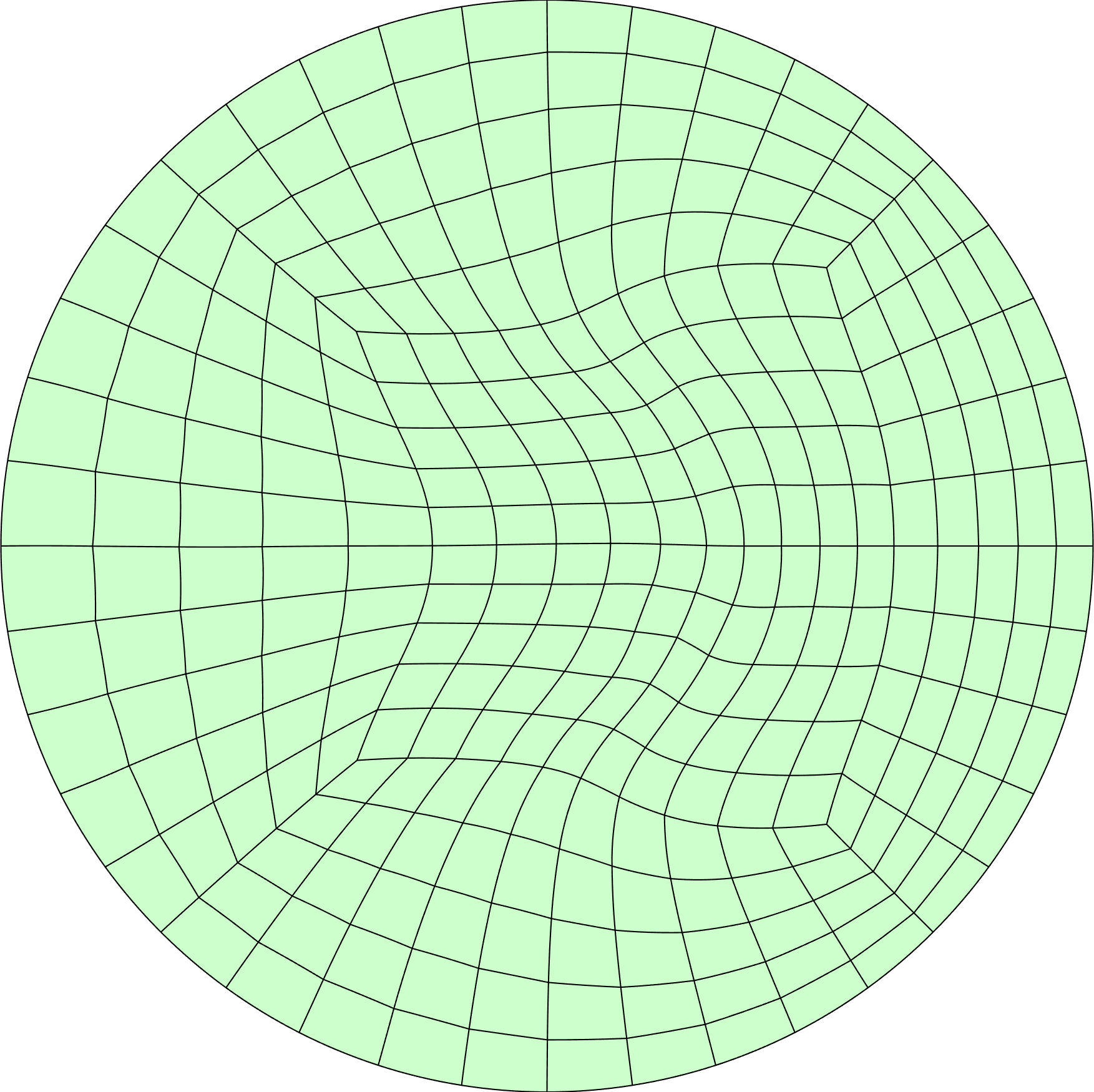} \quad
 \includegraphics[width=0.3\textwidth]{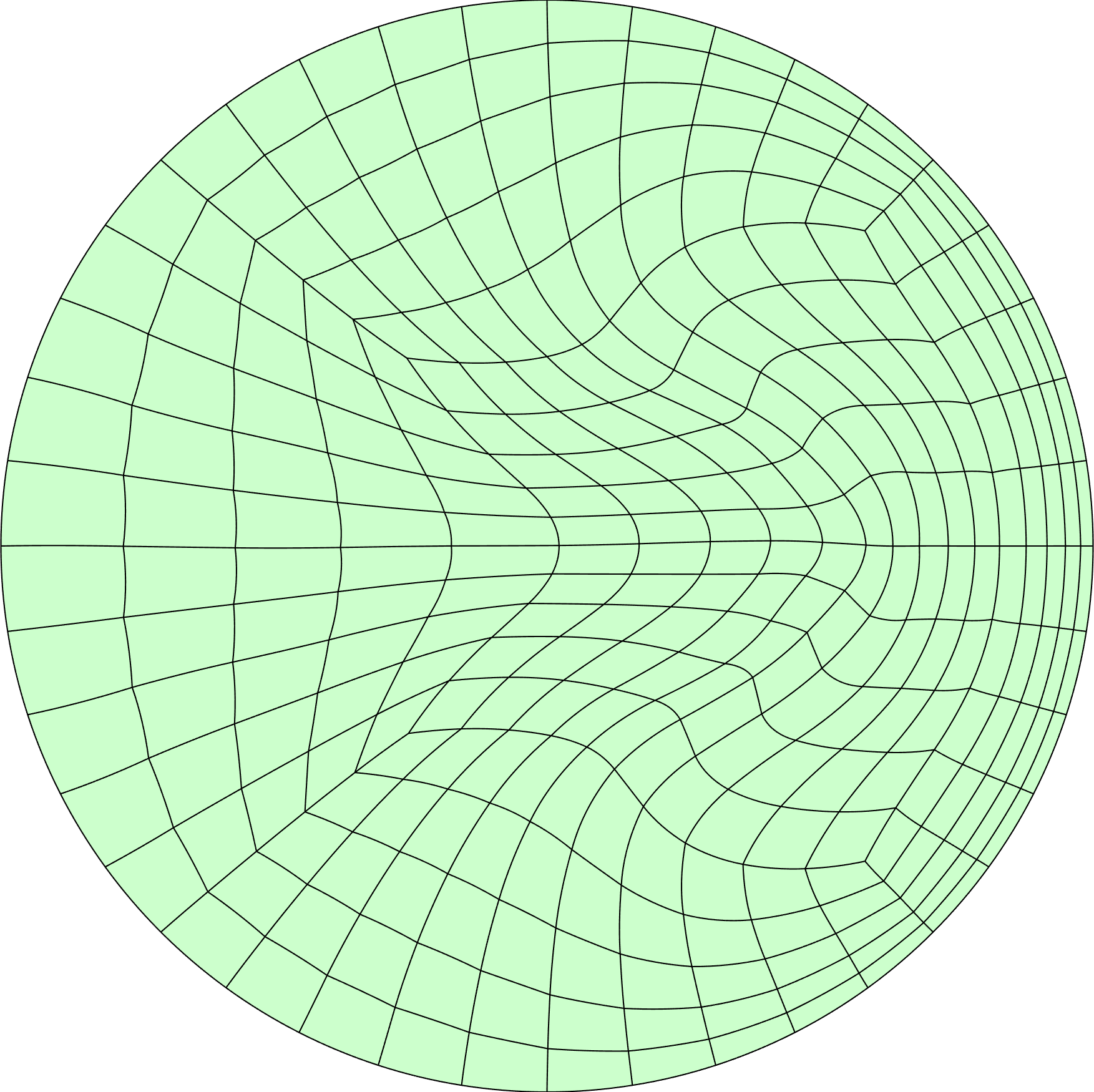}
\caption{Three instances of the domain mapping $\Tcal_s$ applied to a disk
 ($\rdom$): $\Tcal_0(\rdom)=\rdom$ (\textit{left}),
 $\Tcal_{0.5}(\rdom)$ (\textit{middle}), and
 $\Tcal_1(\rdom)=\Tcal(\rdom)$ (\textit{right}).}
\label{fig:demo1_msh}
\end{figure}

To demonstrate the ROM-IFT method (Section~\ref{sec:rom:ift})
and compare it to residual minimization over an affine space
(Section~\ref{sec:rom:sub}), we consider parametrized
advection-diffusion in the disk
$\pdom_\text{d}\coloneqq\{(x,y)\in\Rbb^2\mid x^2+y^2\leq 2\}$
from \cite{eftang_hp_2010}
\begin{equation} \label{eqn:advdiff-ronq}
 \nabla \cdot (\beta u) = \Delta u + 10\quad\text{ in }\Omega_\text{d},
 \qquad
 u = 0 \quad\text{ on }\partial\Omega_\text{d}
\end{equation}
where $(\phi,a)\ni\Dcal\mapsto\beta((\phi,a))\coloneqq a(\cos\phi,\sin\phi)$
is the advection velocity, $\Dcal = [0,\pi]\times [0, 10]$ is the parameter
domain, and $(x,\mu)\in\Omega_\text{d}\times\paramsp\mapsto u(x;\mu)\in\Rbb$
is the solution (Figure~\ref{fig:refmap_ex} for solution at $\mu=(0,10)$).
Furthermore, we define a reference domain $\Omega_0=\Omega_\text{d}$
and a family of domain mappings
$\func{\Tcal_s}{\Omega_\text{d}}{\Omega_\text{d}}$
for $s\in\Rbb$ that linearly interpolate between
the identity mapping and a prescribed mapping
$\func{\Tcal}{\Omega_\text{d}}{\Omega_\text{d}}$
(Figure~\ref{fig:demo1_msh}) defined as
\begin{equation} \label{eqn:disk-map}
 \rcoord\in\Omega_\text{d} \mapsto 
 \Tcal_s(\rcoord) \coloneqq \rcoord + s(\Tcal(\rcoord)-\rcoord).
\end{equation}
We transform the conservation law to the reference domain following
the procedure in Section~\ref{sec:govern:tclaw}; the solution of the
transformed conservation law is
$(\rcoord, \Tcal_s, \mu)\mapsto U(\rcoord;\Tcal_s,\mu)=u(\Tcal_s(\rcoord);\mu)$.
We discretize the transformed conservation law with the compact
discontinuous Galerkin (CDG) method \cite{peraire2008compact} to obtain
the algebraic representation $\Ubm(\Tcal_s, \mu)$.

\begin{figure}
 \centering
 \includegraphics[width=0.3\textwidth]{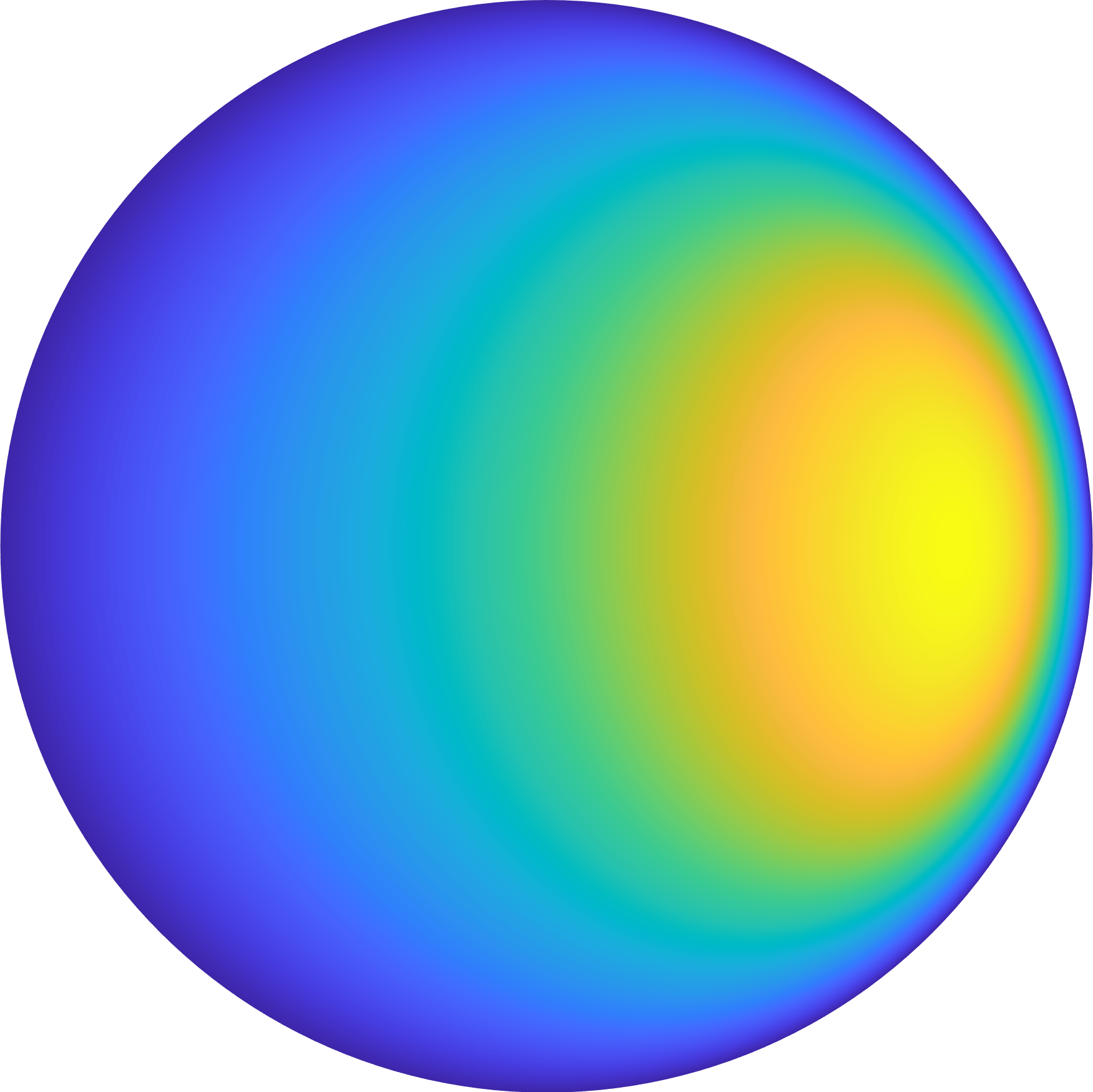} \quad
 \includegraphics[width=0.3\textwidth]{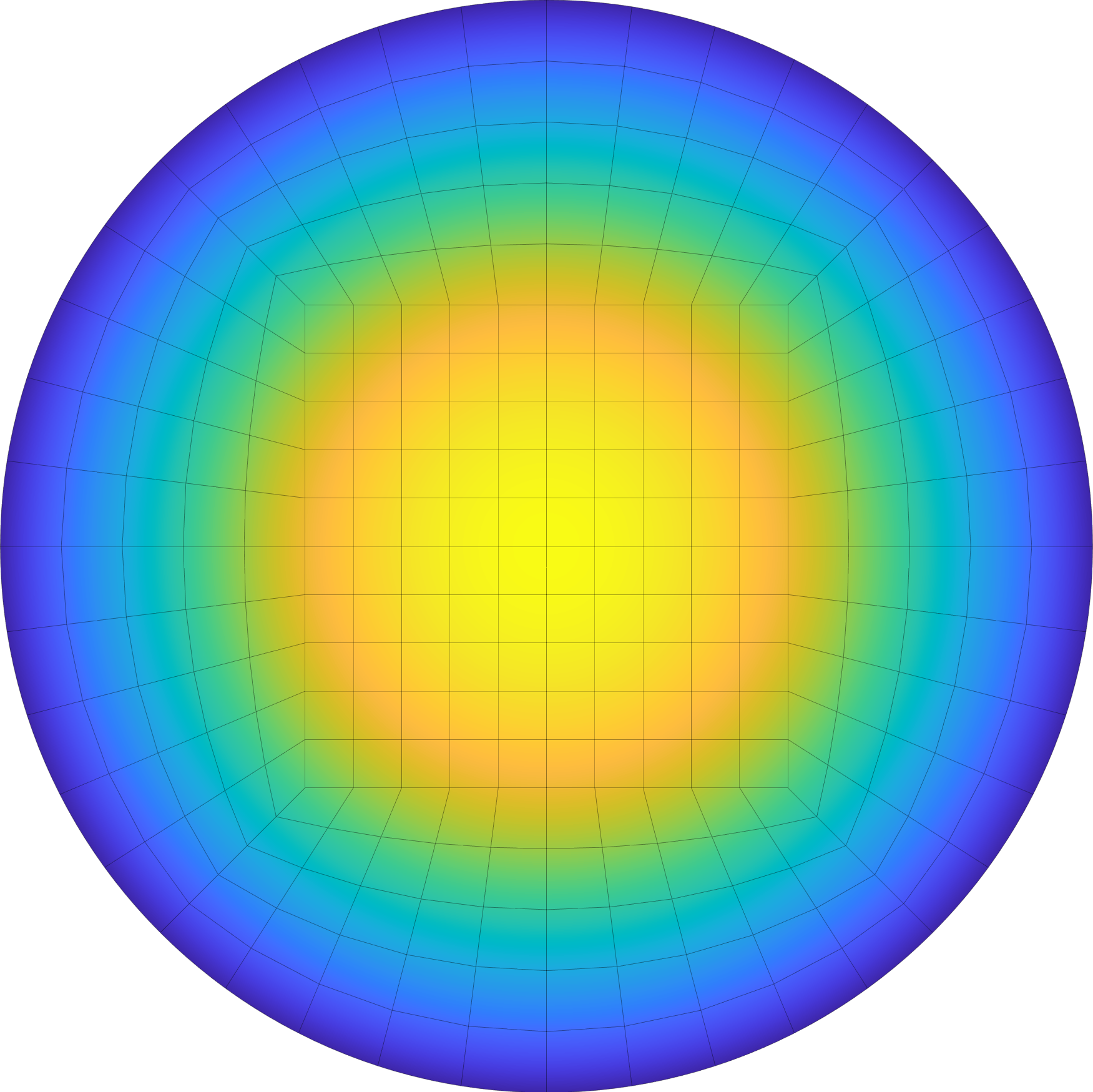} \quad
 \includegraphics[width=0.3\textwidth]{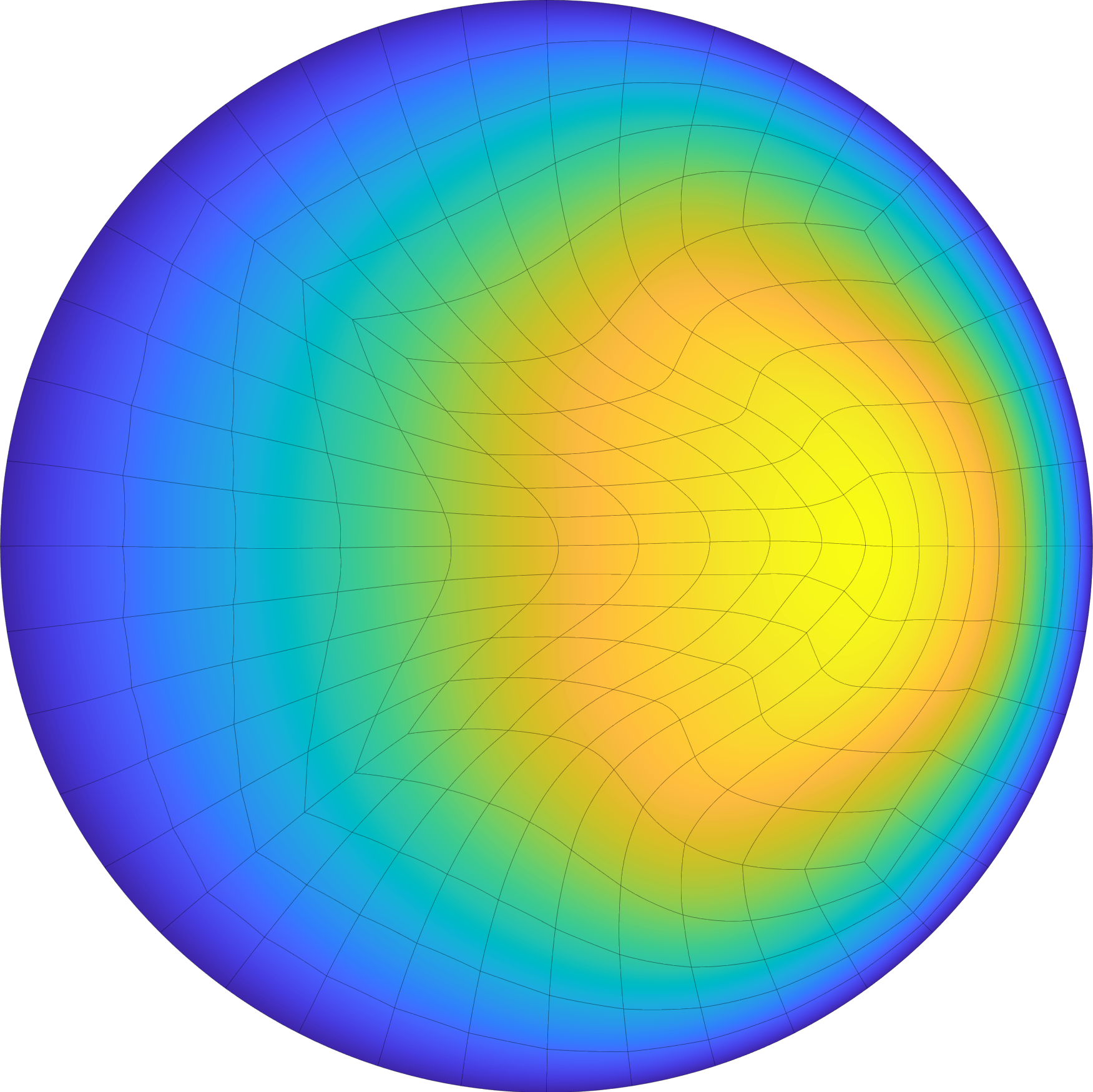}
\caption{Solution of the advection-diffusion equation (\ref{eqn:advdiff-ronq})
         at parameters $\mu=(0,10)$ (\textit{left}) and the corresponding
         approximation using residual minimization over the linear space
         $\Vcal_1$ (\textit{middle}) and nonlinear manifold $\Wcal_1$
         (ROM-IFT) (\textit{right}). The faint edges are included to
         illustrate the domain mapping employed by each method.}
\label{fig:demo1_hdm_vs_rom_vs_ift}
\end{figure}

We construct a one-dimensional linear space $\Vboldcal_1$
by taking the affine offset $\bar\Ubm=\zerobold$ and a
single mode $\Phibold_1$ as the normalized solution
$\Ubm(\text{Id},\mu)$ at $\mu=(0,0)$ (Figure~\ref{fig:dual_interp}).
This subspace is used to construct a traditional minimum-residual
ROM solution by freezing the domain mapping at $\Tcal_0=\text{Id}$,
i.e., $\Phibold_1\vbm_1(\text{Id},\mu)$,
and the ROM-IFT by simultaneously minimizing the HDM
residual over $\Vboldcal_1$ and the family of domain mappings
$\{\Tcal_s\mid s\in\Rbb\}$, i.e., $\Phibold_1\wbm_1(\mu)$
(Figure~\ref{fig:demo1_hdm_vs_rom_vs_ift}).
Recall the connection to the continuous level; traditional residual
minimization corresponds to searching for the solution over the
linear space $\Vcal_1$ (continuous counterpart to $\Vboldcal_1$),
while the ROM-IFT corresponds to searching for the solution over
the nonlinear manifold
$\Wcal_1 = \bigcup_{s\in\Rbb} \Vcal_1\circ\Tcal_s^{-1}$.
In this case, neither solution is highly
accurate; however, the solution from ROM-IFT method is
much better than the one from standard residual minimization
owing to the flexibility to deform the domain to track the
moving feature.

\begin{figure}
 \centering
 \begin{tikzpicture}
\begin{groupplot}[
  group style={
      group size=2 by 1,
      horizontal sep=1.7cm
  },
  width=0.48\textwidth,
  xlabel={$\alpha$},
  ylabel={$s$},
  xmin=10, xmax=80,
  ymin=-0.1, ymax=1.1
]

\nextgroupplot
\addplot graphics [xmin=10, xmax=80, ymin=-0.1, ymax=1.1] {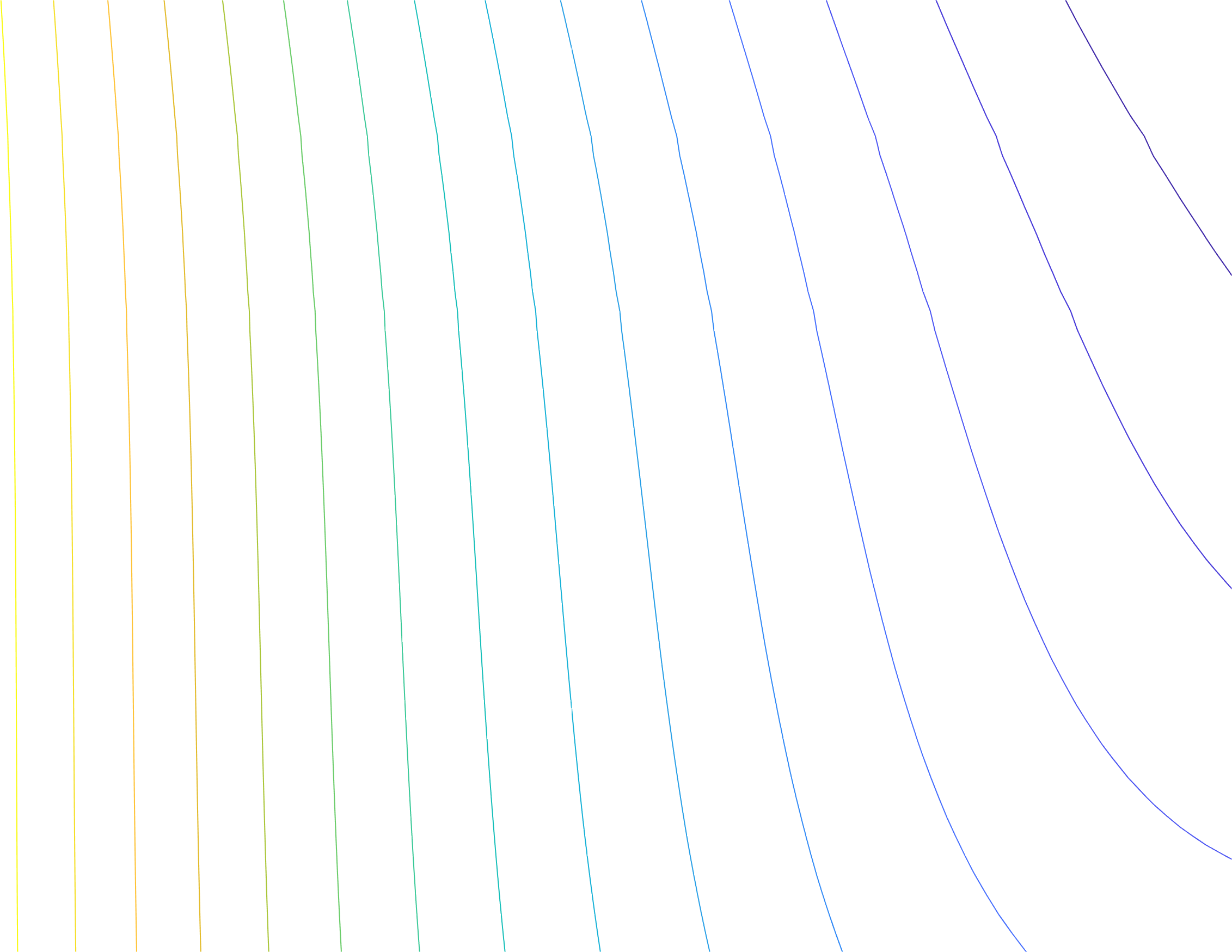};

\addplot [solid, black]
coordinates { (0, 0) (80, 0) }; \label{line:demo1_linsearch}

\addplot [solid, thick, black, mark=*, only marks]
coordinates { ( 18.0828, 0) }; \label{line:demo1_rom}

\addplot [solid, thick, blue, mark=square*, only marks]
coordinates { ( 65.1966, 1) }; \label{line:demo1_romift}

\nextgroupplot
\addplot graphics [xmin=10, xmax=80, ymin=-0.1, ymax=1.1] {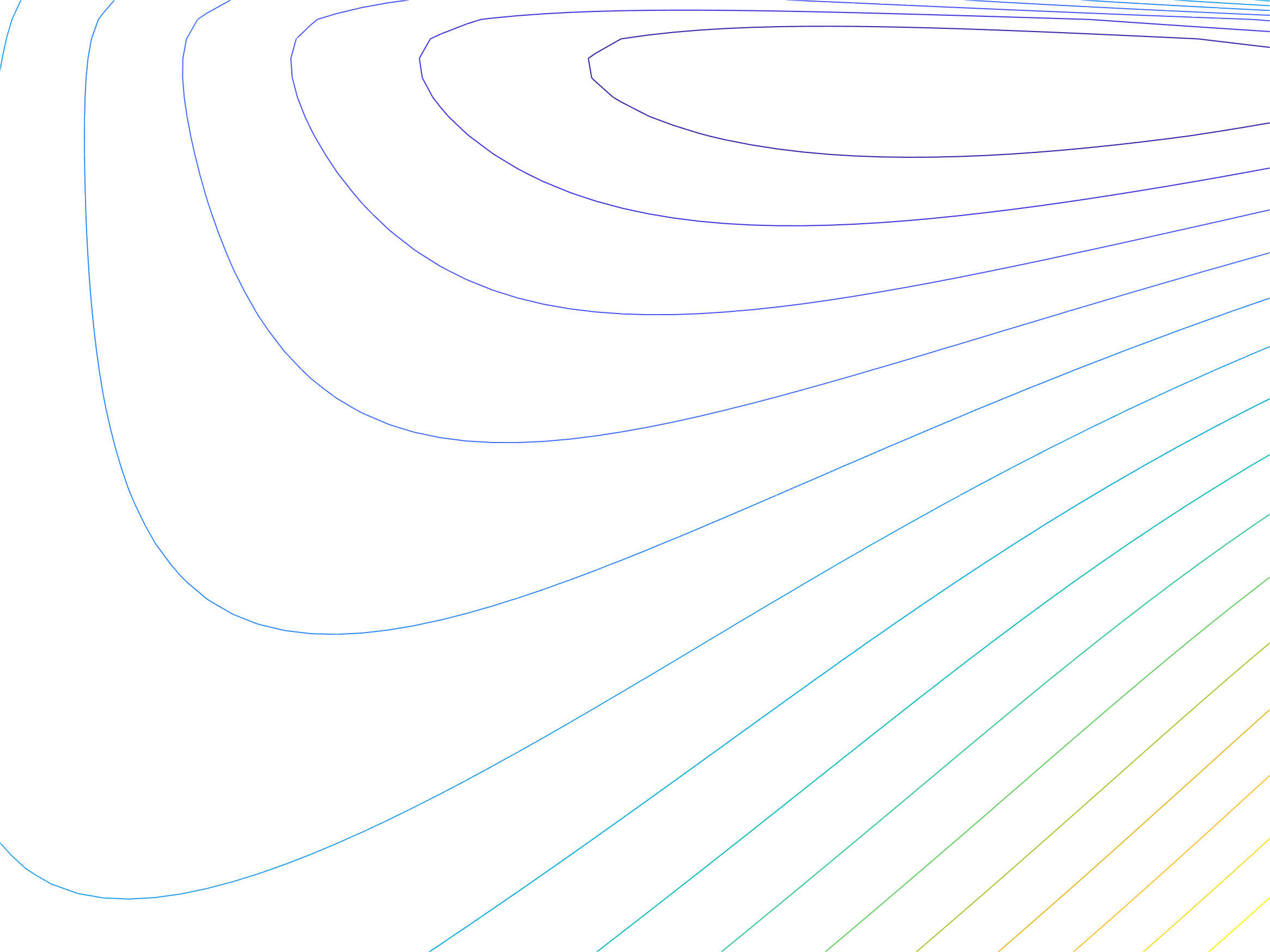};

\addplot [solid, black]
coordinates { (0, 0) (80, 0) };

\addplot [solid, thick, black, mark=*, only marks]
coordinates { ( 18.0828, 0) };

\addplot [solid, thick, blue, mark=square*, only marks]
coordinates { ( 65.1966, 1) };

\end{groupplot}
\end{tikzpicture}
\caption{Contours of $E_1$ (\textit{left}) and $E_2$ (\textit{right}).
 The search space for traditional residual minimization is
 $\{(\alpha,0)\mid \alpha\in\Rbb\}$ (\ref{line:demo1_linsearch}).
 The ROM (\ref{line:demo1_rom}) and ROM-IFT (\ref{line:demo1_romift})
 solutions are indicated on both plots.}
\label{fig:demo1_contour}
\end{figure}
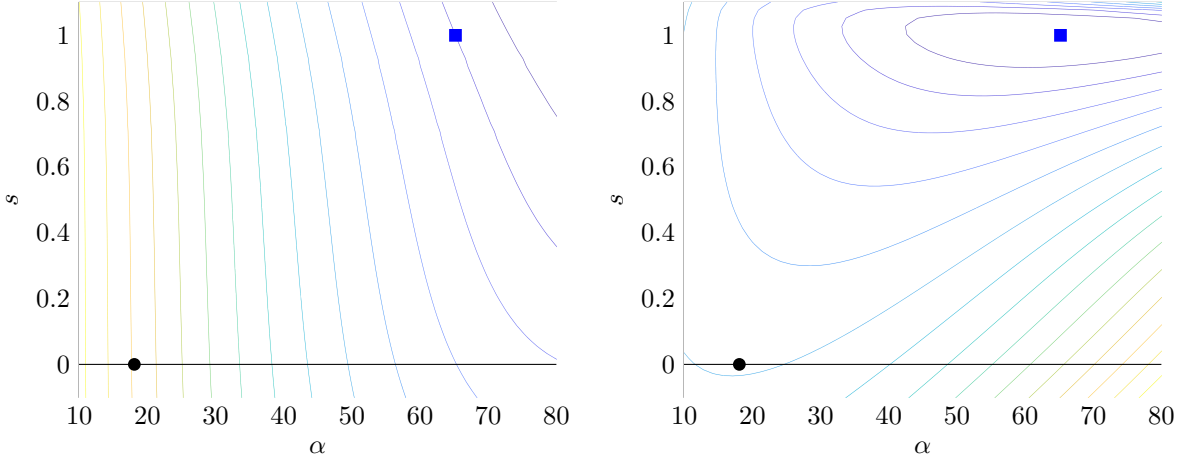

Finally, to further highlight the relationship between residual
minimization over a linear subspace and the proposed ROM-IFT method,
we define $\func{E_1}{\Rbb^2}{\Rbb}$ and
$\func{E_2}{\Rbb^2}{\Rbb}$ as
\begin{equation}
 E_1: (\alpha, s)\mapsto\frac{1}{2}\norm{\Ubm(\Tcal_s,\mu) - \alpha \Phibold_1}_2^2,
 \qquad
 E_2: (\alpha, s)\mapsto\frac{1}{2}\norm{\Rbm(\alpha\Phibold_1,\Tcal_s;\mu)}_2^2,
\end{equation}
for $\mu=(0,10)$. For a fixed $(\alpha,s)\in\Rbb^2$, $E_1(\alpha,s)$ is the
discrete $L^2$ norm of the error associated with a particular approximation
$\alpha\Phibold_1\in\Vboldcal_1$ to $\Ubm(\Tcal_s,\mu)$ and
$E_2(\alpha,s)$ is the corresponding residual norm.
For the case $\kappa=0$, the ROM-IFT method corresponds to
minimization of $E_2$ over $\Rbb^2$, whereas traditional
residual minimization corresponds to minimization of
$E_2(\,\cdot\,, 0)$ over $\Rbb$ (illustrated in
Figure~\ref{fig:demo1_contour} by the contours of
$E_1$ and $E_2$). Neither method locates the point
that minimizes the $E_1$; however, the ROM-IFT solution
is much closer (Figure~\ref{fig:demo1_contour}).


\subsection{Offline training: generation and compression
            of feature-aligned snapshots}
\label{sec:rom:off}
To construct the reduced basis $\Phibold_k$, we use the method of snapshots
\cite{sirovich1987turbulence} by collecting and compressing
solutions of the HDM $\dstvc(\dommap,\param)$ for some values
of $\dommap\in\dommapsp$ and $\param\in\paramsp$.
In this work, we assume an ordered collection of parameters
$\Dcal_\text{tr}=\{\param_1,\dots,\param_M\}$
is defined \textit{a priori}, e.g., based on uniform or Latin Hypercube
sampling; future work will consider efficient sampling strategies, e.g.,
greedy methods \cite{prud2002reliable}, where samples are chosen adaptively.
With the training set identified, snapshots are usually constructed on a
fixed domain defined by $\hat\dommap$ (usually identity), i.e.,
$\dstvc(\hat\dommap,\param_1),\dots,\dstvc(\hat\dommap,\param_M)$.
However, for convection-dominated problems, each of these solutions will, in
general, will have features in different spatial positions and orientations,
making them suboptimal snapshots (e.g., the numerical experiment in
Section~\ref{sec:rom:offcmpr}). Instead,
we choose to construct a domain mapping $\hat\dommap_i$ for each snapshot
$\param_i$ such that, in the reference domain, corresponding features align,
a procedure we will call \textit{snapshot alignment} in the remainder. To
this end, we choose the space of domain mappings to be
$\dommapsp=\dommapsp_{h,q}^\mathrm{b}$ to ensure a sufficiently rich
space of domain mappings to induce alignment of the snapshot features while
preserving the boundaries of the domain. After the first $k\leq M$
snapshots have been aligned, the reduced basis is defined as
\begin{equation} \label{eqn:basis}
 \Phibold_{n_k} \coloneqq
 \text{POD}_{n_k}\left(
 \begin{bmatrix}
  \dstvc(\hat\dommap_1,\param_1)-\bar\Ubm & \cdots &
  \dstvc(\hat\dommap_k,\param_k)-\bar\Ubm
 \end{bmatrix}
 \right),
\end{equation}
where $\bar\Ubm\in\Rbb^N$ is the affine offset (taken to be
the zero vector in this work),  $n_k$ is the number of vectors
retained in the basis after $k$ training points have been processed,
and $\func{\text{POD}_n}{\Rbb^{N\times k}}{\Rbb^{N\times n}}$
applies the singular value decomposition to the input matrix
and extracts the $n$ left singular vectors, i.e., the proper
orthogonal decomposition of the input matrix.

Snapshot alignment is performed by iteratively bootstrapping the
ROM-IFT in (\ref{eqn:romift3}). Because we are only interested
in relative alignment of the features in the snapshots and not the
actual location of these features in the reference domain, we are free
to choose $\hat\dommap_1$ arbitrarily; in this work, we always
take $\rdom=\pdom$ and choose $\hat\dommap_1=\text{Id}$. From
this, we compute the one-dimensional reduced basis according
to (\ref{eqn:basis}) to obtain $\Phibold_{n_1}$. Assuming $k-1$
snapshots have been aligned and the basis $\Phibold_{n_{k-1}}$
has been computed according to (\ref{eqn:basis}), the $k$th snapshot
is aligned by solving the ROM-IFT problem in (\ref{eqn:romift3}) and
$\hat\dommap_k$ is taken to be the resulting domain mapping, i.e.,
\begin{equation} \label{eqn:romift_off1}
 (\hat\wbm_k,\hat\dommap_k) \coloneqq \argmin_{(\wbm,\dommap)\in\Rbb^{n_{k-1}}\times\dommapsp_{h,q}^\mathrm{b}} J(\bar\Ubm+\Phibold_{n_{k-1}}\wbm,\dommap;\param_k).
\end{equation}
The reduced coordinates $\hat\wbm_k$ are discarded and the basis is updated
to $\Phibold_{n_k}$ according to (\ref{eqn:basis}) after computing the HDM
solution.
This procedure is applied recursively to align features in all $M$ snapshots
resulting in the basis $\Phibold_{n_M}$ that is available for use in the
online phase. Because the space $\dommapsp_{h,q}^\mathrm{b}$ is potentially
large and not amenable to efficient online computations, we use the domain
mapping snapshots $\{\hat\dommap_1,\dots,\hat\dommap_M\}$ to define
a low-dimensional space of admissible domain mapping,
$\Gbb_{h,q}^{\mathrm{b,r}}\subset\Gbb_{h,q}^\mathrm{b}$, defined
as
\begin{equation}
 \Gbb_{h,q}^\mathrm{b,r}\coloneqq
 \left\{\Gcal_{h,q}^\mathrm{b}(\,\cdot\,;\hat\ybm_1 + \Psibold_{n'} \cbm) \suchthat \cbm\in\Rbb^{n'}\right\},
\end{equation}
where $\{\hat\ybm_1,\dots,\hat\ybm_M\}\subset\Rbb^{N_\mathrm{u}}$ are the
unconstrained degrees of freedom defining the domain mapping snapshots,
i.e., $\hat\dommap_i = \dommap_{h,q}^\mathrm{b}(\,\cdot\,; \hat\ybm_i)$,
$\Psibold_{n'}\in\Rbb^{N_\mathrm{u}\times n'}$ is a reduced basis
for the \textit{unconstrained} mapping degrees of freedom
\begin{equation} \label{eqn:dombasis}
 \Psibold_{n'} \coloneqq
 \text{POD}_{n'} \left(
 \begin{bmatrix}
  \hat\ybm_2-\hat\ybm_1 & \cdots & \hat\ybm_M-\hat\ybm_1
 \end{bmatrix}
 \right),
\end{equation}
and $n'$ is the number of vectors retained in the basis after
POD truncation. Then, for any element $\dommap\in\dommapsp_{h,q}^\mathrm{b,r}$,
there exists $\cbm\in\Rbb^{n'}$ such that
$\dommap=\dommap_{h,q}^\mathrm{b,r}(\,\cdot\,;\cbm)$, where
\begin{equation}
 \cbm \mapsto \dommap_{h,q}^\mathrm{b,r}(\,\cdot\,;\cbm) \coloneqq
 \Gcal_{h,q}^\mathrm{b}(\,\cdot\,; \hat\ybm_1+\Psibold_{n'}\cbm).
\end{equation}
The low-dimensional space of domain mappings
is constructed by reducing the unconstrained degrees of freedom
and composing with the boundary-preserving mapping $\chibold$ to
ensure all elements of $\Gbb_{h,q}^\mathrm{b,r}$ preserve the
boundaries of the domain (to high-order). The offline training
procedure is summarized in Algorithm~\ref{alg:romift_off}.
\begin{algorithm}
 \caption{Generation and compression of feature-aligned snapshots}
 \label{alg:romift_off}
 \begin{algorithmic}[1]
  \REQUIRE Training parameters $\{\mu_1,\dots,\mu_M\}$
  \ENSURE State basis $\Phibold_{n_M}$, unconstrained domain mapping
   basis $\Psibold_{n'}$
  \STATE \textbf{Initialize state basis}:
   $\Phibold_{n_1} = \Ubm(\mathrm{Id},\mu_1)/\norm{\Ubm(\mathrm{Id},\mu_1)}$
  \FOR{$k=2,\dots,M$}
   \STATE \textbf{Snapshot alignment}: Compute $\hat\Gcal_k$ from
    (\ref{eqn:romift_off1})
   \STATE \textbf{HDM solution}: Compute HDM state corresponding to
    domain mapping $\hat\Gcal_k$, i.e., $\Ubm(\hat\Gcal_k,\mu_k)$
   \label{alg:line:solve_hdm}
   \STATE \textbf{State basis update}: Compute new basis $\Phibold_{n_k}$
    according to (\ref{eqn:basis})
  \ENDFOR
  \STATE \textbf{Domain mapping basis}: Compute unconstrained domain mapping
   basis $\Psibold_{n'}$ according to (\ref{eqn:dombasis})
\end{algorithmic}
\end{algorithm}

\begin{remark}
In this work, we usually choose $n_k = k$ (no truncation) due to the
limited training required to obtain high accuracy and parametric
robustness when snapshots are aligned and the PDE solution is
approximated using the nonlinear manifold based on deforming
modes (Section~\ref{sec:rom:ift}). In this setting, POD is
not required in (\ref{eqn:basis}), the snapshots are simply
normalized to obtain an orthonormal basis.
\end{remark}

\begin{remark}
If the HDM is based on an $r$-adaptive method,
Line~\ref{alg:line:solve_hdm} of Algorithm~\ref{alg:romift_off}
(HDM solve) will produce the PDE state $\Ubm(\mu_k)$ as well as the
corresponding mapping $\Gcal(\mu_k)$. In this case, we replace
$\hat\Gcal_k \leftarrow \Gcal(\mu_k)$ in the construction of
the low-dimensional space of domain mappings.
\end{remark}


\begin{remark}
The proposed offline phase is more expensive than the offline
phase usually associated with the method of snapshots due to
the additional snapshot alignment step, which is required in addition
to standard offline computations (HDM solution, snapshot compression).
However, additional cost of snapshot alignment is relatively small because
it leverages the low-dimensional manifold approximation and ROM-IFT method.
Furthermore, as we show in Section~\ref{sec:numexp}, the ROM-IFT method
requires far fewer training parameters to build an accurate, parametrically
robust surrogate, which further offsets this additional expense.
\end{remark}


\begin{remark}
Other objective functions in place of (\ref{eqn:romift5})
were investigated for use in the offline phase, most notably the error
in the $L^2$ projection of the HDM solution
\begin{equation}
 J_\text{err}(\Vbm,\dommap;\param) =
 \frac{1}{2}\norm{\Phibold_k^T(\dstvc(\dommap,\param)-\bar\dstvc)}_2^2.
\end{equation}
However, this option was more computationally expensive because
each evaluation of $J_\text{err}$ requires a HDM solve and
provided no discernible benefit over (\ref{eqn:romift5}). Thus, we use
(\ref{eqn:romift5}) solely in this work as it is inexpensive and cleanly
unifies the offline and online phases of model reduction.
\end{remark}


\subsection{Comparison: basis construction from non-aligned
            vs. aligned snapshots}
\label{sec:rom:offcmpr}
\begin{figure}
\centering
\input{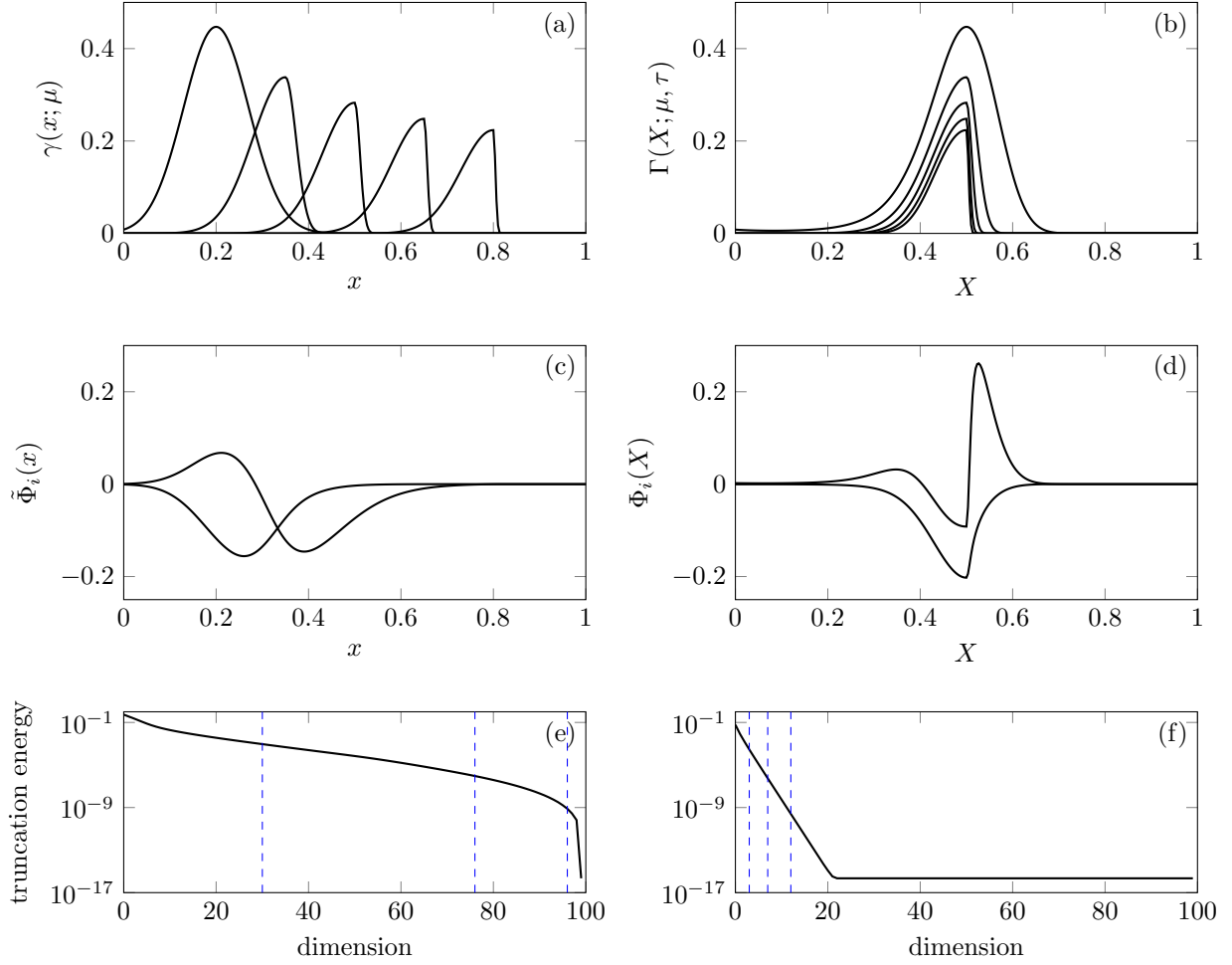}
\caption{Several instances of the steepening Gaussian
 $\gamma(\,\cdot\,;\mu)$ for $\mu\in[0.2,0.8]$ (a) and
 the corresponding mapped functions $\Gamma_\mu(\,\cdot\,\mu)$
 that centers the Gaussian at $X=0.5$ (b).
 The first two POD modes based on the non-aligned snapshots
 $\{\tilde{\Phi}_i\}_{i=1}^2$ (c) and the aligned
 snapshots $\{\Phi_i\}_{i=1}^2$ (d).
 Singular value decay for the non-aligned (e) and aligned
 (f) snapshots; the vertical lines correspond to truncation
 energies of $10^{-3}$, $10^{-6}$, and $10^{-9}$.}
\label{fig:demo1d1_compr}
\end{figure}

To demonstrate the merit of using feature-aligned snapshots
in Section~\ref{sec:rom:off} to construct a reduced basis,
consider the piecewise Gaussian function
$\func{\gamma}{[0,1]\times[0.2,0.8]}{\Rbb_{>0}}$ defined as
\begin{equation}
 (x,\mu)\mapsto
 \gamma(x;\mu) \coloneqq \theta(x;(0.2\mu^{-1/2}, 0.1, \mu)) +
                         \theta(-x;(0.2\mu^{-1/2}, 0.004\mu^{-2}, -\mu)),
\end{equation}
where $\theta$ is the cutoff Gaussian in (\ref{eqn:gausscut}),
that advects and steepens as the parameter $\mu\in[0.2,0.8]$ increases
(Figure~\ref{fig:demo1d1_compr}a).
The first two POD modes $\{\tilde\Phi_i\}_{i=1}^2$
extracted from the collection of $100$ snapshots
$\{\gamma(\,\cdot\,;\mu_i)\}_{i=1}^{100}$,
where the parameters $\{\mu_i\}_{i=1}^{100}$ are uniformly spaced
in $[0.2,0.8]$, target the regions of the solution manifold
corresponding to smooth, large amplitude modes and fail to
capture the advecting feature
(Figure~\ref{fig:demo1d1_compr}c).
Furthermore, the singular value decay associated with the
snapshots is slow, which indicates that the solutions do
not lie in a low-dimensional linear space
(Figure~\ref{fig:demo1d1_compr}e).
Alternatively, consider the
mapped version of the steepening Gaussian
$\func{\Gamma_\tau}{[0,1]\times[0.2,0.8]}{\Rbb_{>0}}$ defined as
\begin{equation}
 (X,\mu)\mapsto
 \Gamma_\tau(X;\mu) \coloneqq
 \gamma(\Lcal_\tau(X);\mu),
\end{equation}
where the domain mapping $\func{\Lcal_\tau}{[0,1]}{[0,1]}$ is defined as
\begin{equation}
 X\mapsto
 \Lcal_\tau(X)\coloneqq
 X + 4(\tau-0.5)X(1-X).
\end{equation}
By choosing $\tau=\mu$, mapped Gaussians will be centered at $X=0.5$;
that is, the advecting feature becomes a standing feature
(Figure~\ref{fig:demo1d1_compr}b).
The first two POD modes $\{\Phi_i\}_{i=1}^2$
extracted from the collection of $100$ snapshots
$\{\Gamma_{\mu_i}(\,\cdot\,;\mu_i)\}_{i=1}^{100}$ effectively
captures the standing feature and smooth variations around it
(Figure~\ref{fig:demo1d1_compr}d).
Furthermore, the singular value decay associated with the
snapshots is rapid indicating the mapped solutions lie in
a low-dimensional linear space
(Figure~\ref{fig:demo1d1_compr}f).
From this simple example, it is preferable to construct
reduced basis in the reference domain from mapped snapshots
where features align.

\subsection{Online computations}
\label{sec:rom:on}
With the basis $\Phibold_k$ ($k=n_M$ for brevity) defined in the offline phase
as well as the low-dimensional space of boundary-preserving domain mappings
$\dommapsp_{h,q}^\mathrm{b,r}$, the online phase uses the ROM-IFT
method (\ref{eqn:romift3}) to efficiently approximate the PDE solution
at any parameter in a test set
$\param\in\Dcal_\text{tst}\subset\paramsp$. That is, for each
$\param\in\Dcal_\text{tst}$, the following ROM-IFT problem is solved
\begin{equation} \label{eqn:romift_on1}
 (\wbm_k(\mu),\dommap_k(\,\cdot\,;\mu)) \coloneqq
 \argmin_{(\wbm,\dommap)\in\Rbb^k\times\dommapsp_{h,q}^\mathrm{b,r}}
 J(\bar\Ubm+\Phibold_k\wbm,\dommap;\param).
\end{equation}
The PDE solution in the physical domain is then approximated as
\begin{equation}
 \pstvc(\,\cdot\,;\param) \approx
 \Upsilon(\bar\Ubm+\Phibold_k\wbm_k(\mu))\circ\dommap_k(\,\cdot\,;\param)^{-1}.
\end{equation}

To apply numerical optimization techniques to solve the ROM-IFT problem,
we first convert it to its fully discrete form.
Because any element $\dommap\in\dommapsp_{\mshsz,\mshord}^\text{b,r}$ can be
expressed as $\dommap=\dommap_{\mshsz,\mshord}^\text{b,r}(\,\cdot\,;\cbm)$
for some $\cbm\in\Rbb^{n'}$, the ROM-IFT problem in (\ref{eqn:romift_on1}) is
equivalent to
\begin{equation}
 \optunc{(\wbm,\cbm)\in\Rbb^k\times\Rbb^{n'}}{\tilde{J}_k(\wbm,\cbm;\param)\coloneqq J(\bar\dstvc+\Phibold_k\wbm,\dommap_{h,q}^\text{b,r}(\,\cdot\,;\cbm);\param)}.
\end{equation}
Because $k+n' \ll N$, this is a small optimization problem compared to the
dimension of the HDM. Finally, we define the fully discrete objective
function in residual form
\begin{equation}
 (\wbm,\cbm;\param) \mapsto \tilde\Fbm_k(\wbm,\cbm;\param)
 \coloneqq \Fbm(\bar\dstvc+\Phibold_k\wbm,\dommap_{h,q}^\text{b,r}(\,\cdot\,;\cbm);\param),
\end{equation}
which will be used extensively in the next section to define the
Levenberg-Marquardt solver, and we denote the optimal solution to the
ROM-IFT problem as
\begin{equation}\label{eqn:romift8}
 (\wbm_k(\param),\cbm_k(\param)) \coloneqq
 \argmin_{(\wbm,\cbm)\in\Rbb^k\times\Rbb^{n'}}{\tilde{J}_k(\wbm,\cbm;\param)=\frac{1}{2}\norm{\tilde\Fbm_k(\wbm,\cbm;\param)}_2^2}
\end{equation}
which is consistent with the definition of $\wbm_k(\param)$ in
(\ref{eqn:romift3}).


\section{Full space solver for implicit feature tracking}
\label{sec:solver}
In this section, we describe a standard Levenberg-Marquardt method
\cite{nocedal_numerical_2006}
for solving the fully discrete ROM-IFT problem (Section~\ref{sec:solver:gn})
for a fixed $\param\in\paramsp$, including practical considerations
such as termination criteria (Section~\ref{sec:solver:optcond}) and
initialization (Section~\ref{sec:solver:init}). We use the notation
of the online ROM-IFT problem (\ref{eqn:romift8}); however, due to
the unifying structure between the offline and online optimization
problems, the same algorithm can be used for the offline snapshot
alignment.

\subsection{Levenberg-Marquardt nonlinear least-squares solver}
\label{sec:solver:gn}
The globally convergent Levenberg-Marquardt method to solve the
nonlinear least-square problem in (\ref{eqn:romift8}) defines a
sequence of iterates $\{\wbm_k^{(n)}\}_{n=0}^\infty\subset\Rbb^{k}$
and $\{\cbm_k^{(n)}\}_{n=0}^\infty\subset\Rbb^{n'}$ as
\begin{equation}
\begin{aligned}
 \wbm_k^{(n+1)} &= \wbm_k^{(n)} + \alpha_{n+1} \Delta\wbm_k^{(n+1)} \\
 \cbm_k^{(n+1)} &= \cbm_k^{(n)} + \alpha_{n+1} \Delta\cbm_k^{(n+1)},
\end{aligned}
\qquad n=0,1,\dots
\end{equation}
where $\wbm_k^{(0)}$ and $\cbm_k^{(0)}$ initialize the
iteration (Section~\ref{sec:solver:init}) and $\alpha_{n+1}\in\Rbb_{>0}$
is the step length chosen to ensure sufficient decrease
to the solution based on the Wolfe conditions \cite{nocedal_numerical_2006}.
The steps $\Delta\wbm_k^{(n+1)}$ and $\Delta\cbm_k^{(n+1)}$
are defined as the solution of the linear least-squares
problem
\begin{equation} \label{eqn:linlstsq}
 (\Delta\wbm_k^{(n+1)}, \Delta\cbm_k^{(n+1)}) \coloneqq
 \argmin_{(\Delta\wbm,\Delta\cbm)\in\Rbb^k\times\Rbb^{n'}}
 \norm{\begin{bmatrix} \tilde\Fbm_k^{(n)} \\ \zerobold \end{bmatrix} +
       \begin{bmatrix} \tilde\Jbm_{\wbm,k}^{(n)} & \tilde\Jbm_{\cbm,k}^{(n)} \\
                       \zerobold & \sqrt\lambda\Ibm_{n'}
       \end{bmatrix}
       \begin{bmatrix} \Delta \wbm \\ \Delta \cbm \end{bmatrix}}_2,
\end{equation}
where
\begin{equation}
 \tilde\Fbm_k^{(n)} \coloneqq
 \tilde\Fbm_k(\wbm_k^{(n)},\cbm_k^{(n)};\param), \qquad
 \tilde\Jbm_{\wbm,k}^{(n)} \coloneqq
 \pder{\tilde\Fbm_k}{\wbm}(\wbm_k^{(n)},\cbm_k^{(n)};\param), \qquad
 \tilde\Jbm_{\cbm,k}^{(n)} \coloneqq
 \pder{\tilde\Fbm_k}{\cbm}(\wbm_k^{(n)},\cbm_k^{(n)};\param)
\end{equation}
and $\lambda\in\Rbb$ is a regularization parameter chosen based
on the strategy in \cite{huang2021robust}.
\begin{remark}
The Levenberg-Marquardt regularization in (\ref{eqn:linlstsq}) is
only applied to the domain mapping variables ($\cbm$), which follows
the approach taken in \cite{corrigan_moving_2019, zahr_implicit_2020}.
\end{remark}

The analytical solution of the linear least-squares problems leads
to the following expression for the step (normal equations)
\begin{equation}
 \tilde\Jbm_k^{(n),T} \tilde\Jbm_k^{(n)}
 \Delta\zbm_k^{(n+1)} = -
 \begin{bmatrix}
  \tilde\Jbm_{\wbm,k}^{(n),T} \tilde\Fbm_k^{(n)} \\
  \tilde\Jbm_{\cbm,k}^{(n),T} \tilde\Fbm_k^{(n)}
 \end{bmatrix},
\end{equation}
where
\begin{equation}
\tilde\Jbm_k^{(n)} =
\begin{bmatrix}
 \tilde\Jbm_{\wbm,k}^{(n)} & \tilde\Jbm_{\cbm,k}^{(n)} \\
 \zerobold & \sqrt{\lambda}\Ibm_{n'}
\end{bmatrix}, \qquad
\Delta \zbm_k^{(n+1)} =
\begin{bmatrix}
 \Delta\wbm_k^{(n+1)} \\
 \Delta\cbm_k^{(n+1)}
\end{bmatrix}.
\end{equation}
Due to the poor conditioning of this linear system, it is preferable
to directly solve the linear least-squares problem in (\ref{eqn:linlstsq})
using the QR factorization \cite{nocedal_numerical_2006}.
Because $\tilde\Jbm_k^{(n)}$ is a small ($k+n'\ll N$), dense matrix,
the cost of computing the step $\Delta\zbm_k^{(n+1)}$ is relatively small.

\subsection{Optimality conditions and termination criteria}
\label{sec:solver:optcond}
The first-order optimality condition of the (unconstrained)
nonlinear least-squares problem in (\ref{eqn:romift8}) requires the
gradient of $J$ to vanish at the solution $(\wbm_k(\param),\cbm_k(\param))$,
i.e.,
\begin{equation}
 \begin{bmatrix}
  \nabla_\wbm \tilde\Fbm_k(\wbm_k(\param),\cbm_k(\param);\param) \\
  \nabla_\cbm \tilde\Fbm_k(\wbm_k(\param),\cbm_k(\param);\param)
 \end{bmatrix} =
 \zerobold,
\end{equation}
which leads to the optimality condition
\begin{equation}
\begin{aligned}
\pder{\tilde\Fbm_k}{\wbm}(\wbm_k(\param),\cbm_k(\param);\param)^T
\tilde\Fbm_k(\wbm_k(\param),\cbm_k(\param);\param) &= \zerobold \\
\pder{\tilde\Fbm_k}{\cbm}(\wbm_k(\param),\cbm_k(\param);\param)^T
\tilde\Fbm_k(\wbm_k(\param),\cbm_k(\param);\param) &= \zerobold.
\end{aligned}
\end{equation}
Thus, given tolerance $\epsilon_1,\epsilon_2>0$,
$\left(\wbm_k^{(n)},\cbm_k^{(n)}\right)$ is a considered a
numerical solution of (\ref{eqn:romift8}) if
\begin{equation}
\norm{\tilde\Jbm_{\wbm,k}^{(n),T}\tilde\Fbm_k^{(n)}} \leq \epsilon_1, \qquad
\norm{\tilde\Jbm_{\cbm,k}^{(n),T}\tilde\Fbm_k^{(n)}} \leq \epsilon_2.
\end{equation}
In this work, we take $\epsilon_1=\epsilon_2= 10^{-8}$ unless otherwise
stated.

\subsection{Solver initialization}
\label{sec:solver:init}
The Levenberg-Marquardt solver is initialized by taking the domain mapping
to be the fixed nominal mapping $\bar\dommap$, i.e., $\cbm_k^{(0)}=\zerobold$,
and the corresponding minimum-residual solution over the linear space
$\Vboldcal_k$ (fixed mapping) for the PDE state, i.e.,
$\wbm_{k,0} = \vbm_{k}(\bar\dommap, \param)$.

\section{Numerical experiments}
\label{sec:numexp}
In this section, we investigate the performance of the proposed ROM-IFT
method relative to standard model reduction techniques on a fixed domain
on three computational fluid dynamics benchmark problems. For the first
two problems (linear advection-reaction and inviscid, compressible flow
through a nozzle), we demonstrate deep convergence of the proposed
ROM-IFT as the reduced basis and parameter space sampling are refined.
The last problem (supersonic flow over a cylinder) demonstrates the ability
of the ROM-IFT method to provide accurate CFD approximations with limited
sampling of the parameter space.

For each problem, we will compare the accuracy of fixed-domain model
reduction and model reduction with implicit feature tracking; the
accuracy of each method will be measured by the relative $L^2(\rdom)$
error with respect to the corresponding HDM solution (or exact solution,
if available).
For any $\mu\in\Dcal$ and $\Gcal\in\Gbb$, the fixed-domain ROM solution
(dimension $k$) $U_k(\,\cdot\,;\Gcal,\mu)$, defined in (\ref{eqn:approx-ref})
and either (\ref{eqn:gal}) (Galerkin) or (\ref{eqn:minres}) (minimum-residual),
approximates the HDM solution $U_h(\,\cdot\,;\Gcal,\mu)$, defined in
(\ref{eqn:hdm})-(\ref{eqn:hdm1}). For all
cases considered, the domain mapping is frozen at the nominal
map $\Gcal=\bar\Gcal$ and the reduced basis is constructed
by applying POD to the HDM snapshots
$\{U_h(\,\cdot\,;\bar\Gcal,\mu) \mid \mu\in\Dcal_\text{tr}\}$,
where $\Dcal_\text{tr}\subset\Dcal$ is a finite training set. In this setting,
the relative $L^2(\rdom)$ error is the following function over $\Dcal$
\begin{equation} \label{eqn:err-rom}
 \mu\in\Dcal\mapsto e_k^\text{rom}(\mu) \coloneqq
 \sqrt{\frac{\int_{\rdom} \norm{U_h(\rcoord;\bar\Gcal(\rcoord),\mu)-U_k(\rcoord;\bar\Gcal(\rcoord),\mu)}_2^2 \, dV}{\int_{\rdom} \norm{U_h(\rcoord;\bar\Gcal(\rcoord),\mu)}_2^2 \, dV}}.
\end{equation}
For the examples considered in this work, fixed-domain ROMs based on a
Galerkin projection and residual minimization perform similarly;
for brevity, we only compare to residual minimization.

For any $\mu\in\Dcal$, the ROM-IFT solution
$(U_k(\,\cdot\,;\mu),\Gcal_k(\,\cdot\,;\mu))$ (\ref{eqn:romift2}) is defined
such that $U_k(\,\cdot\,;\mu)$ approximates the HDM solution
$U_h(\,\cdot\,;\Gcal_k(\,\cdot\,;\mu),\mu)$; unlike fixed-domain
model reduction, the implicit feature tracking method chooses the
most suitable domain deformation for each parameter $\mu$
($\Gcal_k(\,\cdot\,;\mu)$) and thus approximates the corresponding
HDM solution ($U_h(\,\cdot,\,;\Gcal_k(\mu), \mu)$).
In Section~\ref{sec:numexp:advec},  the HDM discretization is chosen such
that the PDE solution is well-resolved for any reasonable domain deformation,
i.e., a good approximation to $U(\,\cdot\,;\Gcal,\mu)$ for any $\Gcal\in\Gbb$
(Remark~\ref{rem:hdm-resolved}), which implies the loss in resolution
associated with evaluating the HDM solution at the mapping
$\Gcal_k(\,\cdot\,;\mu)$ is negligible.
In this setting, the relative $L^2(\rdom)$ error is the following
function over $\Dcal$
\begin{equation} \label{eqn:err-ift}
 \mu\in\Dcal\mapsto e_k^\text{ift}(\mu) \coloneqq
 \sqrt{\frac{\int_{\rdom} \norm{U_h(\rcoord;\Gcal_k(\rcoord;\mu),\mu)-U_k(\rcoord;\mu)}_2^2 \, dV}{\int_{\rdom} \norm{U_h(\rcoord;\Gcal_k(\rcoord;\mu),\mu)}_2^2 \, dV}}.
\end{equation}
In Sections~\ref{sec:numexp:nozzle}-\ref{sec:numexp:euler-cyl0}, the ROM-IFT
method is compared directly to an exact or reference solution rather than
the HDM. For all cases considered, the reduced basis is constructed by
applying POD to \textit{aligned} HDM snapshots (Section~\ref{sec:rom:off})
from a finite training set $\Dcal_\mathrm{tr}\subset\Dcal$.

To assess the parametric performance of fixed-domain model reduction
and implicit feature tracking, we use the maximum relative
$L^2(\rdom)$ error over a finite test set $\Dcal_\mathrm{tst}$, defined as
\begin{equation}
 \Dcal_\mathrm{tst}\subset\Dcal \mapsto E_k^\text{rom}(\Dcal_\mathrm{tst}) \coloneqq
 \max_{\mu\in\Dcal_\mathrm{tst}}~e_k^\text{rom}(\mu), \qquad\qquad
 \Dcal_\mathrm{tst}\subset\Dcal \mapsto E_k^\text{ift}(\Dcal_\mathrm{tst}) \coloneqq
 \max_{\mu\in\Dcal_\mathrm{tst}}~e_k^\text{ift}(\mu).
\end{equation}
The superscript is dropped, i.e., $e_k(\mu)$ or $E_k(\Dcal)$, to denote
the error of either the fixed-domain ROM or ROM-IFT. Finally,
$\sigma_k(\Dcal_\mathrm{tr})$ denotes $k$th singular value associated with the
collection of snapshots associated with the training set $\Dcal_\mathrm{tr}$.


\subsection{Advection-reaction}
\label{sec:numexp:advec}
\begin{figure}
\centering
\begin{tikzpicture}
\begin{groupplot}[
  group style={
      group size=3 by 1,
      horizontal sep=1cm
  },
  width=0.39\textwidth,
  axis equal image,
  xlabel={$x_1$},
  ylabel={$x_2$},
  xtick = {0.0, 0.5, 1.0},
  ytick = {0.0, 0.5, 1.0},
  xmin=0, xmax=1,
  ymin=0, ymax=1
]
\nextgroupplot[title={$\mu=(-\pi/10, 0.3, 60)$}]
\addplot graphics [xmin=0, xmax=1, ymin=0, ymax=1] {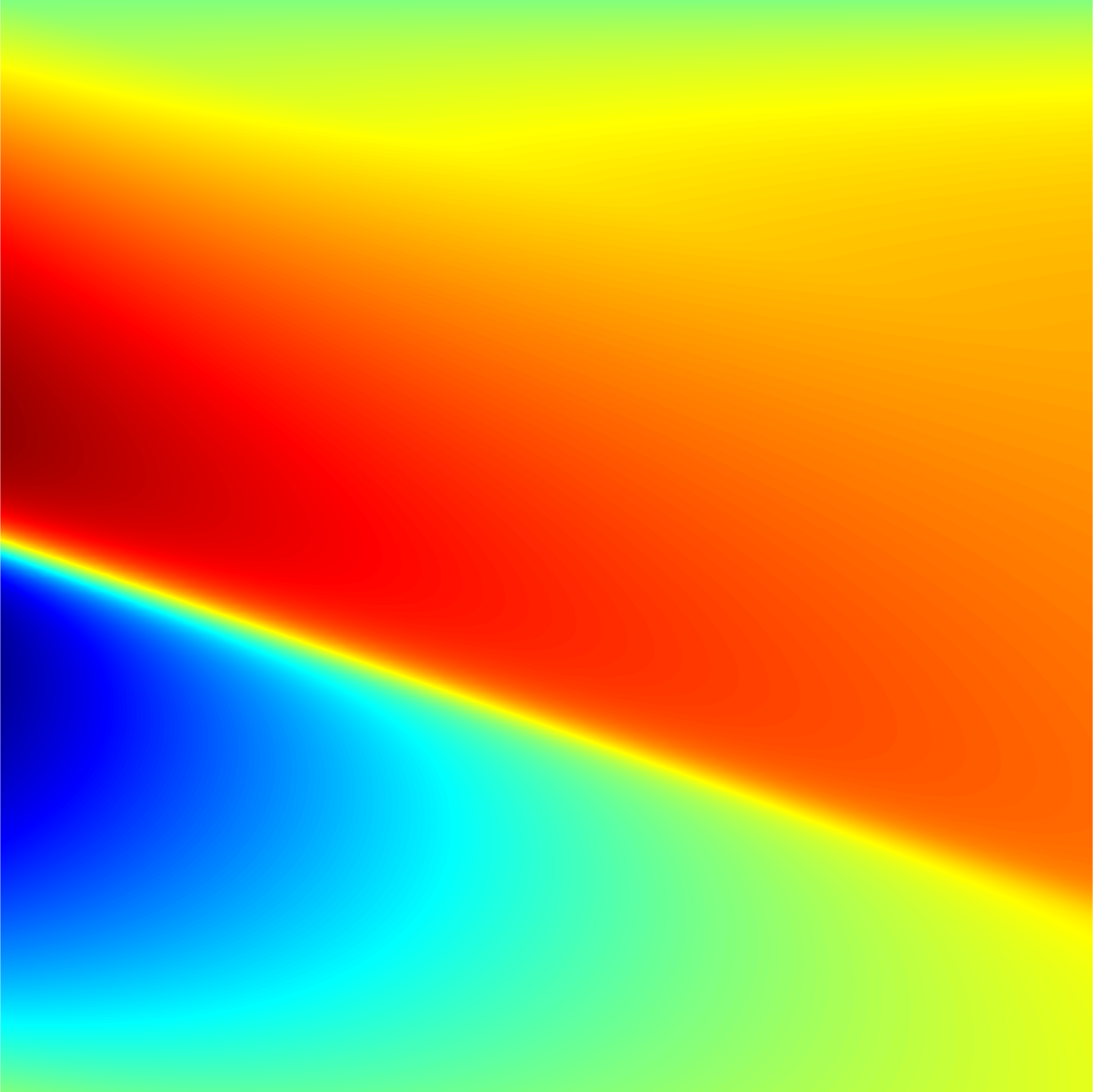};

\nextgroupplot[title={$\mu=(\pi/10, 0.7, 100)$}, ylabel={}, ytick=\empty]
\addplot graphics [xmin=0, xmax=1, ymin=0, ymax=1] {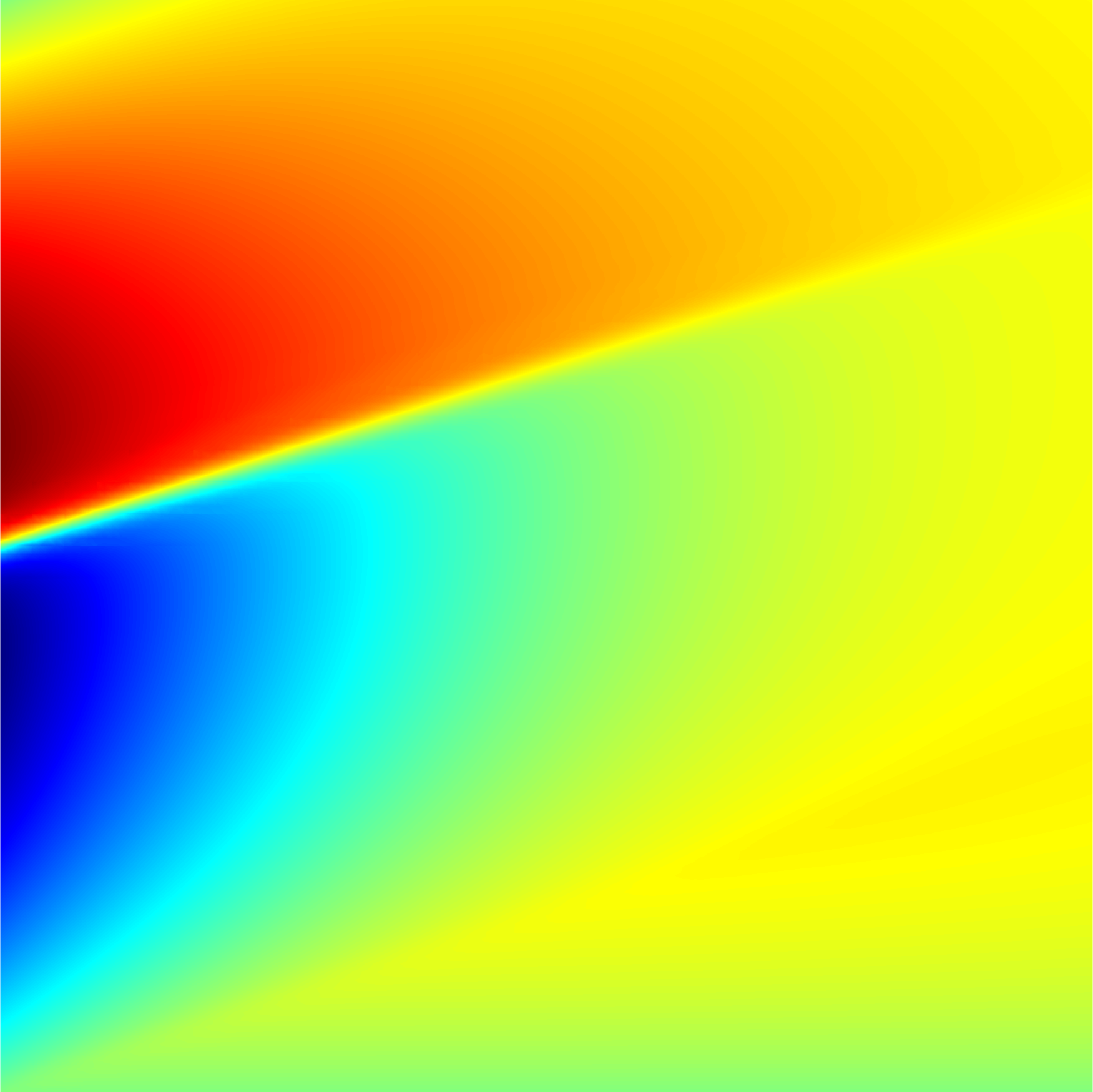};

\nextgroupplot[title={$\mu=(0, 0.55, 80)$}, ylabel={}, ytick=\empty]
\addplot graphics [xmin=0, xmax=1, ymin=0, ymax=1] {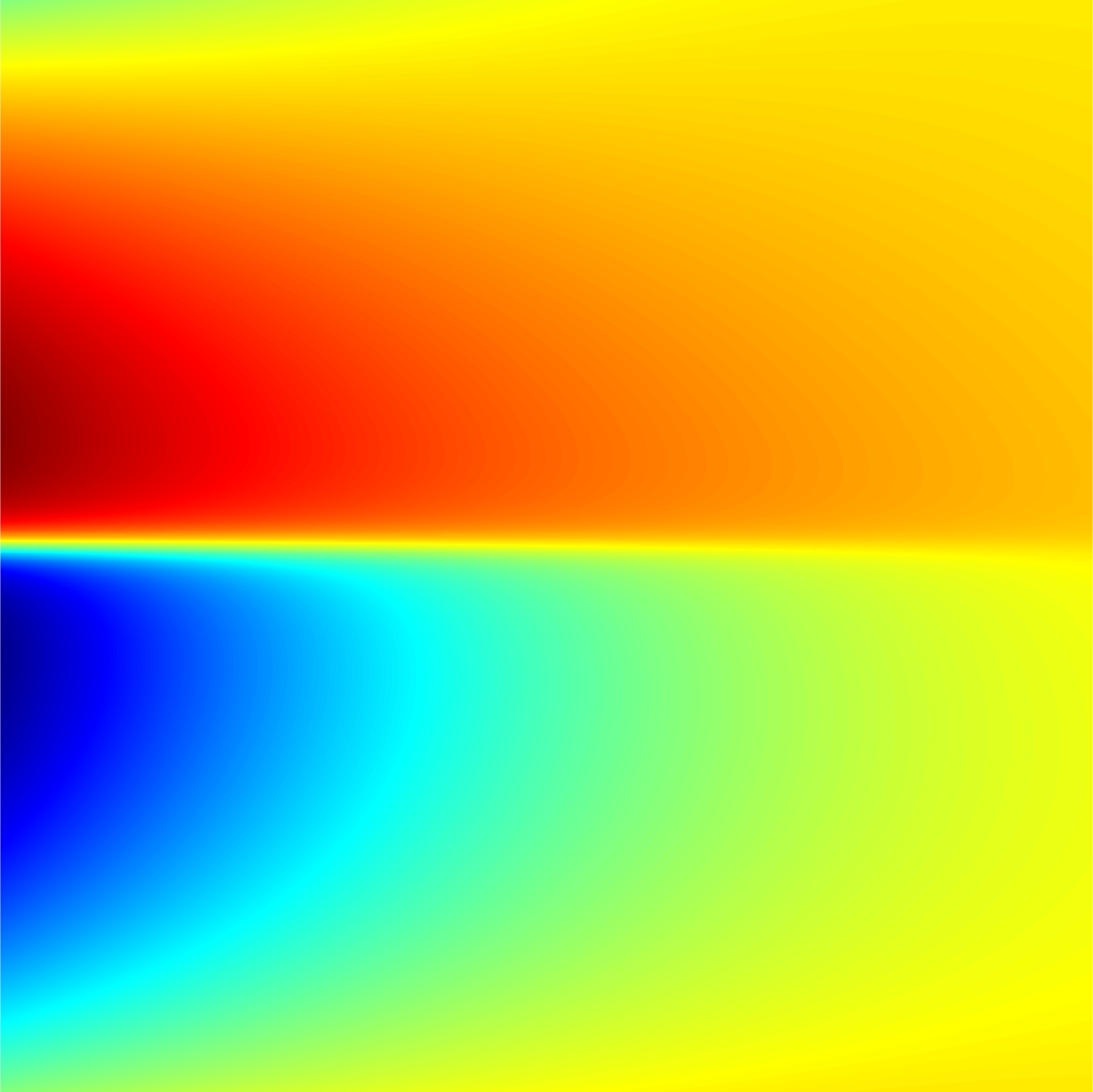};
\end{groupplot}
\end{tikzpicture}
\colorbarMatlabJet{-1.41}{-0.706}{0}{0.706}{1.41}
\caption{Solution of advection-reaction equation (\ref{eqn:advec-react}) for
 three different parameter configurations.}
\label{fig:advec-react}
\end{figure}
We begin by considering the parametrized linear advection-reaction equation in a two-dimensional domain $\Omega \coloneqq [0, 1]^2$ from
\cite{taddei_registration_2020}
\begin{equation} \label{eqn:advec-react}
 \nabla \cdot (\beta u) + \tau u = h\quad\text{in}\quad\Omega, \qquad
 u = \bar{u} \quad\text{on}\quad\partial\Omega_\text{in},
\end{equation}
with solution $u(x;\mu) \in \Rbb$ and parameter vector
$\mu \in \Dcal\coloneqq [-\pi/10,\pi/10]\times[0.3,0.7]\times[60,100]$;
the source term $\func{h}{\Omega}{\Rbb}$, velocity vector
$\func{\beta}{\Dcal}{\Rbb^d}$, reaction coefficient
$\func{\tau}{\Omega\times\Dcal}{\Rbb}$, and boundary
condition $\func{\bar{u}}{\partial\Omega\times\Dcal}{\Rbb}$
are defined as
\begin{equation}
\begin{aligned}
 h(x) &\coloneqq 1+x_1x_2 \\
 \beta((\theta,b,s)) &\coloneqq (\cos\theta,\sin\theta) \\
 \tau(x;(\theta,b,s)) &\coloneqq 1 + b e^{x_1+x_2} \\
 \bar{u}(x;(\theta,b,s)) &\coloneqq
 4\,\text{arctan}\left(s(x_2-1/2)\right)(x_2-x_2^2)
\end{aligned}
\end{equation}
and $\partial\Omega_\text{in}(\mu) \coloneqq
\{x\in\partial\Omega \mid \beta(\mu)\cdot n(x) < 0\}$ is the
parametrized inflow boundary for $\mu\in\Dcal$, where
$\func{n}{\partial\Omega}{\Rbb^d}$ is the unit outward normal.
For this problem, we choose the reference domain to be identical
to the physical domain $\rdom=\pdom=[0,1]^2$ and the nominal
mapping to be the identity mapping $\bar\Gcal=\text{Id}$. The PDE
is discretized using a discontinuous Galerkin method with an upwind
numerical flux on a structured
grid of $2312$ cubic simplex elements ($p=3$) for a total of $N=23120$
degrees of freedom. Because all boundaries of the domain are straight-sided,
we use $q=1$ elements to discretize the domain mapping, which defines
$\Gbb_{h,q}^\text{b}$. Figure~\ref{fig:advec-react} shows the solution
for three parameter configurations; the first parameter ($\theta$) controls
the flow direction and therefore the configuration of the steep gradient,
the primary convection-dominated feature in this problem.

Because the velocity field is constant in space, the primary feature will
be a straight line between $(0, 0.5)$ and $(1, 0.5+\tan\theta)$ for any
$\mu=(\theta,b,s)\in\Dcal$ and there is no need to use a generic domain
parametrization. Instead, we use a parametrization that leverages the
behavior of the solution by defining the mapping $\Gcal$ as
\begin{equation} \label{eqn:onepar}
 (\rcoord, c) \mapsto \Gcal_1(X; c) \coloneqq X_1, \qquad
 (\rcoord, c) \mapsto \Gcal_2(X; c) \coloneqq
 \begin{cases}
  X_2 + 2c X_1 X_2      & \text{if}\quad X_2 < 0.5 \\
  X_2 + 2c X_1(1-X_2) & \text{if}\quad X_2 \geq 0.5,
 \end{cases}
\end{equation}
which leads to a one-dimensional domain parametrization.
This mapping ensures the domain
deformation is  uniformly spaced between the upper (lower)
boundary and the line between $(0, 0.5)$ and $(1, 0.5+c)$. Therefore,
the domain regularization parameters are unnecessary so we take
$\kappa=\lambda=0$, which implies the ROM-IFT method is equivalent
to minimizing the HDM residual over the nonlinear manifold $\Wcal_k$
and the Levenberg-Marquardt solver reduces to a Gauss-Newton solver.
To ensure the optimization problem defining the ROM-IFT method is
solved precisely, we use tight tolerances for the termination criteria,
$\epsilon_1=10^{-15}$ and $\epsilon_2=10^{-10}$, which is usually
achieved in less than $10$ iterations due to the rapid convergence
of the Gauss-Newton method \cite{nocedal_numerical_2006}.
We use this simplified domain deformation to study deep convergence
of the ROM-IFT method; in later section, we use the general domain
parametrization introduced in this work for more complex problems.
Finally, for this problem, we order the training set such that the
centroid of the parameter domain is the first parameter and choose
the domain mapping that initializes the offline training
(Section~\ref{sec:rom:off}) as $\hat\Gcal_1=\mathrm{Id}$.

\subsubsection{Convergence under basis refinement, fixed training set}
\label{sec:numexp:advec:stdy0}
First, we demonstrate rapid, deep convergence of the ROM-IFT method
relative to fixed-domain model reduction. We consider a reduced
parameter space $\tilde\Dcal$ consisting of a one-dimensional
subset of $\Dcal$ where only the first parameter ($\theta$)
varies to focus on variations of the primary convection-dominated
feature
\begin{equation}
 \tilde\Dcal \coloneqq
 \left\{(\theta, 0.55, 80) \mid \theta\in[-\pi/10,\pi/10]\right\}.
\end{equation}
We define a finite subset of $n$ uniformly spaced samples of $\tilde\Dcal$ as
\begin{equation}
 \tilde\Dcal_n \coloneqq 
 \left\{\left(-\frac{\pi}{10} + \left(\frac{i-1}{n-1}\right) \frac{\pi}{5}, 0.55, 80\right) \suchthat i=1,\dots, n
 \right\}.
\end{equation}

\begin{figure}
\centering
\begin{tikzpicture}
\begin{groupplot}[
  group style={
      group size=3 by 2,
      horizontal sep=1cm,
      vertical sep=0.6cm
  },
  width=0.39\textwidth,
  axis equal image,
  xlabel={$x_1$},
  ylabel={$x_2$},
  xtick = {0.0, 0.5, 1.0},
  ytick = {0.0, 0.5, 1.0},
  xmin=0, xmax=1,
  ymin=0, ymax=1
]
\nextgroupplot[xlabel={}, xtick=\empty]
\addplot graphics [xmin=0, xmax=1, ymin=0, ymax=1] {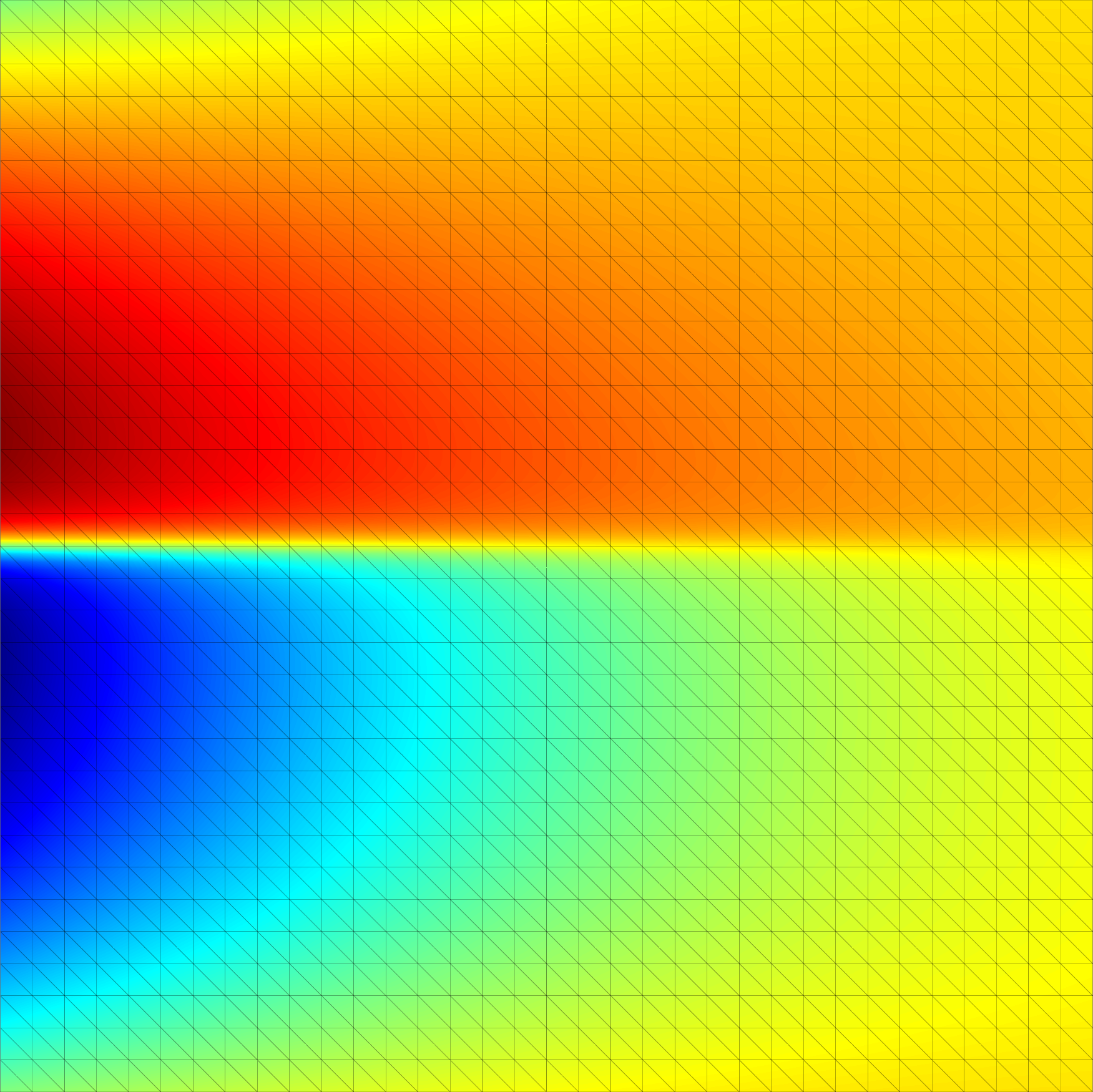};

\nextgroupplot[xlabel={}, xtick=\empty, ylabel={}, ytick=\empty]
\addplot graphics [xmin=0, xmax=1, ymin=0, ymax=1] {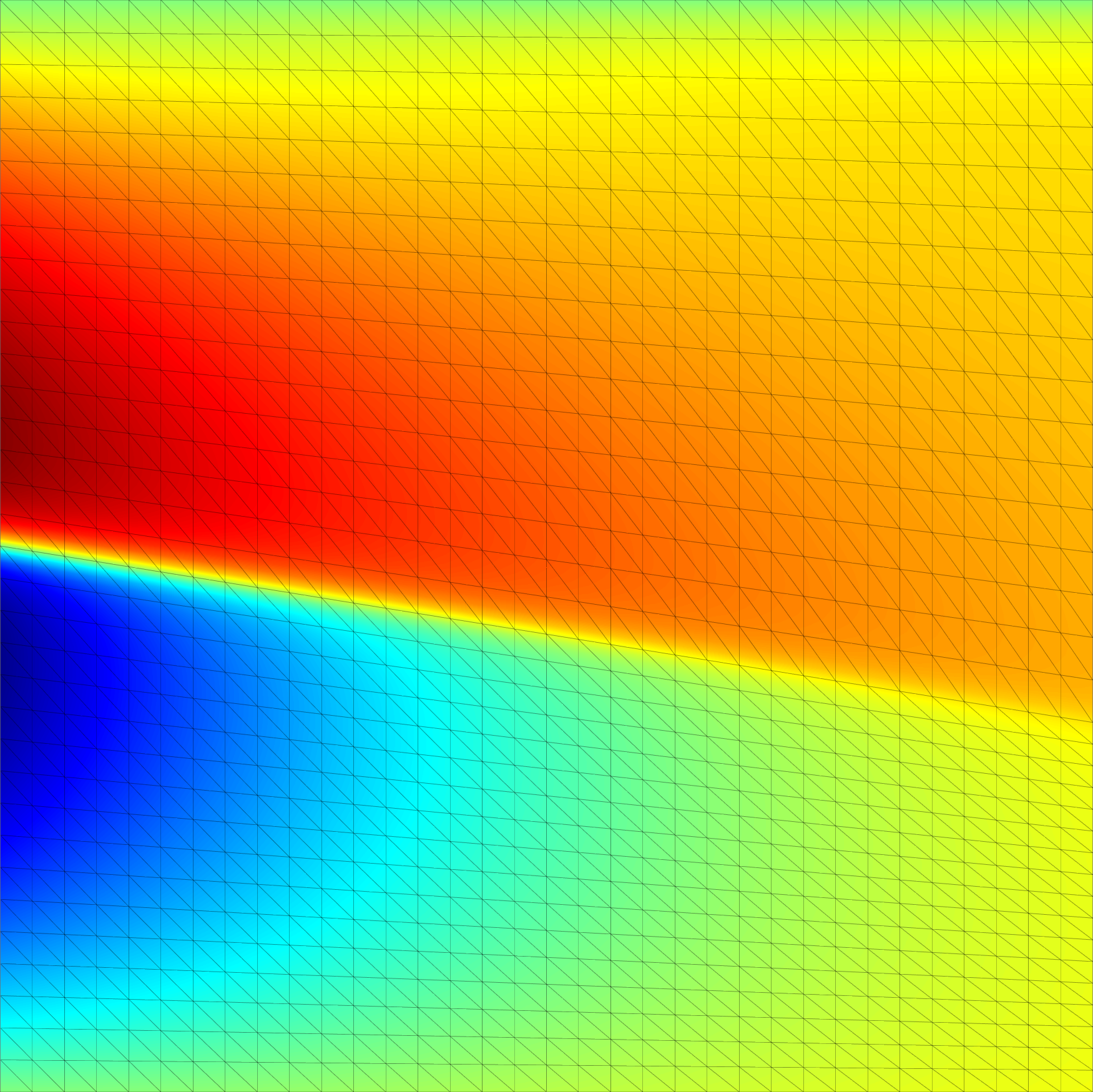};

\nextgroupplot[xlabel={}, xtick=\empty, ylabel={}, ytick=\empty]
\addplot graphics [xmin=0, xmax=1, ymin=0, ymax=1] {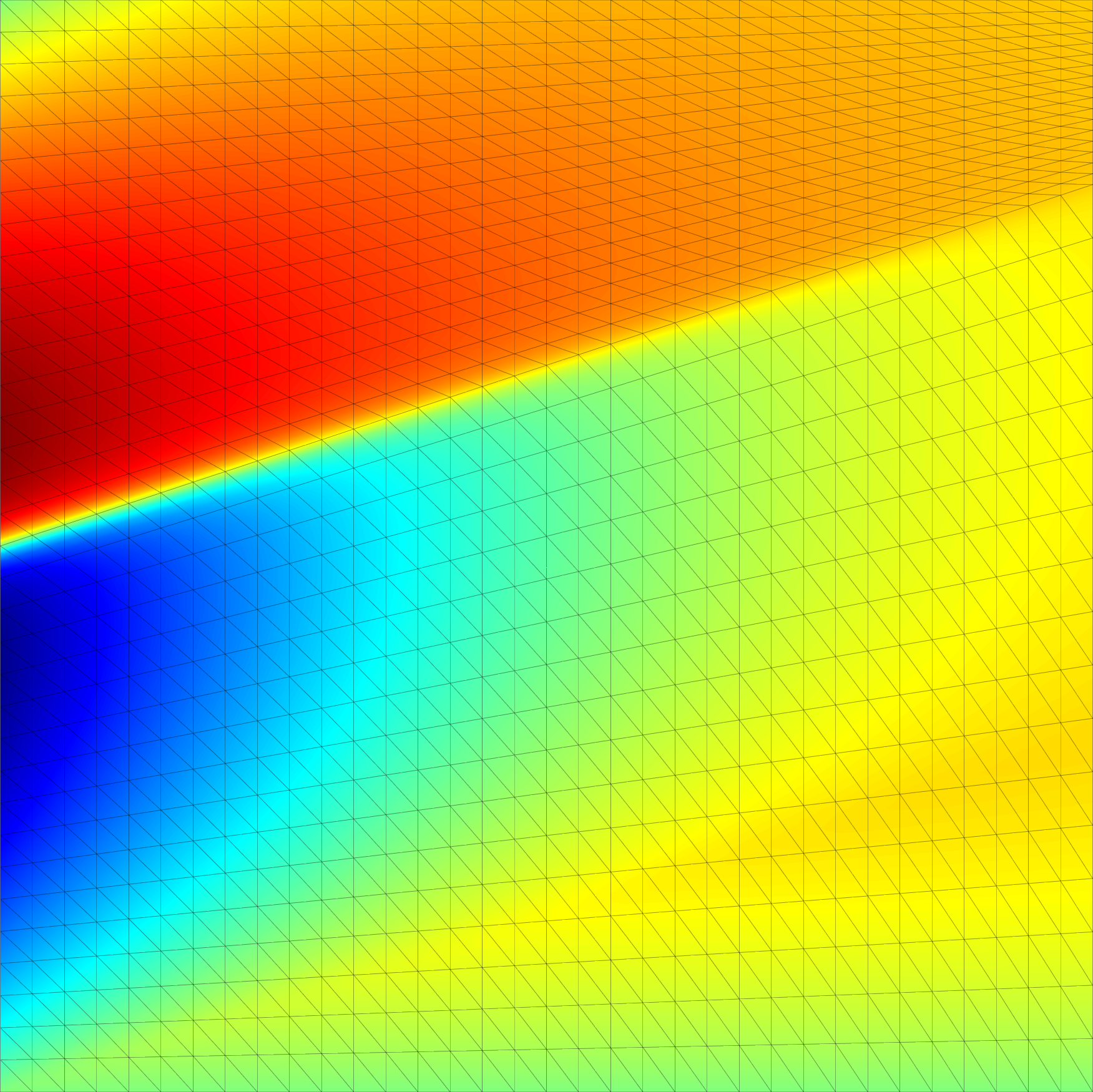};

\nextgroupplot
\addplot graphics [xmin=0, xmax=1, ymin=0, ymax=1] {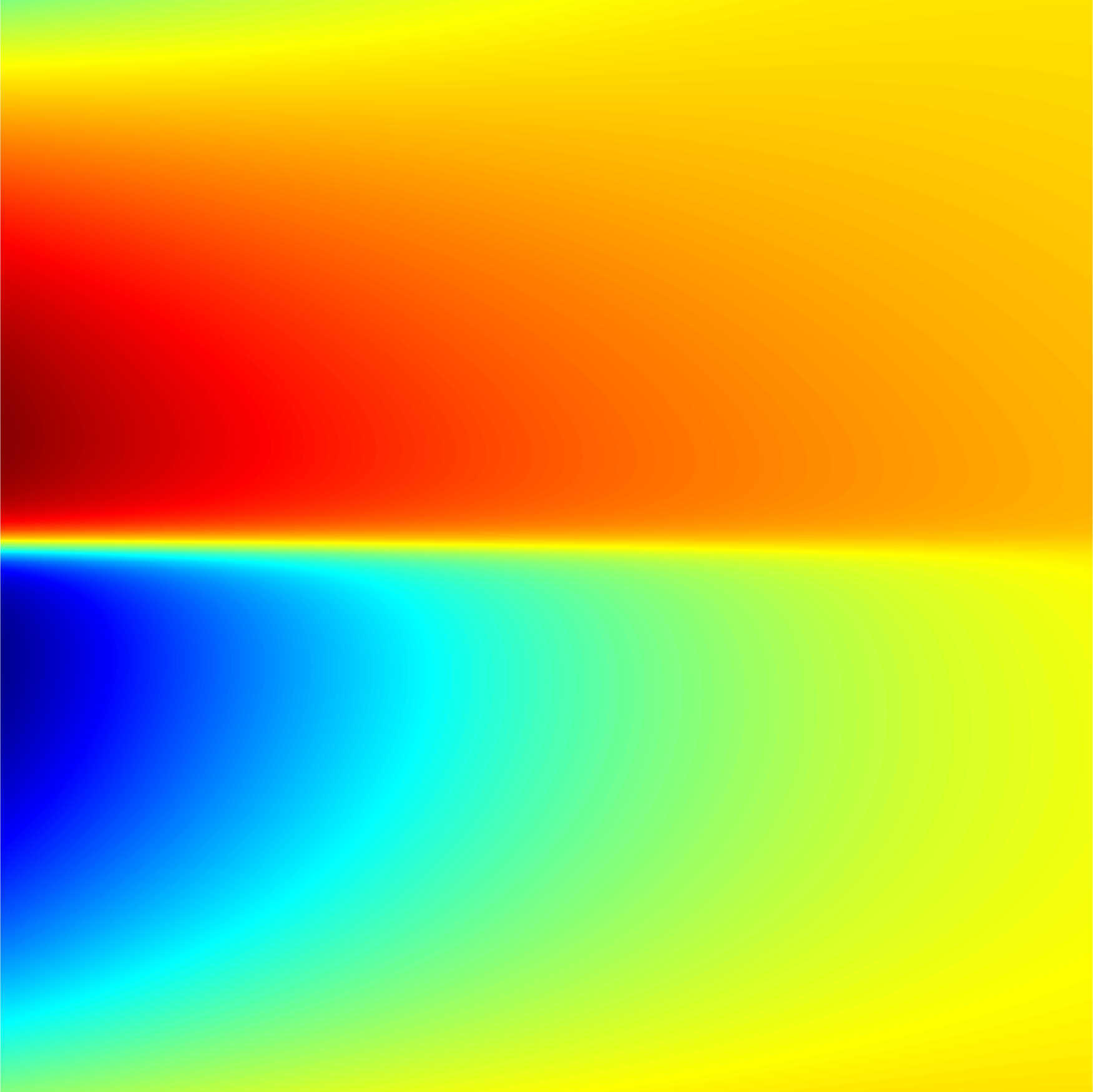};

\nextgroupplot[ylabel={}, ytick=\empty]
\addplot graphics [xmin=0, xmax=1, ymin=0, ymax=1] {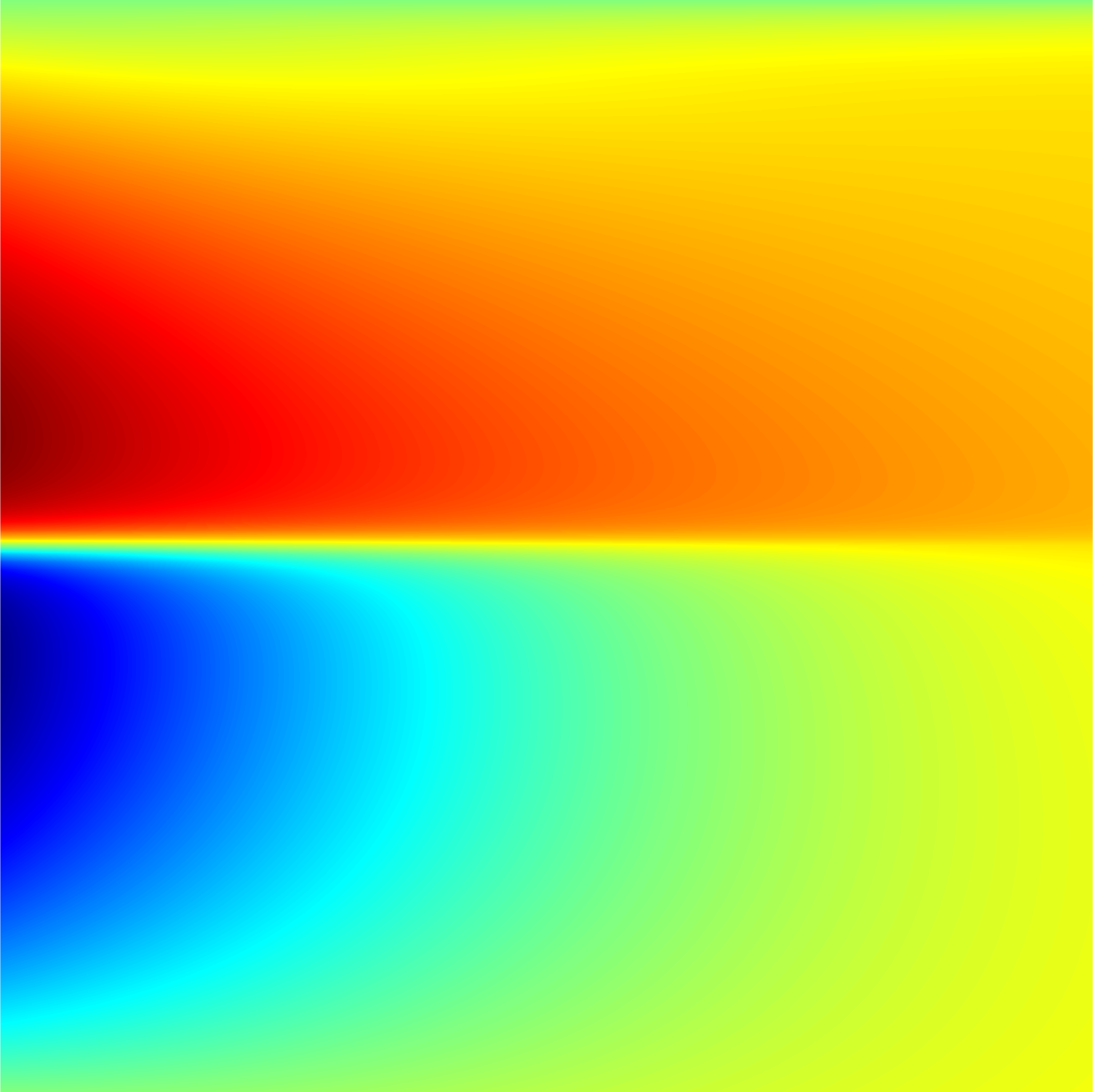};

\nextgroupplot[ylabel={}, ytick=\empty]
\addplot graphics [xmin=0, xmax=1, ymin=0, ymax=1] {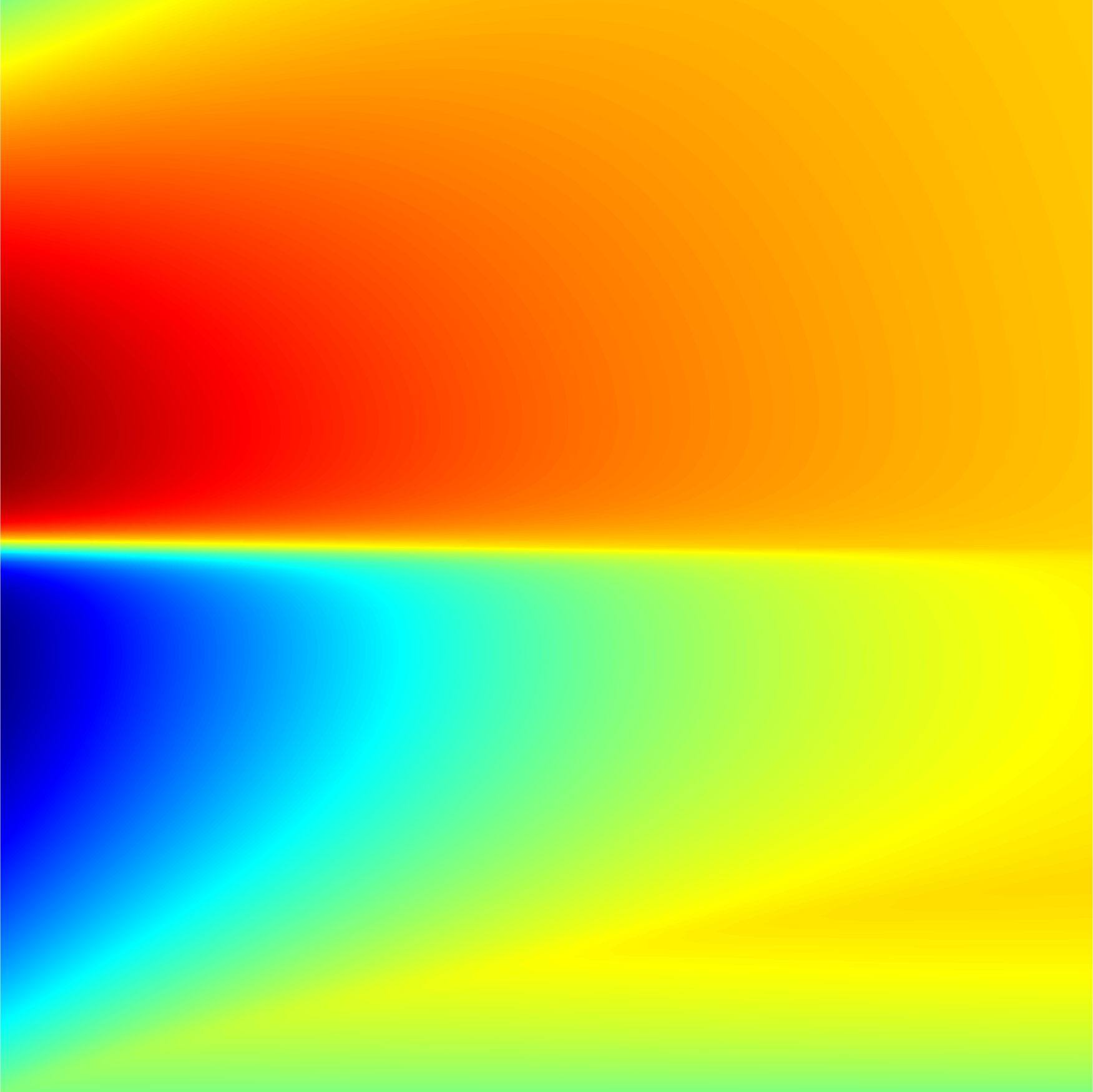};
\end{groupplot}
\end{tikzpicture}
\caption{Selected snapshots generated by the ROM-IFT method corresponding
 to the parameter set $\tilde\Dcal_{101}$. The top row shows the snapshot
 and corresponding domain deformation (mesh edges) in the physical domain,
 while the bottom row shows the corresponding snapshot in the fixed reference
 domain. Colorbar in Figure~\ref{fig:advec-react}.}
\label{fig:advec-react-snap}
\end{figure}

\begin{figure}
\centering
\begin{tikzpicture}
\begin{groupplot}[
  group style={
      group size=3 by 2,
      horizontal sep=1cm,
      vertical sep=0.6cm
  },
  width=0.39\textwidth,
  axis equal image,
  xlabel={$x_1$},
  ylabel={$x_2$},
  xtick = {0.0, 0.5, 1.0},
  ytick = {0.0, 0.5, 1.0},
  xmin=0, xmax=1,
  ymin=0, ymax=1
]
\nextgroupplot[xlabel={}, xtick=\empty]
\addplot graphics [xmin=0, xmax=1, ymin=0, ymax=1] {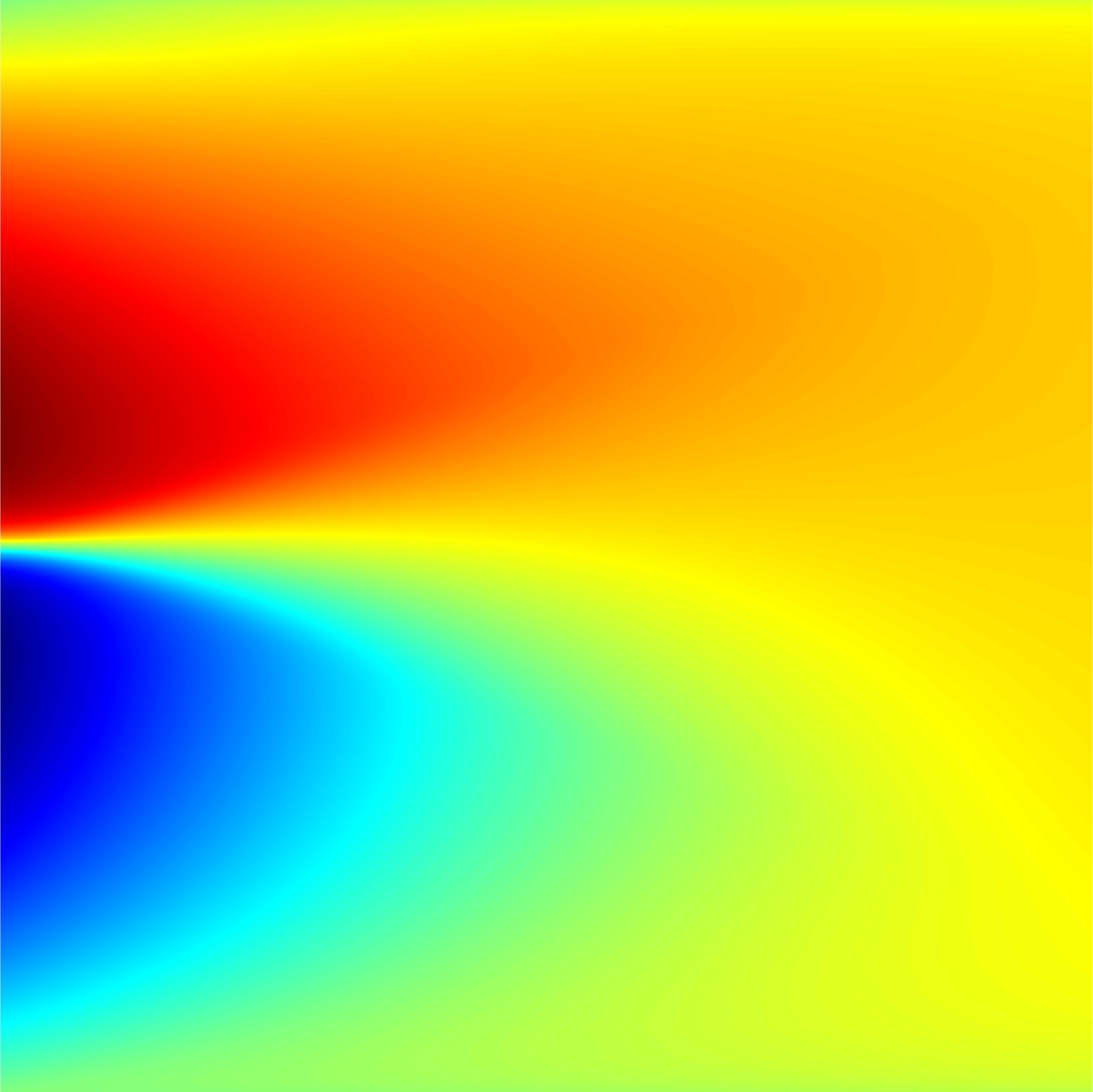};

\nextgroupplot[xlabel={}, xtick=\empty, ylabel={}, ytick=\empty]
\addplot graphics [xmin=0, xmax=1, ymin=0, ymax=1] {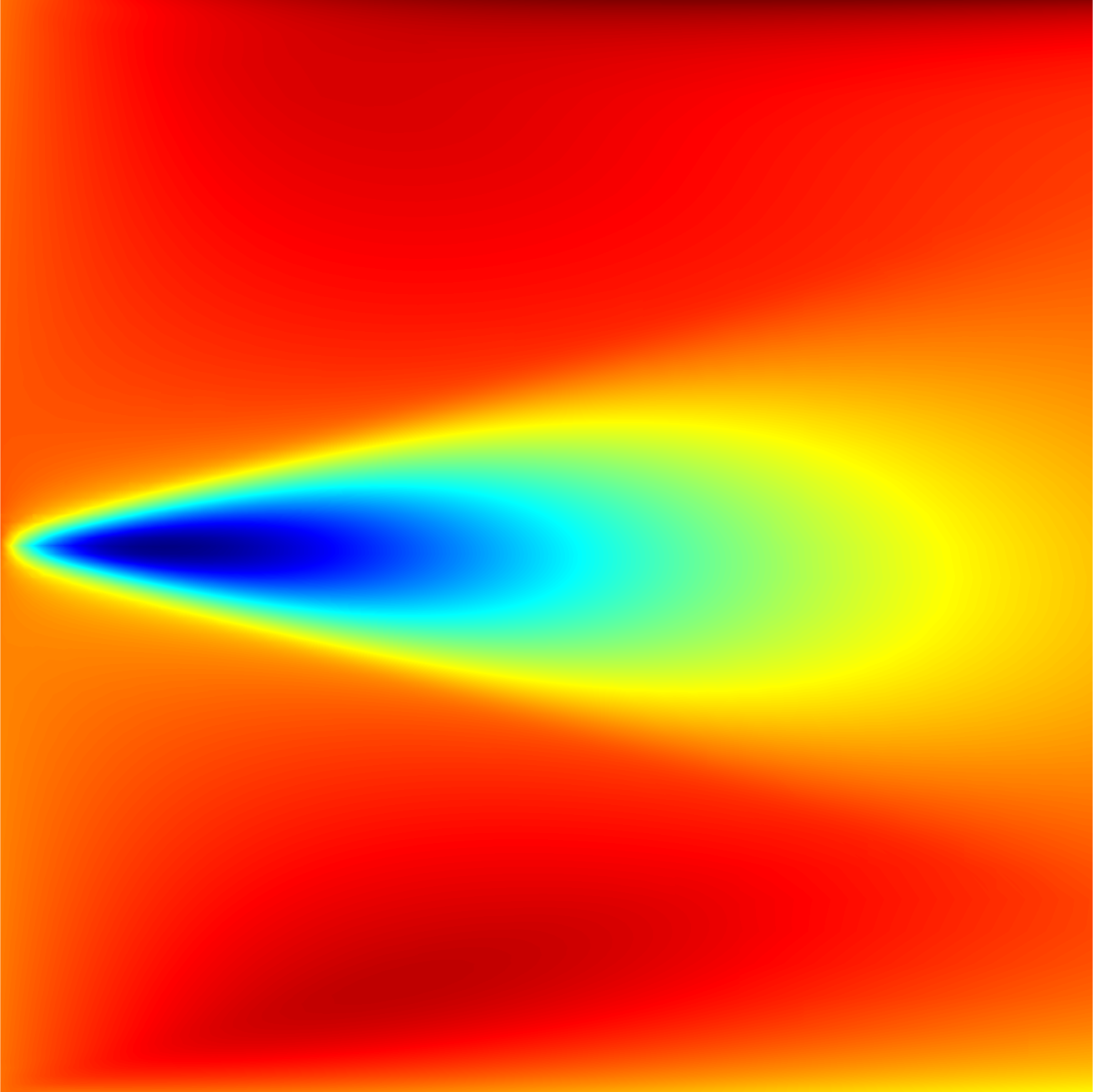};

\nextgroupplot[xlabel={}, xtick=\empty, ylabel={}, ytick=\empty]
\addplot graphics [xmin=0, xmax=1, ymin=0, ymax=1] {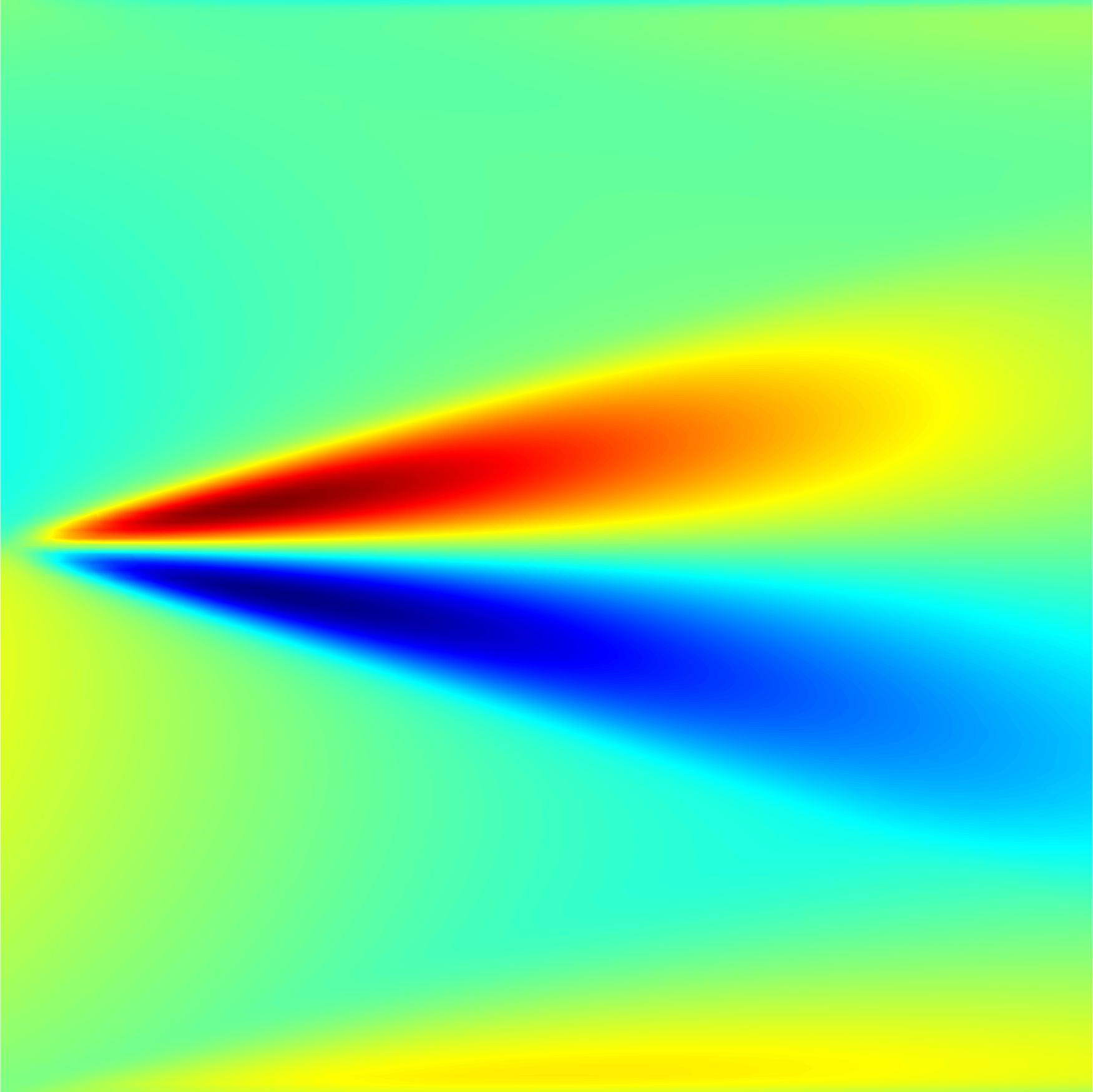};

\nextgroupplot
\addplot graphics [xmin=0, xmax=1, ymin=0, ymax=1] {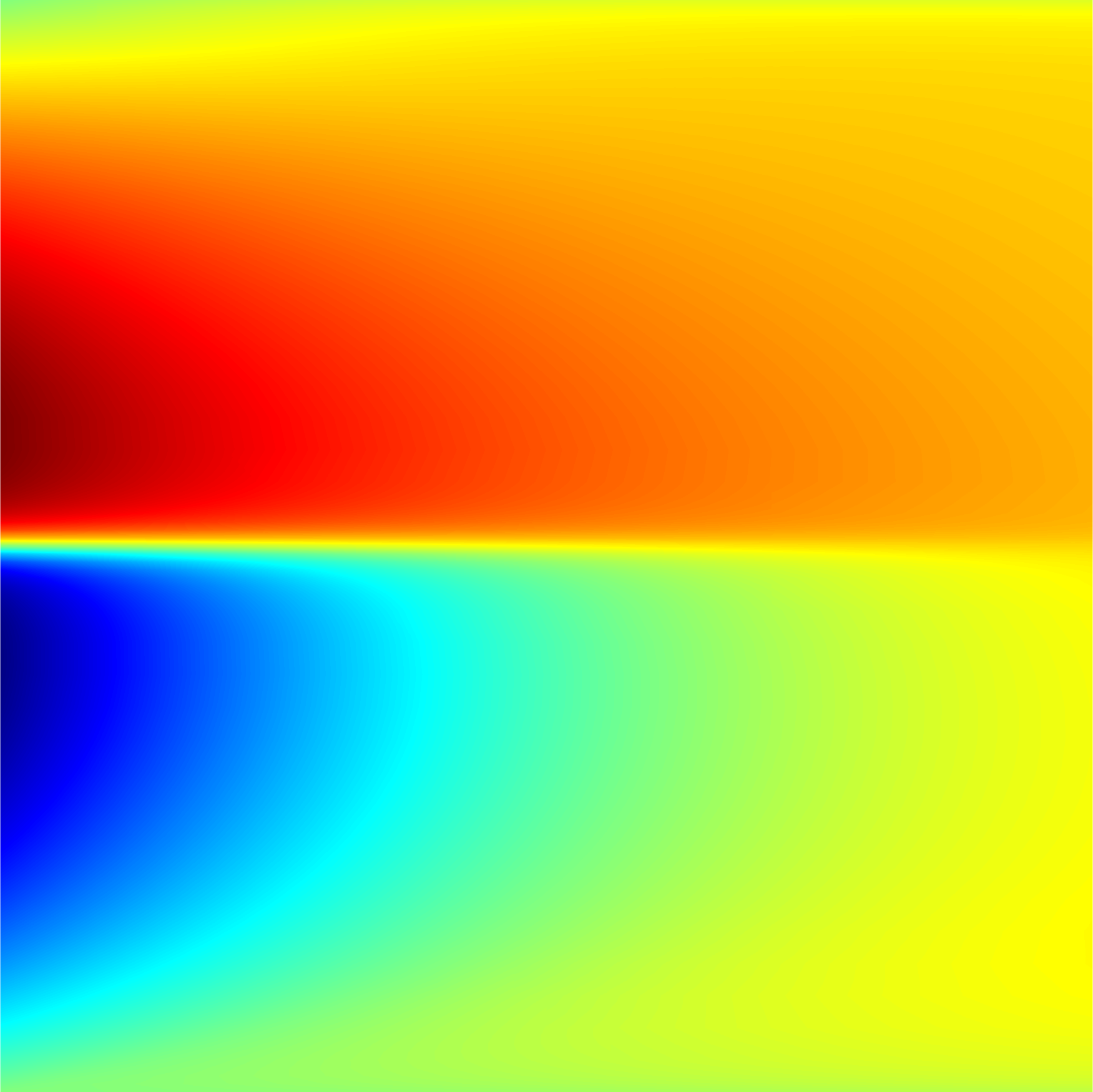};

\nextgroupplot[ylabel={}, ytick=\empty]
\addplot graphics [xmin=0, xmax=1, ymin=0, ymax=1] {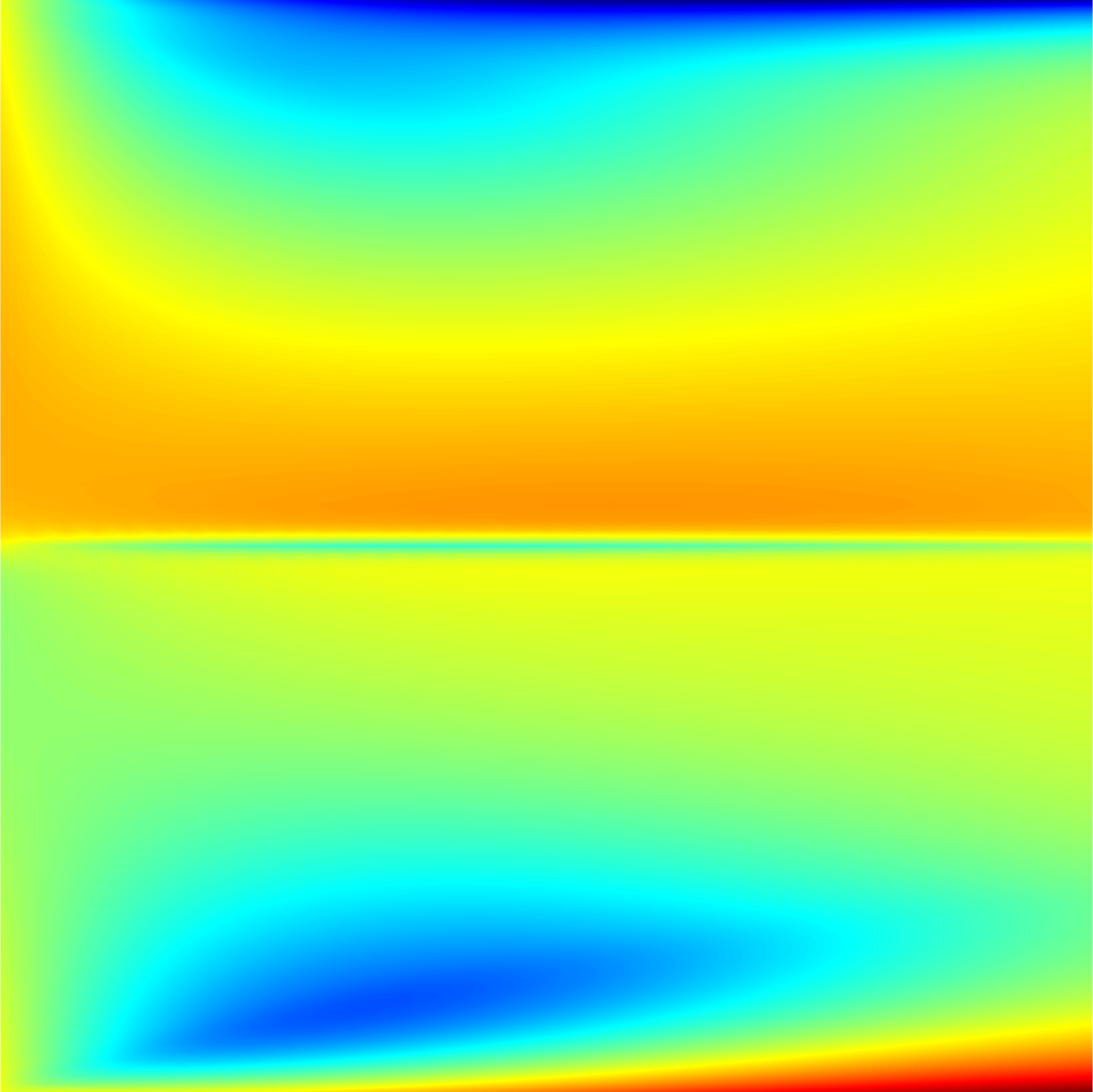};

\nextgroupplot[ylabel={}, ytick=\empty]
\addplot graphics [xmin=0, xmax=1, ymin=0, ymax=1] {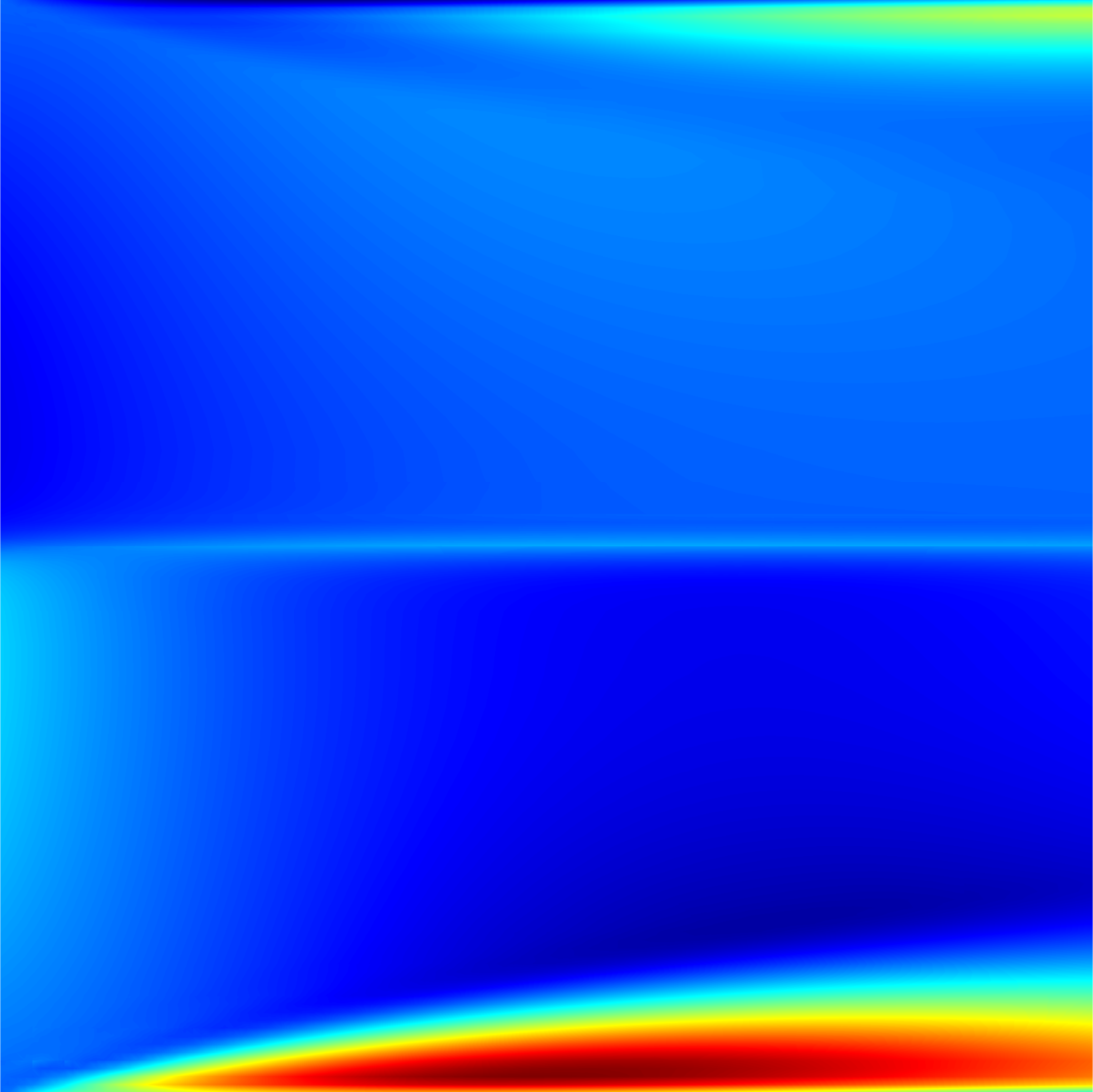};
\end{groupplot}
\end{tikzpicture}
\caption{First three POD modes (\textit{left-to-right}) corresponding to
 non-aligned (\textit{top}) and aligned (\textit{bottom}) snapshots generated
 from the training set $\tilde\Dcal_{101}$. Each figure uses a different
 colorbar, scaled to its range to highlight features in the corresponding
 basis function.}
\label{fig:advec-react-pod}
\end{figure}

To begin, we consider the fixed-domain ROM and ROM-IFT methods
constructed from the training set $\tilde\Dcal_{101}$ with the
corresponding basis truncated to dimension $k$.
The large sample size of $M=101$ is
chosen because the distance from any $\mu\in\tilde\Dcal$ to
the closest point in $\tilde\Dcal_{101}$ is relatively small
(less than $1\%$), so $\tilde\Dcal_{101}$ is a reasonable
approximation to $\tilde\Dcal$. The snapshot alignment of the ROM-IFT
method causes the singular values of the corresponding snapshot matrix
(\ref{eqn:basis}) to decay more rapidly than without alignment
(Figure~\ref{fig:advec-react-study0-cnvg}).
To demonstrate deep convergence, we take the testing set to be
identical to the training set ($\tilde\Dcal_{101}$); the $L^2(\rdom)$
error of both the ROM and ROM-IFT method converge toward machine
zero as $k\rightarrow 101$ with the ROM-IFT method converging more rapidly
(Figure~\ref{fig:advec-react-study0-cnvg}).
For a fixed basis dimension $k$, the ROM-IFT method has a smaller
$L^2(\rdom)$ error for all points in the test set, and for some basis
sizes, it is several orders of magnitude lower
(Figure~\ref{fig:advec-react-study0-cnvg}-\ref{fig:advec-react-study0-sweep}).
In the reference domain, the primary feature is in the same
configuration for all aligned snapshots generated by the
ROM-IFT method (Figure~\ref{fig:advec-react-snap}), effectively
removing the convection-dominated nature of the solution and
enhancing compressibility as demonstrated by the singular value
decay (Figure~\ref{fig:advec-react-study0-cnvg}).
Without alignment, the dominant POD modes of the snapshot matrix
emphasize the parametrized primary feature, which is necessary to
capture its motion; however, with alignment, the first POD mode
captures the primary feature and subsequent modes resolve smoother
regions of the flow (Figure~\ref{fig:advec-react-pod}).

\begin{figure}
\centering
\begin{tikzpicture}
\begin{groupplot} [
group style={group size = 2 by 1, horizontal sep = 3cm}]
\nextgroupplot[width=0.4\textwidth, xlabel=$k$, ymax=10, ylabel={$\sigma_k(\tilde\Dcal_{101})/\sigma_1(\tilde\Dcal_{101})$}, ymode=log, ymin=1e-16]
\addplot [black, thick, mark options={solid, thin}, mark=*, mark size=1.5, mark repeat={5}]
coordinates {
( 0.00000000e+00,  1.00000000e+00)
( 1.00000000e+00,  2.86718571e-01)
( 2.00000000e+00,  1.30194943e-01)
( 3.00000000e+00,  6.96517903e-02)
( 4.00000000e+00,  4.63900339e-02)
( 5.00000000e+00,  3.23529439e-02)
( 6.00000000e+00,  2.35881738e-02)
( 7.00000000e+00,  1.83209565e-02)
( 8.00000000e+00,  1.40508663e-02)
( 9.00000000e+00,  1.15461568e-02)
( 1.00000000e+01,  9.12622065e-03)
( 1.10000000e+01,  7.70984217e-03)
( 1.20000000e+01,  6.30677150e-03)
( 1.30000000e+01,  5.39516979e-03)
( 1.40000000e+01,  4.58319830e-03)
( 1.50000000e+01,  3.88111301e-03)
( 1.60000000e+01,  3.49081633e-03)
( 1.70000000e+01,  2.88117347e-03)
( 1.80000000e+01,  2.65175234e-03)
( 1.90000000e+01,  2.16712851e-03)
( 2.00000000e+01,  1.98541576e-03)
( 2.10000000e+01,  1.64184083e-03)
( 2.20000000e+01,  1.47890882e-03)
( 2.30000000e+01,  1.24665519e-03)
( 2.40000000e+01,  1.11144447e-03)
( 2.50000000e+01,  9.63931803e-04)
( 2.60000000e+01,  8.30591334e-04)
( 2.70000000e+01,  7.41781684e-04)
( 2.80000000e+01,  6.30236616e-04)
( 2.90000000e+01,  5.69011761e-04)
( 3.00000000e+01,  4.81481468e-04)
( 3.10000000e+01,  4.37327158e-04)
( 3.20000000e+01,  3.67324823e-04)
( 3.30000000e+01,  3.35766395e-04)
( 3.40000000e+01,  2.75677484e-04)
( 3.50000000e+01,  2.57315070e-04)
( 3.60000000e+01,  2.01969990e-04)
( 3.70000000e+01,  1.95237778e-04)
( 3.80000000e+01,  1.48966246e-04)
( 3.90000000e+01,  1.44649124e-04)
( 4.00000000e+01,  1.12090395e-04)
( 4.10000000e+01,  1.04837936e-04)
( 4.20000000e+01,  8.31542343e-05)
( 4.30000000e+01,  7.61853139e-05)
( 4.40000000e+01,  5.89289285e-05)
( 4.50000000e+01,  5.35443628e-05)
( 4.60000000e+01,  4.09398105e-05)
( 4.70000000e+01,  3.57276880e-05)
( 4.80000000e+01,  3.00382415e-05)
( 4.90000000e+01,  2.34203180e-05)
( 5.00000000e+01,  1.90881504e-05)
( 5.10000000e+01,  1.53708603e-05)
( 5.20000000e+01,  1.27995549e-05)
( 5.30000000e+01,  9.49166456e-06)
( 5.40000000e+01,  7.74438298e-06)
( 5.50000000e+01,  5.71564795e-06)
( 5.60000000e+01,  4.57035334e-06)
( 5.70000000e+01,  3.38464636e-06)
( 5.80000000e+01,  2.65806664e-06)
( 5.90000000e+01,  1.86156812e-06)
( 6.00000000e+01,  1.46009712e-06)
( 6.10000000e+01,  1.05364918e-06)
( 6.20000000e+01,  7.59536191e-07)
( 6.30000000e+01,  5.46083211e-07)
( 6.40000000e+01,  3.73737522e-07)
( 6.50000000e+01,  2.88987824e-07)
( 6.60000000e+01,  1.75277278e-07)
( 6.70000000e+01,  1.38283367e-07)
( 6.80000000e+01,  7.85966018e-08)
( 6.90000000e+01,  6.80944912e-08)
( 7.00000000e+01,  3.32316616e-08)
( 7.10000000e+01,  2.99434105e-08)
( 7.20000000e+01,  1.37755304e-08)
( 7.30000000e+01,  1.27563334e-08)
( 7.40000000e+01,  5.61601836e-09)
( 7.50000000e+01,  4.58824548e-09)
( 7.60000000e+01,  2.27860185e-09)
( 7.70000000e+01,  1.57710322e-09)
( 7.80000000e+01,  8.61894350e-10)
( 7.90000000e+01,  5.23896516e-10)
( 8.00000000e+01,  2.99805556e-10)
( 8.10000000e+01,  1.55005588e-10)
( 8.20000000e+01,  1.01149683e-10)
( 8.30000000e+01,  3.91931049e-11)
( 8.40000000e+01,  3.08252482e-11)
( 8.50000000e+01,  9.01961372e-12)
( 8.60000000e+01,  8.49194136e-12)
( 8.70000000e+01,  2.34552115e-12)
( 8.80000000e+01,  2.11455800e-12)
( 8.90000000e+01,  5.64850701e-13)
( 9.00000000e+01,  4.48239530e-13)
( 9.10000000e+01,  1.35974800e-13)
( 9.20000000e+01,  7.17980091e-14)
( 9.30000000e+01,  2.69036015e-14)
( 9.40000000e+01,  1.08174308e-14)
( 9.50000000e+01,  4.61142716e-15)
( 9.60000000e+01,  1.47925237e-15)
( 9.70000000e+01,  6.86183617e-16)
( 9.80000000e+01,  1.77641999e-16)
( 9.90000000e+01,  1.50773290e-16)
( 1.00000000e+02,  1.34564666e-16)};\label{line:advec1v2d_taddei0_study0_svals_romfix}

\addplot [blue, thick, mark options={solid, thin}, mark=square*, mark size=1.5, mark repeat={5}]
coordinates {
( 0.00000000e+00,  1.00000000e+00)
( 1.00000000e+00,  9.84020220e-02)
( 2.00000000e+00,  2.76045452e-02)
( 3.00000000e+00,  1.63943822e-02)
( 4.00000000e+00,  8.75969863e-03)
( 5.00000000e+00,  5.35727326e-03)
( 6.00000000e+00,  4.09179681e-03)
( 7.00000000e+00,  3.27664915e-03)
( 8.00000000e+00,  2.22220962e-03)
( 9.00000000e+00,  1.45720985e-03)
( 1.00000000e+01,  1.41692057e-03)
( 1.10000000e+01,  9.99889425e-04)
( 1.20000000e+01,  7.14938891e-04)
( 1.30000000e+01,  6.24428447e-04)
( 1.40000000e+01,  5.30479712e-04)
( 1.50000000e+01,  4.06714765e-04)
( 1.60000000e+01,  3.75796784e-04)
( 1.70000000e+01,  2.96152596e-04)
( 1.80000000e+01,  2.24637628e-04)
( 1.90000000e+01,  2.14789178e-04)
( 2.00000000e+01,  1.53066742e-04)
( 2.10000000e+01,  1.41700449e-04)
( 2.20000000e+01,  1.08462992e-04)
( 2.30000000e+01,  9.31406082e-05)
( 2.40000000e+01,  7.80801969e-05)
( 2.50000000e+01,  6.69075189e-05)
( 2.60000000e+01,  5.00985553e-05)
( 2.70000000e+01,  4.96436903e-05)
( 2.80000000e+01,  3.60416853e-05)
( 2.90000000e+01,  3.23103805e-05)
( 3.00000000e+01,  2.51962430e-05)
( 3.10000000e+01,  1.87616943e-05)
( 3.20000000e+01,  1.71246157e-05)
( 3.30000000e+01,  1.18693379e-05)
( 3.40000000e+01,  1.08225161e-05)
( 3.50000000e+01,  8.30952353e-06)
( 3.60000000e+01,  5.76831797e-06)
( 3.70000000e+01,  5.52714159e-06)
( 3.80000000e+01,  3.53452924e-06)
( 3.90000000e+01,  2.85823626e-06)
( 4.00000000e+01,  2.25656029e-06)
( 4.10000000e+01,  1.36145586e-06)
( 4.20000000e+01,  1.31436895e-06)
( 4.30000000e+01,  7.85156591e-07)
( 4.40000000e+01,  5.61937222e-07)
( 4.50000000e+01,  4.10883705e-07)
( 4.60000000e+01,  2.28678438e-07)
( 4.70000000e+01,  2.00053272e-07)
( 4.80000000e+01,  9.10628542e-08)
( 4.90000000e+01,  8.88340456e-08)
( 5.00000000e+01,  3.83900487e-08)
( 5.10000000e+01,  3.22776970e-08)
( 5.20000000e+01,  1.38171625e-08)
( 5.30000000e+01,  1.07828096e-08)
( 5.40000000e+01,  4.29326999e-09)
( 5.50000000e+01,  3.17512327e-09)
( 5.60000000e+01,  1.50090636e-09)
( 5.70000000e+01,  8.87342502e-10)
( 5.80000000e+01,  5.37344579e-10)
( 5.90000000e+01,  2.84444510e-10)
( 6.00000000e+01,  2.57293004e-10)
( 6.10000000e+01,  1.26312943e-10)
( 6.20000000e+01,  1.11657303e-10)
( 6.30000000e+01,  5.35409053e-11)
( 6.40000000e+01,  5.15045710e-11)
( 6.50000000e+01,  2.41368735e-11)
( 6.60000000e+01,  2.17934400e-11)
( 6.70000000e+01,  1.13550416e-11)
( 6.80000000e+01,  9.52731202e-12)
( 6.90000000e+01,  5.55911825e-12)
( 7.00000000e+01,  4.55635786e-12)
( 7.10000000e+01,  4.10647961e-12)
( 7.20000000e+01,  1.77540972e-12)
( 7.30000000e+01,  1.70069468e-12)
( 7.40000000e+01,  8.04331925e-13)
( 7.50000000e+01,  6.03961451e-13)
( 7.60000000e+01,  3.50163020e-13)
( 7.70000000e+01,  2.34101434e-13)
( 7.80000000e+01,  1.18822256e-13)
( 7.90000000e+01,  9.95851812e-14)
( 8.00000000e+01,  6.51825346e-14)
( 8.10000000e+01,  3.81705786e-14)
( 8.20000000e+01,  3.63796659e-14)
( 8.30000000e+01,  2.60420067e-14)
( 8.40000000e+01,  1.35479361e-14)
( 8.50000000e+01,  1.25427252e-14)
( 8.60000000e+01,  5.52185453e-15)
( 8.70000000e+01,  3.54318229e-15)
( 8.80000000e+01,  3.11361466e-15)
( 8.90000000e+01,  1.25100846e-15)
( 9.00000000e+01,  8.96982328e-16)
( 9.10000000e+01,  8.75371918e-16)
( 9.20000000e+01,  6.19201715e-16)
( 9.30000000e+01,  4.96279288e-16)
( 9.40000000e+01,  4.59119692e-16)
( 9.50000000e+01,  4.42328300e-16)
( 9.60000000e+01,  4.26707200e-16)
( 9.70000000e+01,  4.13660944e-16)
( 9.80000000e+01,  4.07259606e-16)
( 9.90000000e+01,  4.04141903e-16)
( 1.00000000e+02,  3.35970531e-16)};\label{line:advec1v2d_taddei0_study0_svals_romtrk}

\nextgroupplot[width=0.4\textwidth, xlabel=$k$, ymax=10, ylabel={$E_k(\tilde\Dcal_{101})$}, ymode=log, ymin=1e-12]
\addplot [black, thick, mark options={solid, thin}, mark=*, mark size=1.5, mark repeat={1}]
coordinates {
( 1.00000000e+00,  4.81234229e-01)
( 2.00000000e+00,  2.91336991e-01)
( 3.00000000e+00,  1.94590646e-01)
( 4.00000000e+00,  1.36787008e-01)
( 5.00000000e+00,  1.06738694e-01)
( 1.00000000e+01,  3.64233711e-02)
( 2.00000000e+01,  8.66403148e-03)
( 3.00000000e+01,  1.86277637e-03)
( 4.00000000e+01,  5.01037838e-04)
( 5.00000000e+01,  8.01892681e-05)
( 6.00000000e+01,  5.75724568e-06)
( 7.00000000e+01,  1.20488803e-07)};\label{line:advec1v2d_taddei0_study0_L2err_romfix}

\addplot [blue, thick, mark options={solid, thin}, mark=square*, mark size=1.5, mark repeat={1}]
coordinates {
( 1.00000000e+00,  1.67467931e-01)
( 2.00000000e+00,  6.12701981e-02)
( 3.00000000e+00,  4.74232092e-02)
( 4.00000000e+00,  2.53043599e-02)
( 5.00000000e+00,  2.03931006e-02)
( 1.00000000e+01,  6.18861635e-03)
( 2.00000000e+01,  7.26840457e-04)
( 3.00000000e+01,  9.72679224e-05)
( 4.00000000e+01,  6.04093012e-06)
( 5.00000000e+01,  1.34115083e-07)
( 6.00000000e+01,  1.06457274e-09)
( 7.00000000e+01,  5.78661035e-11)};\label{line:advec1v2d_taddei0_study0_L2err_romfix}

\end{groupplot}\end{tikzpicture}
\caption{\textit{Left}: Convergence of the singular values of the
 non-aligned (\ref{line:advec1v2d_taddei0_study0_svals_romfix}) and
 aligned (\ref{line:advec1v2d_taddei0_study0_svals_romtrk}) snapshot
 matrices associated with the training set $\tilde\Dcal_{101}$.
 \textit{Right}: Convergence of the maximum relative $L^2(\rdom)$
 error over the test set $\tilde\Dcal_{101}$ for the fixed-domain ROM
 (\ref{line:advec1v2d_taddei0_study0_svals_romfix}) and
 ROM-IFT (\ref{line:advec1v2d_taddei0_study0_svals_romtrk})
 with the size of the reduced basis ($k$) (training set:
 $\tilde\Dcal_{101}$).}
\label{fig:advec-react-study0-cnvg}
\end{figure}

\begin{figure}
\centering
\input{_py/advec1v2d_taddei0_study0_sweep.tikz}
\caption{The relative $L^2(\rdom)$ error across all parameters in the
 test set $\tilde\Dcal_{101}$ for the fixed-domain ROM with
 basis size
$k=1$ (\ref{line:advec1v2d_taddei0_study0_L2err_romfix_sweep_nrob1}),
$k=10$ (\ref{line:advec1v2d_taddei0_study0_L2err_romfix_sweep_nrob10}),
$k=30$ (\ref{line:advec1v2d_taddei0_study0_L2err_romfix_sweep_nrob30})
and the ROM-IFT with basis size
$k=1$ (\ref{line:advec1v2d_taddei0_study0_L2err_romtrk_sweep_nrob1}),
$k=10$ (\ref{line:advec1v2d_taddei0_study0_L2err_romtrk_sweep_nrob10}),
$k=30$ (\ref{line:advec1v2d_taddei0_study0_L2err_romtrk_sweep_nrob30}).}
\label{fig:advec-react-study0-sweep}
\end{figure}

\subsubsection{Convergence under training set refinement}
\label{sec:numexp:advec:stdy1}
Next, we consider the relative performance of fixed-domain
model reduction and the ROM-IFT method as the training set
is refined. 
In particular, we take the training set to be
$\tilde\Dcal_n$ for $1\leq n\leq 41$ and the testing set to
be $\tilde\Dcal_{101}$. In all cases, we do not truncate the
POD modes, which leads to a basis of dimension $n$. The
smallest singular value of the aligned snapshot matrix
associated with the parameter set $\tilde\Dcal_n$ decays
significantly more rapidly than the non-aligned snapshot
matrix (Figure~\ref{fig:advec-react-study1-cnvg}),
suggesting that amount of new information in
additional snapshots is rapidly diminishing when
the snapshots are aligned as described in
Section~\ref{sec:rom:off}.
In addition, we observe
that the ROM-IFT method produces solutions with
significantly lower error than the fixed-domain
ROM method for small training sets; the improvement
diminishes as the training set is enriched
(Figure~\ref{fig:advec-react-study1-cnvg}).
For example, for $n=3$ ($n=5$), the maximum
relative $L^2(\rdom)$ error over the testing set is
$20.8\%$ ($13.8\%$) for the fixed-domain ROM and only
$6.0\%$ ($2.4\%$) for the ROM-IFT. Whereas for $n=19$,
the error is $1.1\%$ for the fixed-domain ROM and
$0.53\%$ for the ROM-IFT. This suggests that the
largest benefit from the ROM-IFT method comes when
the training is limited.
\begin{figure}
\centering
\begin{tikzpicture}
\begin{groupplot} [
group style={group size = 2 by 3, horizontal sep = 3cm, vertical sep = 1.5cm}]
\nextgroupplot[width=0.4\textwidth, xlabel=$n$, ymax=1.1, ylabel={$\sigma_n(\tilde\Dcal_n)/\sigma_1(\tilde\Dcal_n)$}, ymode=log, ymin=1e-08]
\addplot [black, thick, mark options={solid, thin}, mark=*, mark size=1.5]
coordinates {
( 3.00000000e+00,  1.96229542e-01)
( 5.00000000e+00,  4.80154817e-02)
( 7.00000000e+00,  2.27763317e-02)
( 9.00000000e+00,  1.39768912e-02)
( 1.10000000e+01,  1.00266413e-02)
( 1.30000000e+01,  7.90456840e-03)
( 1.50000000e+01,  6.28463366e-03)
( 1.70000000e+01,  4.48261098e-03)
( 1.90000000e+01,  3.17286342e-03)
( 2.10000000e+01,  2.26113459e-03)
( 3.10000000e+01,  3.82682232e-04)
( 4.10000000e+01,  3.78037516e-05)};\label{line:advec1v2d_taddei0_study1_svals_romfix}

\addplot [blue, thick, mark options={solid, thin}, mark=square*, mark size=1.5]
coordinates {
( 3.00000000e+00,  4.82973415e-02)
( 5.00000000e+00,  1.28223469e-02)
( 7.00000000e+00,  3.93017512e-03)
( 9.00000000e+00,  9.57800210e-04)
( 1.10000000e+01,  2.90819564e-04)
( 1.30000000e+01,  1.24336364e-04)
( 1.50000000e+01,  7.06037531e-05)
( 1.70000000e+01,  4.36243324e-05)
( 1.90000000e+01,  2.75971748e-05)
( 2.10000000e+01,  1.75176053e-05)
( 3.10000000e+01,  1.47897653e-06)
( 4.10000000e+01,  3.00981263e-08)};\label{line:advec1v2d_taddei0_study1_svals_romtrk}

\nextgroupplot[width=0.4\textwidth, ytick={1e-5, 1e-4, 1e-3, 1e-2, 1e-1, 1}, xlabel=$n$, ymax=0.3, ylabel={$E_n(\tilde\Dcal_{101})$}, ymode=log, ymin=0.0008]
\addplot [black, thick, mark options={solid, thin}, mark=*, mark size=1.5]
coordinates {
( 3.00000000e+00,  2.07625191e-01)
( 5.00000000e+00,  1.37616212e-01)
( 7.00000000e+00,  7.75930147e-02)
( 9.00000000e+00,  5.03061428e-02)
( 1.10000000e+01,  3.49747393e-02)
( 1.30000000e+01,  2.53446751e-02)
( 1.50000000e+01,  1.86224972e-02)
( 1.70000000e+01,  1.45002649e-02)
( 1.90000000e+01,  1.10462446e-02)
( 2.10000000e+01,  8.70535507e-03)
( 3.10000000e+01,  2.87047881e-03)
( 4.10000000e+01,  1.30601747e-03)};\label{line:advec1v2d_taddei0_study1_L2err_romfix}

\addplot [blue, thick, mark options={solid, thin}, mark=square*, mark size=1.5]
coordinates {
( 3.00000000e+00,  5.96855075e-02)
( 5.00000000e+00,  2.39786321e-02)
( 7.00000000e+00,  1.77488593e-02)
( 9.00000000e+00,  1.42531372e-02)
( 1.10000000e+01,  1.21854441e-02)
( 1.30000000e+01,  1.01051532e-02)
( 1.50000000e+01,  8.16069128e-03)
( 1.70000000e+01,  6.43820223e-03)
( 1.90000000e+01,  5.26645278e-03)
( 2.10000000e+01,  4.65738247e-03)
( 3.10000000e+01,  2.12248805e-03)
( 4.10000000e+01,  9.61396199e-04)};\label{line:advec1v2d_taddei0_study1_L2err_romtrk}

\end{groupplot}\end{tikzpicture}
\caption{\textit{Left}: Convergence of the smallest singular value of the
 non-aligned (\ref{line:advec1v2d_taddei0_study1_svals_romfix}) and
 aligned (\ref{line:advec1v2d_taddei0_study1_svals_romtrk}) snapshot
 matrices associated with the training set $\tilde\Dcal_n$.
 \textit{Right}: Convergence of the maximum relative $L^2(\rdom)$
 error over the test set $\tilde\Dcal_{101}$ for the fixed-domain ROM
 (\ref{line:advec1v2d_taddei0_study1_svals_romfix}) and
 ROM-IFT (\ref{line:advec1v2d_taddei0_study1_svals_romtrk})
 with the size of the training set ($n$) (training set: $\tilde\Dcal_n$).}
\label{fig:advec-react-study1-cnvg}
\end{figure}

\subsubsection{Performance with limited training}
\begin{figure}
\centering
\begin{tikzpicture}
\begin{groupplot}[
  group style={
      group size=3 by 1,
      horizontal sep=1cm
  },
  width=0.39\textwidth,
  axis equal image,
  xlabel={$x_1$},
  ylabel={$x_2$},
  xtick = {0.0, 0.5, 1.0},
  ytick = {0.0, 0.5, 1.0},
  xmin=0, xmax=1,
  ymin=0, ymax=1
]
\nextgroupplot
\addplot graphics [xmin=0, xmax=1, ymin=0, ymax=1] {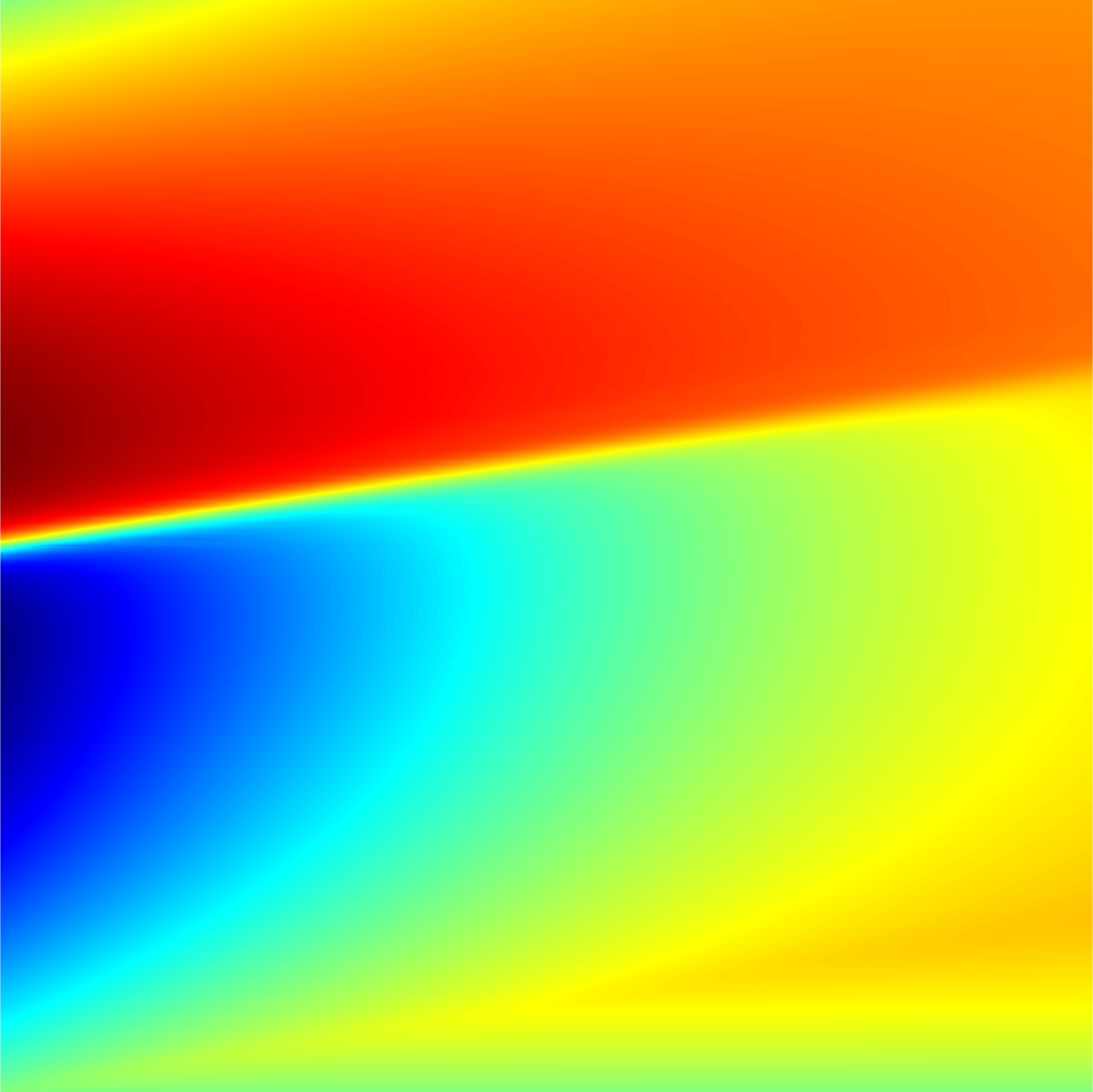};

\nextgroupplot[ylabel={}, ytick=\empty]
\addplot graphics [xmin=0, xmax=1, ymin=0, ymax=1] {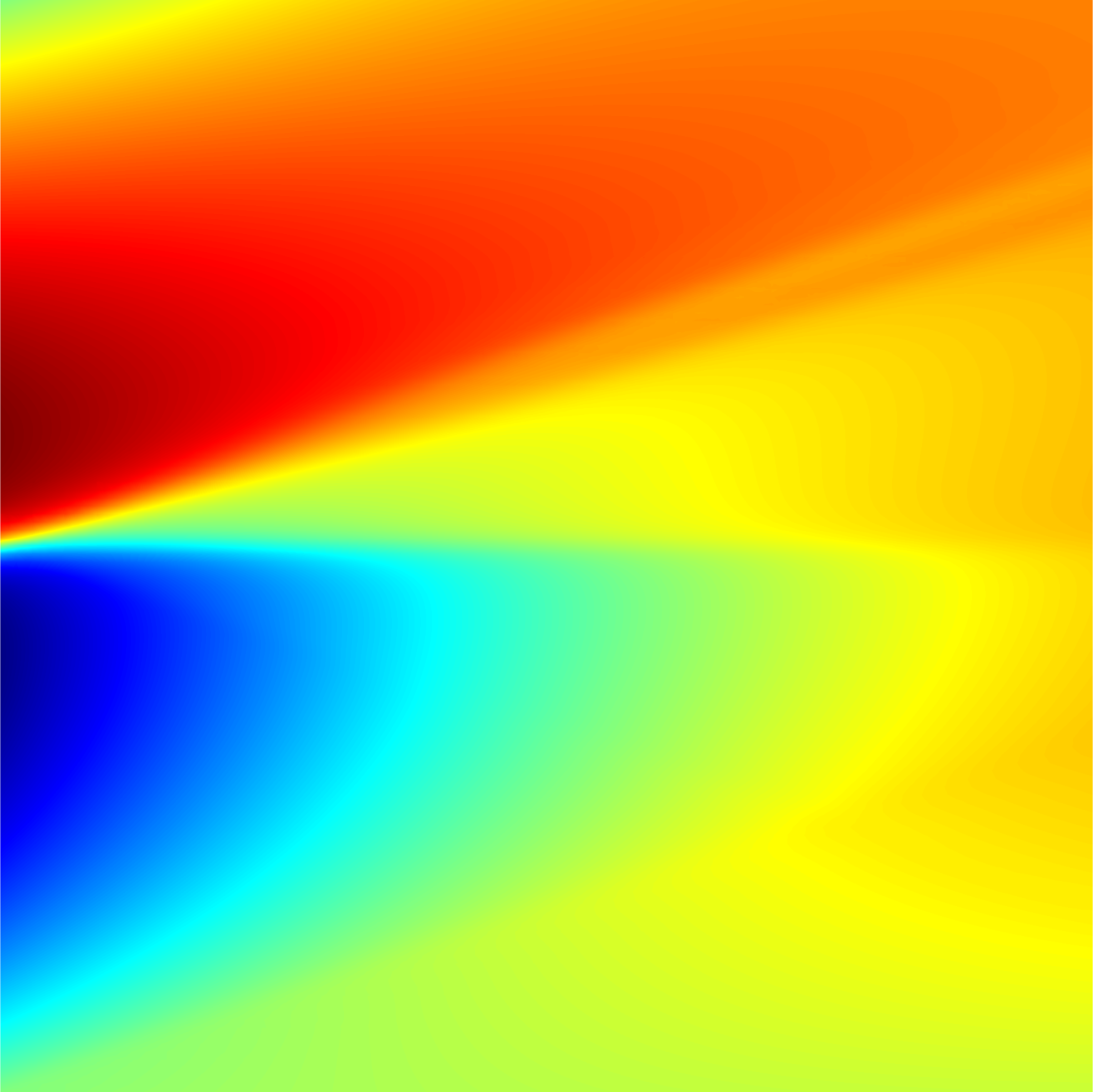};

\nextgroupplot[ylabel={}, ytick=\empty]
\addplot graphics [xmin=0, xmax=1, ymin=0, ymax=1] {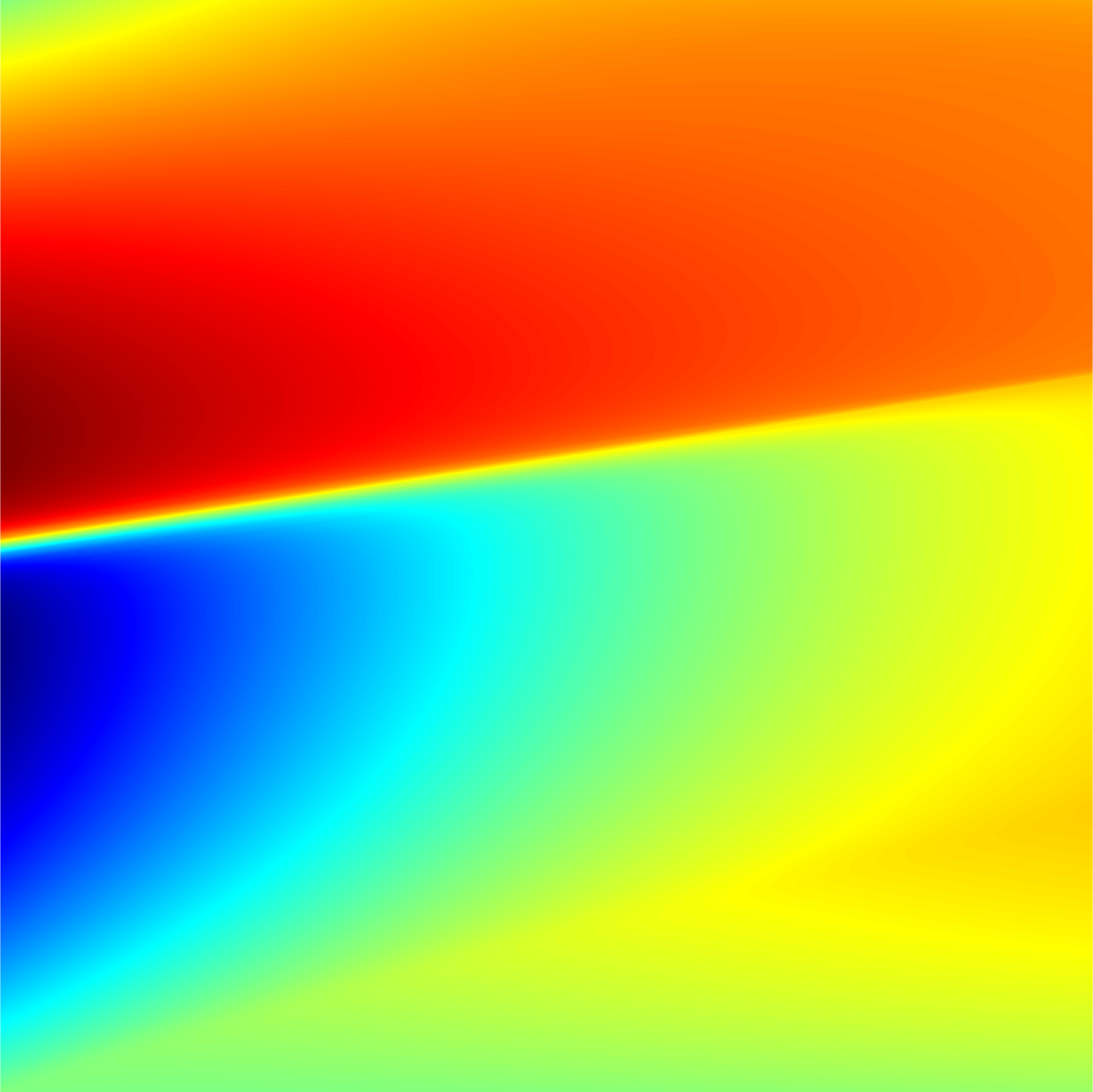};
\end{groupplot}
\end{tikzpicture}
\caption{Advection-reaction solution at $\mu=(\pi/20, 0.3, 100)$
 using the HDM (\textit{left}), fixed-domain ROM (\textit{middle}),
 and ROM-IFT (\textit{right}). The ROM and ROM-IFT use a bases of
 size $k=9$ constructed from the training set $\Dcal_2$ using their
 respective offline procedure. Colorbar in Figure~\ref{fig:advec-react}.}
\label{fig:advec-react-stdy2}
\end{figure}

Finally, we return to the original problem of
predicting the solution of (\ref{eqn:advec-react})
over the three-dimensional parameter space $\Dcal$.
To this end, define $\Dcal_n$ as the union of centroid
of $\Dcal$, $\bar\mu\coloneqq (0, 0.55, 80)$, and the
uniform sampling of $\Dcal$ with $n$ points per dimension
for a total of $n^3$ ($n^3+1$) points if $n$ is odd (even).
We consider training sets $\Dcal_\mathrm{tst}=\Dcal_n$ for $n=1,2$
and a single test set $\Dcal_\mathrm{tst}=\Dcal_5$. Due to the limited
training, we do not truncate the POD modes, leading to reduced
bases of sizes $k=1, 9$ for $n=1, 2$, respectively. For a training
set of size $n=1$, both methods have relatively large worst-case
error over the test set with the ROM-IFT error less than half of
the fixed-domain ROM error: $E_1^\text{rom}(\Dcal_\mathrm{tst}) = 59.0\%$
and $E_1^\text{ift}(\Dcal_\mathrm{tst})=26.0\%$. For a training
set of size $n=9$, the worst-case ROM-IFT error has dropped
to only $E_9^\text{ift}(\Dcal_\mathrm{tst})=5.5\%$, whereas the
fixed-domain ROM error is still large $E_9^\text{rom}(\Dcal_\mathrm{tst})=24.3\%$.
Similar to the conclusion from Section~\ref{sec:numexp:advec:stdy1},
this shows that with limited training the ROM-IFT produces
accurate approximations to the solution manifold of
(\ref{eqn:advec-react}); in contrast, the worst-case error
of fixed-domain ROMs is large with limited training. Even
in the case with $n=2$, the fixed-domain ROM struggles to
represent the primary feature in the right half of the
domain (where its motion is largest), whereas it is
well-approximated by the ROM-IFT solution;
Figure~\ref{fig:advec-react-stdy2} illustrates this
for $\mu=(\pi/20,0.3,100)\in\Dcal_\mathrm{tst}$ and $n=2$.

\subsection{Transonic flow through nozzle}
\label{sec:numexp:nozzle}
Next, we consider inviscid, transonic flow of an ideal gas through an
area variation $A(x)$, commonly used to model flow through ducts, pipes,
shock tubes, and nozzles \cite{toro2013riemann}. Consider a conservation law
of the form (\ref{eqn:claw-phys}) defined over the one-dimensional domain
$\Omega\coloneqq[0, L]$ with solution
$\func{u}{\Omega\times\Dcal}{\Rbb^3}$ with components
\begin{equation}\label{eqn:nozzle-stvc}
u = \begin{bmatrix} A \rho \\ A \rho v \\ A \rho E \end{bmatrix},
\end{equation}
where
$\func{\rho}{\Omega\times\Dcal}{\Rbb_{>0}}$,
$\func{v}{\Omega\times\Dcal}{\Rbb}$,
$\func{E}{\Omega\times\Dcal}{\Rbb_{>0}}$
are the density, velocity, and energy of the fluid, respectively,
and $\func{A}{\Omega\times\Dcal}{\Rbb}$ is the area variation.
The flux function and source term are defined as
\begin{equation} \label{eqn:nozzle-srcflux}
 f(u) \coloneqq \begin{bmatrix} A\rho v \\ A (\rho v^2+P(u)) \\ Av(\rho E+P(u))\end{bmatrix},
 \qquad
 s(u;A) \coloneqq \begin{bmatrix} 0 \\ P(u) A_{,x} \\ 0 \end{bmatrix},
\end{equation}
where $u\mapsto P(u) \coloneqq (\gamma-1)(\rho E - \rho v^2/2)$ is the
pressure of the fluid and $\gamma\in\Rbb_{>0}$ is the ratio of specific
heats. Following the work in \cite{nair_transported_2019}, we take
$L=10$ and $\gamma=1.4$; the nozzle is a parabola with parametrized
throat area at $x=L/2$
\begin{equation}
 (x,\mu)\in\Omega\times\Dcal \mapsto
 A(x;\mu)\coloneqq 3+4(\mu-3) (x/L) (1-x/L),
\end{equation}
where the nozzle throat area varies in $\Dcal\coloneqq[0.5,1.625]$.
At the left boundary ($x=0$), we have a subsonic inlet with a
prescribed density $\rho_0=1$ and pressure $P_0=1$. At the right
boundary ($x=L$), we have a subsonic outlet with a prescribed
pressure $P_L=0.7$. The solution has a strong shock in the right half
of the domain that moves as the throat area varies
(Figure~\ref{fig:eulernozii1d_stdy1_snaps_basis}, \textit{top left}).
We order all training sets $\Dcal_\mathrm{tr}\subset\Dcal$ such that
its elements are monotonically increasing.

The system of conservation laws is discretized using the High-Order
Implicit Shock Tracking (HOIST) method
\cite{zahr_implicit_2020,huang2021robust} using $200$ elements of
polynomial degree $p=2$ for a total of $N=1800$ degrees of freedom.
The Roe flux \cite{roe1981approximate} with Harten-Hyman entropy fix
\cite{harten1983self} is used for the inviscid numerical flux.
We choose the reference domain to be identical to the physical domain
$\rdom=\pdom=[0, L]$ and the nominal mapping to be the identity mapping
$\bar\Gcal=\text{Id}$. The HOIST method is used to define a mesh of the
domain $[0,L]$ that conforms to the shock at $\mu=0.5$ (the first parameter
of any training set). Finally, because there is no benefit to high-order
approximations of the geometry in $d=1$, we take $q=1$, which defines
$\Gbb_{h,q}^\text{b}$ as the collection of functions $\Gcal$
such that $\Gcal(0) = 0$ and $\Gcal(L) = L$ (to preserve the
boundaries).

In the remainder, we repeat the studies from the previous section to
demonstrate the benefits of ROM-IFT hold for this nonlinear problem.
Unlike the previous section, we compare the ROM and ROM-IFT solutions
obtained at a given parameter to a highly refined reference solution of
the governing equations (\ref{eqn:nozzle-srcflux}), rather than the HDM on
a fixed mesh. This is because the HDM is an implicit shock tracking method
that naturally adapts the mesh.
As such, we do not expect the ROM approximation error to tend to
zero; at best, it will plateau at the truncation error associated with
the snapshots, which is $\Ocal(10^{-6})$ in this case. Finally, because
this problem contains a discontinuity (as opposed to the previous problem
that contained a steep feature), we use the $L^1$ norm to avoid excessive
penalization at the location of the shock, which is common practice
\cite{lee1999spurious,powers2006exact,bonfiglioli2014convergence,zahr_optimization-based_2018,zahr_implicit_2020}.

\subsubsection{Convergence under basis refinement, fixed training set}
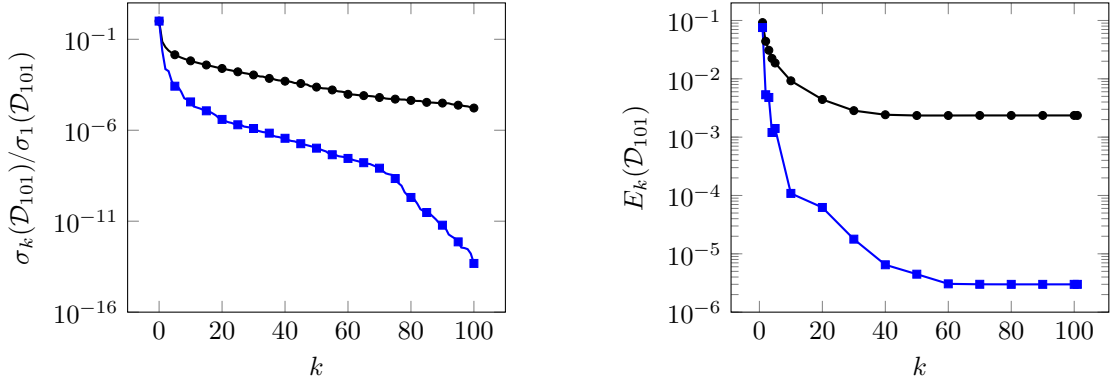
\begin{figure}[t]
 \centering
 \begin{tikzpicture}
\begin{groupplot} [
group style={group size = 2 by 1, horizontal sep = 3cm}]
\nextgroupplot[width=0.4\textwidth, xlabel=$k$, ymax=10, ylabel={$\sigma_k(\Dcal_{101})/\sigma_1(\Dcal_{101})$}, ymode=log, ymin=1e-16]
\addplot [black, thick, mark options={solid, thin}, mark=*, mark size=1.5, mark repeat={5}]
coordinates {
( 0.00000000e+00,  1.00000000e+00)
( 1.00000000e+00,  7.39814638e-02)
( 2.00000000e+00,  3.71621985e-02)
( 3.00000000e+00,  2.40718174e-02)
( 4.00000000e+00,  1.79083557e-02)
( 5.00000000e+00,  1.41118947e-02)
( 6.00000000e+00,  1.16342483e-02)
( 7.00000000e+00,  9.82726107e-03)
( 8.00000000e+00,  8.46802021e-03)
( 9.00000000e+00,  7.40272940e-03)
( 1.00000000e+01,  6.53773075e-03)
( 1.10000000e+01,  5.81067486e-03)
( 1.20000000e+01,  5.21277861e-03)
( 1.30000000e+01,  4.70854624e-03)
( 1.40000000e+01,  4.27205251e-03)
( 1.50000000e+01,  3.86466983e-03)
( 1.60000000e+01,  3.51360800e-03)
( 1.70000000e+01,  3.23057466e-03)
( 1.80000000e+01,  2.92480646e-03)
( 1.90000000e+01,  2.69731799e-03)
( 2.00000000e+01,  2.48884765e-03)
( 2.10000000e+01,  2.26929653e-03)
( 2.20000000e+01,  2.08277155e-03)
( 2.30000000e+01,  1.93549978e-03)
( 2.40000000e+01,  1.76218002e-03)
( 2.50000000e+01,  1.62997105e-03)
( 2.60000000e+01,  1.51746116e-03)
( 2.70000000e+01,  1.38005909e-03)
( 2.80000000e+01,  1.29961893e-03)
( 2.90000000e+01,  1.17760747e-03)
( 3.00000000e+01,  1.09379214e-03)
( 3.10000000e+01,  1.02584886e-03)
( 3.20000000e+01,  9.03053157e-04)
( 3.30000000e+01,  8.88256988e-04)
( 3.40000000e+01,  7.79287398e-04)
( 3.50000000e+01,  7.07512435e-04)
( 3.60000000e+01,  6.85168676e-04)
( 3.70000000e+01,  5.93862369e-04)
( 3.80000000e+01,  5.49113817e-04)
( 3.90000000e+01,  5.36621489e-04)
( 4.00000000e+01,  4.96933721e-04)
( 4.10000000e+01,  4.84113707e-04)
( 4.20000000e+01,  4.52081328e-04)
( 4.30000000e+01,  3.96523259e-04)
( 4.40000000e+01,  3.85145941e-04)
( 4.50000000e+01,  3.70523871e-04)
( 4.60000000e+01,  3.36920713e-04)
( 4.70000000e+01,  3.24463192e-04)
( 4.80000000e+01,  2.63642940e-04)
( 4.90000000e+01,  2.51929727e-04)
( 5.00000000e+01,  2.35051041e-04)
( 5.10000000e+01,  2.17804917e-04)
( 5.20000000e+01,  2.01366104e-04)
( 5.30000000e+01,  1.97789058e-04)
( 5.40000000e+01,  1.88867649e-04)
( 5.50000000e+01,  1.63787699e-04)
( 5.60000000e+01,  1.47676039e-04)
( 5.70000000e+01,  1.36760893e-04)
( 5.80000000e+01,  1.31579063e-04)
( 5.90000000e+01,  1.22364421e-04)
( 6.00000000e+01,  9.44432345e-05)
( 6.10000000e+01,  9.20806975e-05)
( 6.20000000e+01,  9.05136514e-05)
( 6.30000000e+01,  8.61558189e-05)
( 6.40000000e+01,  8.24791353e-05)
( 6.50000000e+01,  7.99399993e-05)
( 6.60000000e+01,  7.74533707e-05)
( 6.70000000e+01,  7.52585371e-05)
( 6.80000000e+01,  7.28004854e-05)
( 6.90000000e+01,  6.68730444e-05)
( 7.00000000e+01,  6.30133700e-05)
( 7.10000000e+01,  6.18178946e-05)
( 7.20000000e+01,  5.52355115e-05)
( 7.30000000e+01,  5.32491014e-05)
( 7.40000000e+01,  5.27608093e-05)
( 7.50000000e+01,  5.13642575e-05)
( 7.60000000e+01,  5.05081637e-05)
( 7.70000000e+01,  4.86216671e-05)
( 7.80000000e+01,  4.84911202e-05)
( 7.90000000e+01,  4.46769990e-05)
( 8.00000000e+01,  4.36070642e-05)
( 8.10000000e+01,  4.15961363e-05)
( 8.20000000e+01,  4.06571605e-05)
( 8.30000000e+01,  4.04951561e-05)
( 8.40000000e+01,  3.67094649e-05)
( 8.50000000e+01,  3.41449609e-05)
( 8.60000000e+01,  3.31888395e-05)
( 8.70000000e+01,  3.30045158e-05)
( 8.80000000e+01,  3.19453742e-05)
( 8.90000000e+01,  3.11058453e-05)
( 9.00000000e+01,  3.10820603e-05)
( 9.10000000e+01,  3.04444103e-05)
( 9.20000000e+01,  3.01018892e-05)
( 9.30000000e+01,  2.68277031e-05)
( 9.40000000e+01,  2.54509814e-05)
( 9.50000000e+01,  2.38215226e-05)
( 9.60000000e+01,  2.29835193e-05)
( 9.70000000e+01,  2.13881722e-05)
( 9.80000000e+01,  2.04726979e-05)
( 9.90000000e+01,  1.83281623e-05)
( 1.00000000e+02,  1.66116654e-05)};\label{line:eulernozii1d_study0_svals_romfix}

\addplot [blue, thick, mark options={solid, thin}, mark=square*, mark size=1.5, mark repeat={5}]
coordinates {
( 0.00000000e+00,  1.00000000e+00)
( 1.00000000e+00,  2.40501157e-02)
( 2.00000000e+00,  2.29679246e-03)
( 3.00000000e+00,  1.77735553e-03)
( 4.00000000e+00,  5.04409503e-04)
( 5.00000000e+00,  2.64174749e-04)
( 6.00000000e+00,  2.45207915e-04)
( 7.00000000e+00,  1.32717219e-04)
( 8.00000000e+00,  5.71072766e-05)
( 9.00000000e+00,  4.77135183e-05)
( 1.00000000e+01,  3.60690660e-05)
( 1.10000000e+01,  2.01912894e-05)
( 1.20000000e+01,  1.77897710e-05)
( 1.30000000e+01,  1.56239897e-05)
( 1.40000000e+01,  1.32363489e-05)
( 1.50000000e+01,  1.16888209e-05)
( 1.60000000e+01,  1.05646751e-05)
( 1.70000000e+01,  8.45049966e-06)
( 1.80000000e+01,  6.61880487e-06)
( 1.90000000e+01,  4.61043264e-06)
( 2.00000000e+01,  3.94529838e-06)
( 2.10000000e+01,  3.21949543e-06)
( 2.20000000e+01,  2.76319918e-06)
( 2.30000000e+01,  2.53515679e-06)
( 2.40000000e+01,  2.15376766e-06)
( 2.50000000e+01,  1.99358700e-06)
( 2.60000000e+01,  1.97432244e-06)
( 2.70000000e+01,  1.58968115e-06)
( 2.80000000e+01,  1.49346245e-06)
( 2.90000000e+01,  1.40412358e-06)
( 3.00000000e+01,  1.24460795e-06)
( 3.10000000e+01,  9.88358241e-07)
( 3.20000000e+01,  8.76559858e-07)
( 3.30000000e+01,  7.90819347e-07)
( 3.40000000e+01,  7.62684504e-07)
( 3.50000000e+01,  6.91929504e-07)
( 3.60000000e+01,  5.13655280e-07)
( 3.70000000e+01,  4.48460927e-07)
( 3.80000000e+01,  4.22425825e-07)
( 3.90000000e+01,  3.81390257e-07)
( 4.00000000e+01,  3.61717737e-07)
( 4.10000000e+01,  2.87741988e-07)
( 4.20000000e+01,  2.57438494e-07)
( 4.30000000e+01,  2.45176033e-07)
( 4.40000000e+01,  2.13228506e-07)
( 4.50000000e+01,  1.81595357e-07)
( 4.60000000e+01,  1.65975745e-07)
( 4.70000000e+01,  1.49947682e-07)
( 4.80000000e+01,  1.34890745e-07)
( 4.90000000e+01,  1.21838952e-07)
( 5.00000000e+01,  1.03524938e-07)
( 5.10000000e+01,  9.15378956e-08)
( 5.20000000e+01,  7.84064908e-08)
( 5.30000000e+01,  6.84065372e-08)
( 5.40000000e+01,  5.14900599e-08)
( 5.50000000e+01,  4.49331241e-08)
( 5.60000000e+01,  4.41390136e-08)
( 5.70000000e+01,  3.74214749e-08)
( 5.80000000e+01,  3.52667043e-08)
( 5.90000000e+01,  3.32908280e-08)
( 6.00000000e+01,  2.82507830e-08)
( 6.10000000e+01,  2.53060353e-08)
( 6.20000000e+01,  2.22628638e-08)
( 6.30000000e+01,  2.01772373e-08)
( 6.40000000e+01,  1.80104368e-08)
( 6.50000000e+01,  1.64020831e-08)
( 6.60000000e+01,  1.46682206e-08)
( 6.70000000e+01,  1.38361780e-08)
( 6.80000000e+01,  1.16759618e-08)
( 6.90000000e+01,  9.82688911e-09)
( 7.00000000e+01,  8.20620945e-09)
( 7.10000000e+01,  5.91801461e-09)
( 7.20000000e+01,  4.27356761e-09)
( 7.30000000e+01,  3.80846350e-09)
( 7.40000000e+01,  3.52282478e-09)
( 7.50000000e+01,  2.21328772e-09)
( 7.60000000e+01,  1.25107563e-09)
( 7.70000000e+01,  8.96239255e-10)
( 7.80000000e+01,  3.89517831e-10)
( 7.90000000e+01,  2.54266550e-10)
( 8.00000000e+01,  2.01935019e-10)
( 8.10000000e+01,  1.62367703e-10)
( 8.20000000e+01,  9.32338671e-11)
( 8.30000000e+01,  4.57144052e-11)
( 8.40000000e+01,  3.60621008e-11)
( 8.50000000e+01,  2.99740321e-11)
( 8.60000000e+01,  2.56102768e-11)
( 8.70000000e+01,  1.78389142e-11)
( 8.80000000e+01,  1.22106206e-11)
( 8.90000000e+01,  7.59484472e-12)
( 9.00000000e+01,  6.03835821e-12)
( 9.10000000e+01,  4.10321623e-12)
( 9.20000000e+01,  1.85208757e-12)
( 9.30000000e+01,  1.42309463e-12)
( 9.40000000e+01,  1.11466430e-12)
( 9.50000000e+01,  7.30679313e-13)
( 9.60000000e+01,  3.56483785e-13)
( 9.70000000e+01,  3.34048906e-13)
( 9.80000000e+01,  2.91275959e-13)
( 9.90000000e+01,  1.65823929e-13)
( 1.00000000e+02,  4.76668471e-14)};\label{line:eulernozii1d_study0_svals_romtrk}

\nextgroupplot[width=0.4\textwidth, xlabel=$k$, ymax=0.2, ylabel={$E_k(\Dcal_{101})$}, ymode=log, ymin=1e-06]
\addplot [black, thick, mark options={solid, thin}, mark=*, mark size=1.5, mark repeat={1}]
coordinates {
( 1.00000000e+00,  9.20086856e-02)
( 2.00000000e+00,  4.39616386e-02)
( 3.00000000e+00,  3.08331986e-02)
( 4.00000000e+00,  2.23638363e-02)
( 5.00000000e+00,  1.86529271e-02)
( 1.00000000e+01,  9.24231654e-03)
( 2.00000000e+01,  4.41873231e-03)
( 3.00000000e+01,  2.84596538e-03)
( 4.00000000e+01,  2.42008159e-03)
( 5.00000000e+01,  2.34843090e-03)
( 6.00000000e+01,  2.34914277e-03)
( 7.00000000e+01,  2.35492360e-03)
( 8.00000000e+01,  2.35500244e-03)
( 9.00000000e+01,  2.35498336e-03)
( 1.00000000e+02,  2.35498645e-03)
( 1.01000000e+02,  2.35498739e-03)};\label{line:eulernozii1d_study0_L2err_romfix}

\addplot [blue, thick, mark options={solid, thin}, mark=square*, mark size=1.5, mark repeat={1}]
coordinates {
( 1.00000000e+00,  7.55804081e-02)
( 2.00000000e+00,  5.33181193e-03)
( 3.00000000e+00,  4.78741474e-03)
( 4.00000000e+00,  1.20175718e-03)
( 5.00000000e+00,  1.40747428e-03)
( 1.00000000e+01,  1.08541378e-04)
( 2.00000000e+01,  6.23034086e-05)
( 3.00000000e+01,  1.78400522e-05)
( 4.00000000e+01,  6.48683694e-06)
( 5.00000000e+01,  4.47427743e-06)
( 6.00000000e+01,  3.05630638e-06)
( 7.00000000e+01,  3.00348896e-06)
( 8.00000000e+01,  2.99513195e-06)
( 9.00000000e+01,  2.99513204e-06)
( 1.00000000e+02,  2.99513203e-06)
( 1.01000000e+02,  2.99513203e-06)};\label{line:eulernozii1d_study0_L2err_romtrk}

\end{groupplot}\end{tikzpicture}
 \caption{\textit{Left}: Convergence of the singular values of the non-aligned
 (\ref{line:eulernozii1d_study0_svals_romfix}) and aligned
 (\ref{line:eulernozii1d_study0_svals_romtrk}) snapshot matrices
 associated with the training set $\Dcal_{101}$.
 \textit{Right}: Convergence of the maximum relative
 $L^1(\Omega_0)$ error over the test set $\Dcal_{101}$ for the fixed-domain
 ROM (\ref{line:eulernozii1d_study0_L2err_romfix}) and ROM-IFT
 (\ref{line:eulernozii1d_study0_L2err_romtrk}) with the size of the reduced
 basis ($k$)
 (training set: $\Dcal_{101}$).}
 \label{fig:eulernozii1d_stdy0_cnvg}
\end{figure}

First, we demonstrate rapid convergence of the ROM-IFT method relative
to fixed-domain model reduction in the idealized setting where the training
and testing set are identical. To this end, we construct a ROM and ROM-IFT
from a collection of $M=101$ parameters ($\Dcal_{101}$). Similar to the
previous problem, we define $\Dcal_n$ to be the subset of $\Dcal$
containing $n$ uniformly spaced samples. The snapshot
alignment procedure of the ROM-IFT method causes the singular values
of the corresponding snapshot matrix to decay significantly
faster than without alignment. The decrease in the maximum
$L^1(\Omega_0)$ error over the test set ($\Dcal_{101}$) is also significantly
faster for the ROM-IFT method (Figure~\ref{fig:eulernozii1d_stdy0_cnvg}).
For a fixed basis dimension $k$, the ROM-IFT method has a much smaller
$L^1(\Omega_0)$ error than a standard reduced-order model.

\begin{remark}
The fixed-domain ROM method plateaus at an error around $10^{-2}$ because
a shock capturing method \cite{persson2006sub} with $200$ quartic elements
is used (discontinuity represented as a steep feature, which causes the
increased truncation error). Implicit shock tracking would
require interpolating all snapshots onto a unique mesh prior to compression,
which is expensive given the number of snapshots considered.
\end{remark}

\begin{figure}
 \centering
 \input{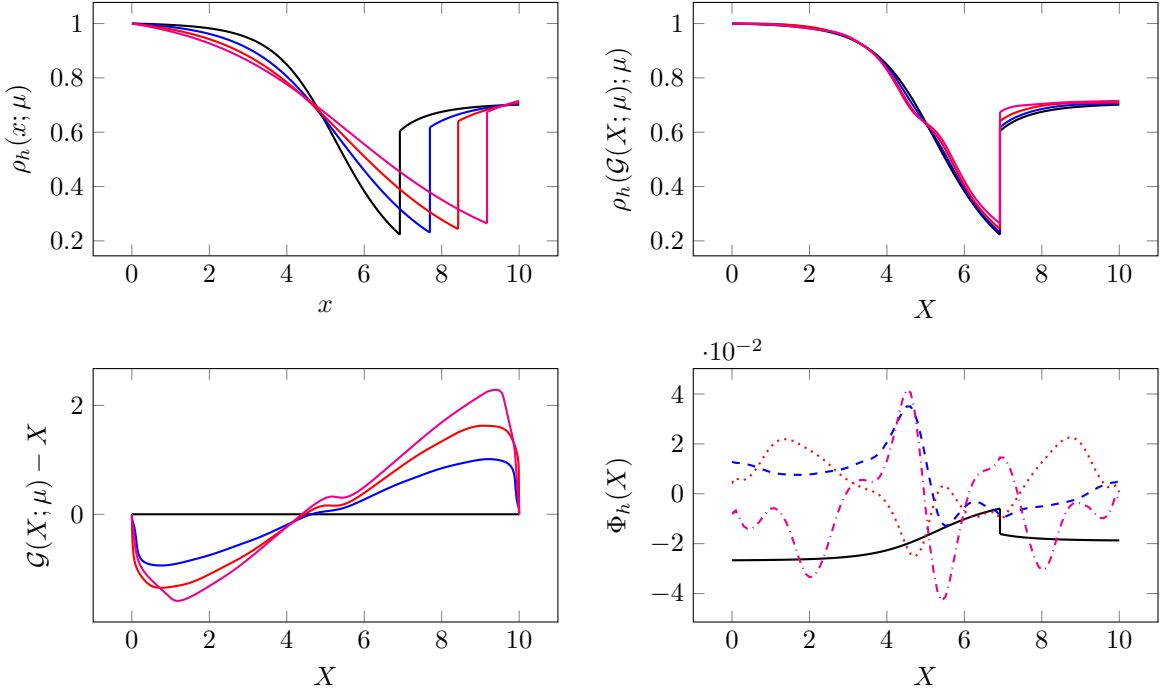}
 \caption{Four snapshots of flow through the nozzle (density) associated
 with the training set $\Dcal_4$ in the physical domain (\textit{top left}),
 reference domain (\textit{top right}), and the domain deformation that
 causes the shocks to align in the reference domain (\textit{bottom left});
 the color corresponds to the same parameter value ($\mu$) across the three
 plots. The aligned snapshots are compressed and the resulting modes are
 shown in the \textit{bottom right}: first (\ref{line:eulernozii1d_basis0}),
 second (\ref{line:eulernozii1d_basis1}),
 third (\ref{line:eulernozii1d_basis2}), and
 fourth (\ref{line:eulernozii1d_basis3}) mode.}
 \label{fig:eulernozii1d_stdy1_snaps_basis}
\end{figure}

In the reference domain, the shock is located at the same spatial coordinate
($x\approx 6.92$) for all aligned snapshots generated by the ROM-IFT method
(induced by a parameter-dependent domain deformation), despite the
parameter-dependent position of the shock in the physical domain 
(Figure~\ref{fig:eulernozii1d_stdy1_snaps_basis}). This effectively
removes the convection-dominated nature of the solution, which enhances
its compressibility as seen in Figure~\ref{fig:eulernozii1d_stdy0_cnvg}.
Furthermore, this allows the first POD mode of the aligned snapshots to
represent the shock and all remaining modes resolve features on either
side of it (Figure~\ref{fig:eulernozii1d_stdy1_snaps_basis}).

\subsubsection{Convergence under training set refinement}
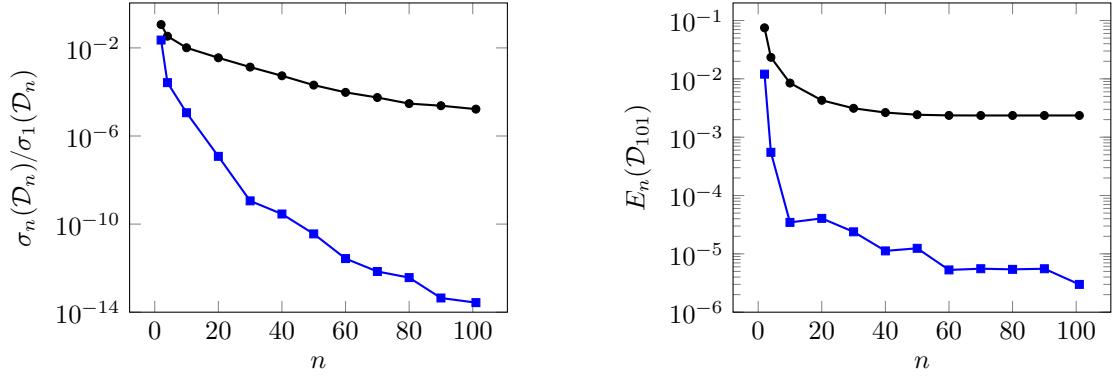
\begin{figure}[t]
 \centering
 \begin{tikzpicture}
\begin{groupplot} [
group style={group size = 2 by 3, horizontal sep = 3cm, vertical sep = 1.5cm}]
\nextgroupplot[width=0.4\textwidth, xlabel=$n$, ymax=1.1, ylabel={$\sigma_n(\Dcal_n)/\sigma_1(\Dcal_n)$}, ymode=log, ymin=1e-14]
\addplot [black, thick, mark options={solid, thin}, mark=*, mark size=1.5]
coordinates {
( 2.00000000e+00,  1.15150730e-01)
( 4.00000000e+00,  3.39643925e-02)
( 1.00000000e+01,  1.01459518e-02)
( 2.00000000e+01,  3.56084705e-03)
( 3.00000000e+01,  1.34566297e-03)
( 4.00000000e+01,  5.41110301e-04)
( 5.00000000e+01,  2.04523536e-04)
( 6.00000000e+01,  9.58051085e-05)
( 7.00000000e+01,  5.59443792e-05)
( 8.00000000e+01,  2.95747104e-05)
( 9.00000000e+01,  2.37330825e-05)
( 1.01000000e+02,  1.66116654e-05)};\label{line:eulernozii1d_study1_svals_romfix}

\addplot [blue, thick, mark options={solid, thin}, mark=square*, mark size=1.5]
coordinates {
( 2.00000000e+00,  2.27109350e-02)
( 4.00000000e+00,  2.64634416e-04)
( 1.00000000e+01,  1.13770444e-05)
( 2.00000000e+01,  1.19452997e-07)
( 3.00000000e+01,  1.13407298e-09)
( 4.00000000e+01,  2.87564828e-10)
( 5.00000000e+01,  3.58416388e-11)
( 6.00000000e+01,  2.71472229e-12)
( 7.00000000e+01,  7.06055465e-13)
( 8.00000000e+01,  3.71602811e-13)
( 9.00000000e+01,  4.41431863e-14)
( 1.01000000e+02,  2.71346551e-14)};\label{line:eulernozii1d_study1_svals_romtrk}

\nextgroupplot[width=0.4\textwidth, xlabel=$n$, ymax=0.2, ylabel={$E_n(\Dcal_{101})$}, ymode=log, ymin=1e-06]
\addplot [black, thick, mark options={solid, thin}, mark=*, mark size=1.5]
coordinates {
( 2.00000000e+00,  7.49747340e-02)
( 4.00000000e+00,  2.33805101e-02)
( 1.00000000e+01,  8.46816805e-03)
( 2.00000000e+01,  4.28571679e-03)
( 3.00000000e+01,  3.12679586e-03)
( 4.00000000e+01,  2.64372895e-03)
( 5.00000000e+01,  2.42394980e-03)
( 6.00000000e+01,  2.36219908e-03)
( 7.00000000e+01,  2.35498739e-03)
( 8.00000000e+01,  2.35498739e-03)
( 9.00000000e+01,  2.35498739e-03)
( 1.01000000e+02,  2.35498739e-03)};\label{line:eulernozii1d_study1_L2err_romfix}

\addplot [blue, thick, mark options={solid, thin}, mark=square*, mark size=1.5]
coordinates {
( 2.00000000e+00,  1.20369258e-02)
( 4.00000000e+00,  5.50418100e-04)
( 1.00000000e+01,  3.46099159e-05)
( 2.00000000e+01,  4.05652249e-05)
( 3.00000000e+01,  2.38504121e-05)
( 4.00000000e+01,  1.12373318e-05)
( 5.00000000e+01,  1.24533163e-05)
( 6.00000000e+01,  5.30186875e-06)
( 7.00000000e+01,  5.55963262e-06)
( 8.00000000e+01,  5.42606716e-06)
( 9.00000000e+01,  5.55336464e-06)
( 1.01000000e+02,  2.99513035e-06)};\label{line:eulernozii1d_study1_L2err_romtrk}

\end{groupplot}\end{tikzpicture}
 \caption{\textit{Left}: Convergence of the smallest singular value of the
  non-aligned (\ref{line:eulernozii1d_study1_svals_romfix}) and aligned
  (\ref{line:eulernozii1d_study1_svals_romtrk}) snapshot matrices associated
  with the training set $\Dcal_n$. \textit{Right}: Convergence of the maximum
  relative $L^1(\Omega_0)$ error over the test set $\Dcal_{101}$ for the
  fixed-domain ROM (\ref{line:eulernozii1d_study1_L2err_romfix}) and ROM-IFT
  (\ref{line:eulernozii1d_study1_L2err_romtrk}) with the size of the
  training set ($n$) (training set: $\Dcal_n$).}
 \label{fig:eulernozii1d_stdy1_cnvg}
\end{figure}

\begin{figure}
 \centering
 \input{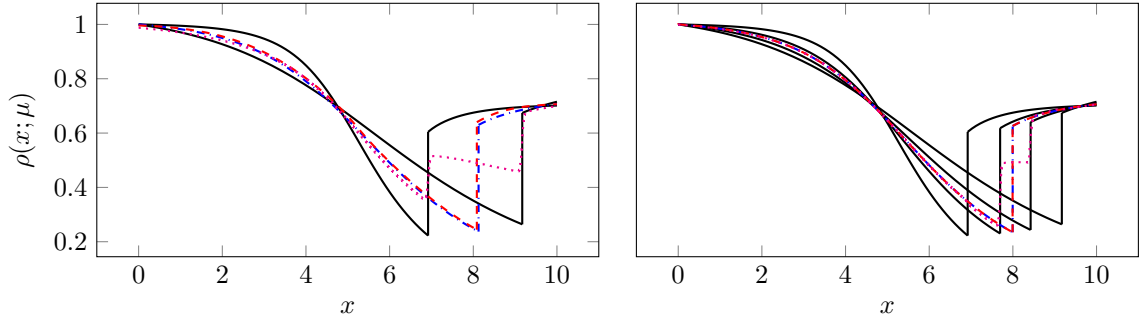}
 \caption{Comparison of the HDM, $L^2$ projection (without snapshot alignment),
 and ROM-IFT on the nozzle flow problem (density) with two
 (\textit{left}) and four (\textit{right}) training parameters. The test
 set includes $101$ parameter configuration ($\Dcal_{101}$), although only
 the parameter that results in the largest error for the ROM-IFT method
 is shown. Legend: training snapshots (\ref{line:eulernozii1d_states_snaps}),
 HDM solution at test parameter (\ref{line:eulernozii1d_states_hdm}),
 $L^2$ projection of test parameter onto reduced basis
 (\ref{line:eulernozii1d_states_l2}), and ROM-IFT solution
 at test parameter (\ref{line:eulernozii1d_states_ift}).}
 \label{fig:eulernozii1d_stdy1_states}
\end{figure}

Next, we consider the relative performance of fixed-domain model
reduction and the ROM-IFT method as the training set is refined.
In particular, we take the training set to be $\Dcal_n$ for
$1\leq n \leq 101$ and the testing set to be $\Dcal_{101}$. In
all cases, we do not truncate the POD modes, which leads to a
basis of dimension $n$. Similar to the previous problem, we
see the smallest singular value of the aligned snapshot matrix
decays significantly more rapidly than those associated with
the non-aligned snapshot matrix (Figure~\ref{fig:eulernozii1d_stdy1_cnvg}),
suggesting most of the information of the solution manifold over $\Dcal$
is captured with relatively few training points when snapshots
are aligned. In addition, we observe that the ROM-IFT method
requires little training to obtain small errors over the
entire testing set (Figure~\ref{fig:eulernozii1d_stdy1_cnvg}).
For example, for $n=2$ ($n=4$), the maximum relative $L^1(\Omega_0)$
error over the testing set is 7.5\% (2.34\%) for the fixed-domain ROM and
only 1.2\% (0.055\%) for the ROM-IFT. While even the fixed-domain ROM
errors seem reasonable, the density profiles are non-physical, whereas
the ROM-IFT approximations are highly accurate, even with as little as
$n=2$ training points (Figure~\ref{fig:eulernozii1d_stdy1_states}).




\subsection{Supersonic, inviscid flow over cylinder}
\label{sec:numexp:euler-cyl0}
\begin{figure}
 \centering
 \begin{tikzpicture}
\begin{axis}[
axis equal image,
axis line style={gray},
axis x line*=bottom,
axis y line*=left,
width=0.7\textwidth,
xtick={-8, -1, 0, 1, 8},
ytick={0, 1, 4},
grid=major,
ymax=4.3,
xmax=8.3,
xmin=-8.3,
ymin=-0.3]
\addplot [opacity=0.6, fill=black!30!white, opacity=0.6, forget plot]
coordinates {
( 1.00000000e+00,  0.00000000e+00)
( 9.99496542e-01,  3.17279335e-02)
( 9.97986676e-01,  6.34239197e-02)
( 9.95471923e-01,  9.50560433e-02)
( 9.91954813e-01,  1.26592454e-01)
( 9.87438889e-01,  1.58001396e-01)
( 9.81928697e-01,  1.89251244e-01)
( 9.75429787e-01,  2.20310533e-01)
( 9.67948701e-01,  2.51147987e-01)
( 9.59492974e-01,  2.81732557e-01)
( 9.50071118e-01,  3.12033446e-01)
( 9.39692621e-01,  3.42020143e-01)
( 9.28367933e-01,  3.71662456e-01)
( 9.16108457e-01,  4.00930535e-01)
( 9.02926538e-01,  4.29794912e-01)
( 8.88835449e-01,  4.58226522e-01)
( 8.73849377e-01,  4.86196736e-01)
( 8.57983413e-01,  5.13677392e-01)
( 8.41253533e-01,  5.40640817e-01)
( 8.23676581e-01,  5.67059864e-01)
( 8.05270258e-01,  5.92907929e-01)
( 7.86053095e-01,  6.18158986e-01)
( 7.66044443e-01,  6.42787610e-01)
( 7.45264450e-01,  6.66769001e-01)
( 7.23734038e-01,  6.90079011e-01)
( 7.01474888e-01,  7.12694171e-01)
( 6.78509412e-01,  7.34591709e-01)
( 6.54860734e-01,  7.55749574e-01)
( 6.30552667e-01,  7.76146464e-01)
( 6.05609687e-01,  7.95761841e-01)
( 5.80056910e-01,  8.14575952e-01)
( 5.53920064e-01,  8.32569855e-01)
( 5.27225468e-01,  8.49725430e-01)
( 5.00000000e-01,  8.66025404e-01)
( 4.72271075e-01,  8.81453363e-01)
( 4.44066613e-01,  8.95993774e-01)
( 4.15415013e-01,  9.09631995e-01)
( 3.86345126e-01,  9.22354294e-01)
( 3.56886222e-01,  9.34147860e-01)
( 3.27067963e-01,  9.45000819e-01)
( 2.96920375e-01,  9.54902241e-01)
( 2.66473814e-01,  9.63842159e-01)
( 2.35758936e-01,  9.71811568e-01)
( 2.04806668e-01,  9.78802446e-01)
( 1.73648178e-01,  9.84807753e-01)
( 1.42314838e-01,  9.89821442e-01)
( 1.10838200e-01,  9.93838464e-01)
( 7.92499569e-02,  9.96854776e-01)
( 4.75819158e-02,  9.98867339e-01)
( 1.58659638e-02,  9.99874128e-01)
(-1.58659638e-02,  9.99874128e-01)
(-4.75819158e-02,  9.98867339e-01)
(-7.92499569e-02,  9.96854776e-01)
(-1.10838200e-01,  9.93838464e-01)
(-1.42314838e-01,  9.89821442e-01)
(-1.73648178e-01,  9.84807753e-01)
(-2.04806668e-01,  9.78802446e-01)
(-2.35758936e-01,  9.71811568e-01)
(-2.66473814e-01,  9.63842159e-01)
(-2.96920375e-01,  9.54902241e-01)
(-3.27067963e-01,  9.45000819e-01)
(-3.56886222e-01,  9.34147860e-01)
(-3.86345126e-01,  9.22354294e-01)
(-4.15415013e-01,  9.09631995e-01)
(-4.44066613e-01,  8.95993774e-01)
(-4.72271075e-01,  8.81453363e-01)
(-5.00000000e-01,  8.66025404e-01)
(-5.27225468e-01,  8.49725430e-01)
(-5.53920064e-01,  8.32569855e-01)
(-5.80056910e-01,  8.14575952e-01)
(-6.05609687e-01,  7.95761841e-01)
(-6.30552667e-01,  7.76146464e-01)
(-6.54860734e-01,  7.55749574e-01)
(-6.78509412e-01,  7.34591709e-01)
(-7.01474888e-01,  7.12694171e-01)
(-7.23734038e-01,  6.90079011e-01)
(-7.45264450e-01,  6.66769001e-01)
(-7.66044443e-01,  6.42787610e-01)
(-7.86053095e-01,  6.18158986e-01)
(-8.05270258e-01,  5.92907929e-01)
(-8.23676581e-01,  5.67059864e-01)
(-8.41253533e-01,  5.40640817e-01)
(-8.57983413e-01,  5.13677392e-01)
(-8.73849377e-01,  4.86196736e-01)
(-8.88835449e-01,  4.58226522e-01)
(-9.02926538e-01,  4.29794912e-01)
(-9.16108457e-01,  4.00930535e-01)
(-9.28367933e-01,  3.71662456e-01)
(-9.39692621e-01,  3.42020143e-01)
(-9.50071118e-01,  3.12033446e-01)
(-9.59492974e-01,  2.81732557e-01)
(-9.67948701e-01,  2.51147987e-01)
(-9.75429787e-01,  2.20310533e-01)
(-9.81928697e-01,  1.89251244e-01)
(-9.87438889e-01,  1.58001396e-01)
(-9.91954813e-01,  1.26592454e-01)
(-9.95471923e-01,  9.50560433e-02)
(-9.97986676e-01,  6.34239197e-02)
(-9.99496542e-01,  3.17279335e-02)
(-1.00000000e+00,  1.22464680e-16)
(-8.00000000e+00,  0.00000000e+00)
(-8.00000000e+00,  4.00000000e+00)
( 8.00000000e+00,  4.00000000e+00)
( 8.00000000e+00,  0.00000000e+00)
( 1.00000000e+00,  0.00000000e+00)};

\addplot [thick, color=black]
coordinates {
( 1.00000000e+00,  0.00000000e+00)
( 9.99496542e-01,  3.17279335e-02)
( 9.97986676e-01,  6.34239197e-02)
( 9.95471923e-01,  9.50560433e-02)
( 9.91954813e-01,  1.26592454e-01)
( 9.87438889e-01,  1.58001396e-01)
( 9.81928697e-01,  1.89251244e-01)
( 9.75429787e-01,  2.20310533e-01)
( 9.67948701e-01,  2.51147987e-01)
( 9.59492974e-01,  2.81732557e-01)
( 9.50071118e-01,  3.12033446e-01)
( 9.39692621e-01,  3.42020143e-01)
( 9.28367933e-01,  3.71662456e-01)
( 9.16108457e-01,  4.00930535e-01)
( 9.02926538e-01,  4.29794912e-01)
( 8.88835449e-01,  4.58226522e-01)
( 8.73849377e-01,  4.86196736e-01)
( 8.57983413e-01,  5.13677392e-01)
( 8.41253533e-01,  5.40640817e-01)
( 8.23676581e-01,  5.67059864e-01)
( 8.05270258e-01,  5.92907929e-01)
( 7.86053095e-01,  6.18158986e-01)
( 7.66044443e-01,  6.42787610e-01)
( 7.45264450e-01,  6.66769001e-01)
( 7.23734038e-01,  6.90079011e-01)
( 7.01474888e-01,  7.12694171e-01)
( 6.78509412e-01,  7.34591709e-01)
( 6.54860734e-01,  7.55749574e-01)
( 6.30552667e-01,  7.76146464e-01)
( 6.05609687e-01,  7.95761841e-01)
( 5.80056910e-01,  8.14575952e-01)
( 5.53920064e-01,  8.32569855e-01)
( 5.27225468e-01,  8.49725430e-01)
( 5.00000000e-01,  8.66025404e-01)
( 4.72271075e-01,  8.81453363e-01)
( 4.44066613e-01,  8.95993774e-01)
( 4.15415013e-01,  9.09631995e-01)
( 3.86345126e-01,  9.22354294e-01)
( 3.56886222e-01,  9.34147860e-01)
( 3.27067963e-01,  9.45000819e-01)
( 2.96920375e-01,  9.54902241e-01)
( 2.66473814e-01,  9.63842159e-01)
( 2.35758936e-01,  9.71811568e-01)
( 2.04806668e-01,  9.78802446e-01)
( 1.73648178e-01,  9.84807753e-01)
( 1.42314838e-01,  9.89821442e-01)
( 1.10838200e-01,  9.93838464e-01)
( 7.92499569e-02,  9.96854776e-01)
( 4.75819158e-02,  9.98867339e-01)
( 1.58659638e-02,  9.99874128e-01)
(-1.58659638e-02,  9.99874128e-01)
(-4.75819158e-02,  9.98867339e-01)
(-7.92499569e-02,  9.96854776e-01)
(-1.10838200e-01,  9.93838464e-01)
(-1.42314838e-01,  9.89821442e-01)
(-1.73648178e-01,  9.84807753e-01)
(-2.04806668e-01,  9.78802446e-01)
(-2.35758936e-01,  9.71811568e-01)
(-2.66473814e-01,  9.63842159e-01)
(-2.96920375e-01,  9.54902241e-01)
(-3.27067963e-01,  9.45000819e-01)
(-3.56886222e-01,  9.34147860e-01)
(-3.86345126e-01,  9.22354294e-01)
(-4.15415013e-01,  9.09631995e-01)
(-4.44066613e-01,  8.95993774e-01)
(-4.72271075e-01,  8.81453363e-01)
(-5.00000000e-01,  8.66025404e-01)
(-5.27225468e-01,  8.49725430e-01)
(-5.53920064e-01,  8.32569855e-01)
(-5.80056910e-01,  8.14575952e-01)
(-6.05609687e-01,  7.95761841e-01)
(-6.30552667e-01,  7.76146464e-01)
(-6.54860734e-01,  7.55749574e-01)
(-6.78509412e-01,  7.34591709e-01)
(-7.01474888e-01,  7.12694171e-01)
(-7.23734038e-01,  6.90079011e-01)
(-7.45264450e-01,  6.66769001e-01)
(-7.66044443e-01,  6.42787610e-01)
(-7.86053095e-01,  6.18158986e-01)
(-8.05270258e-01,  5.92907929e-01)
(-8.23676581e-01,  5.67059864e-01)
(-8.41253533e-01,  5.40640817e-01)
(-8.57983413e-01,  5.13677392e-01)
(-8.73849377e-01,  4.86196736e-01)
(-8.88835449e-01,  4.58226522e-01)
(-9.02926538e-01,  4.29794912e-01)
(-9.16108457e-01,  4.00930535e-01)
(-9.28367933e-01,  3.71662456e-01)
(-9.39692621e-01,  3.42020143e-01)
(-9.50071118e-01,  3.12033446e-01)
(-9.59492974e-01,  2.81732557e-01)
(-9.67948701e-01,  2.51147987e-01)
(-9.75429787e-01,  2.20310533e-01)
(-9.81928697e-01,  1.89251244e-01)
(-9.87438889e-01,  1.58001396e-01)
(-9.91954813e-01,  1.26592454e-01)
(-9.95471923e-01,  9.50560433e-02)
(-9.97986676e-01,  6.34239197e-02)
(-9.99496542e-01,  3.17279335e-02)
(-1.00000000e+00,  1.22464680e-16)};\label{line:cyl0:wall}

\addplot [thick, color=red]
coordinates {
(-1.00000000e+00,  1.22464680e-16)
(-8.00000000e+00,  0.00000000e+00)};\label{line:cyl0:outf}

\addplot [thick, color=blue]
coordinates {
(-8.00000000e+00,  0.00000000e+00)
(-8.00000000e+00,  4.00000000e+00)};\label{line:cyl0:farf}

\addplot [thick, color=blue, forget plot]
coordinates {
(-8.00000000e+00,  4.00000000e+00)
( 8.00000000e+00,  4.00000000e+00)};

\addplot [thick, color=blue, forget plot]
coordinates {
( 8.00000000e+00,  4.00000000e+00)
( 8.00000000e+00,  0.00000000e+00)};

\addplot [thick, color=red, forget plot]
coordinates {
( 8.00000000e+00,  0.00000000e+00)
( 1.00000000e+00,  0.00000000e+00)};

\end{axis}
\end{tikzpicture}
 \caption{Geometry and boundary conditions for the cylinder test case.
 Boundary conditions: slip walls (\ref{line:cyl0:wall}), Mach $M_\infty$
 supersonic inflow (\ref{line:cyl0:farf}), and supersonic outflow
 (\ref{line:cyl0:outf}).}
 \label{fig:cyl0_geom}
\end{figure}
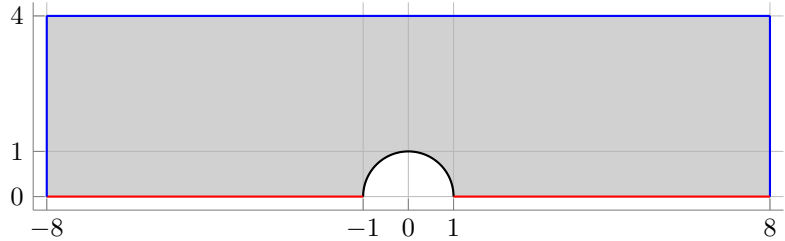

\begin{figure}
 \centering
 \includegraphics[width=0.124\textwidth,angle=-90]{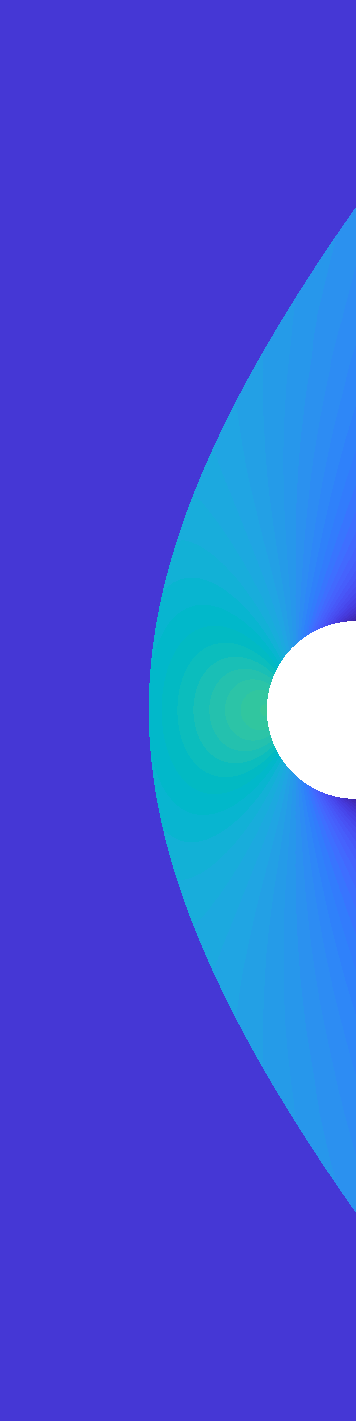}
 \includegraphics[width=0.124\textwidth,angle=-90]{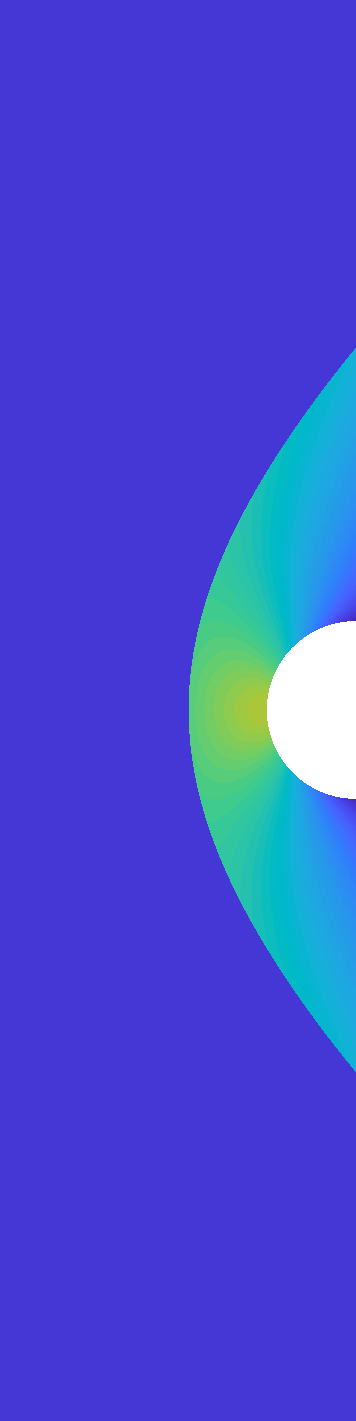} \\
 \includegraphics[width=0.124\textwidth,angle=-90]{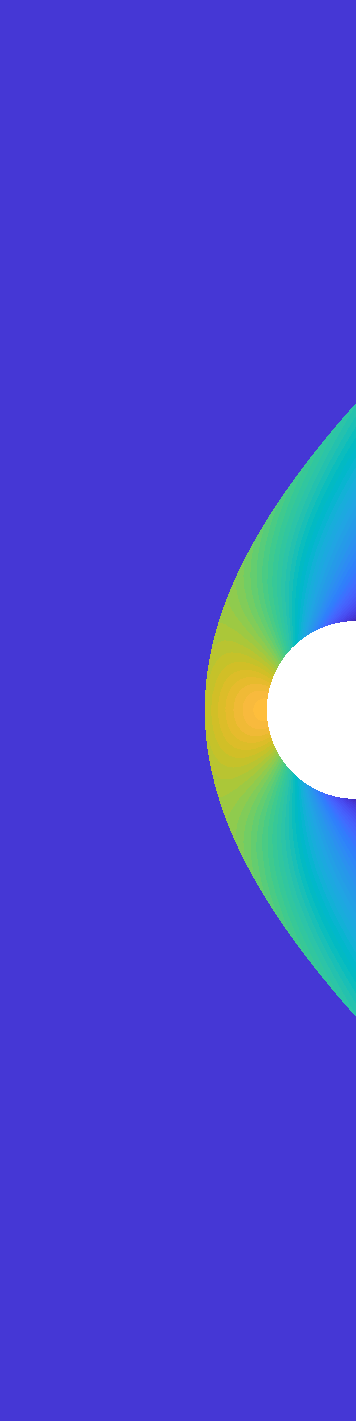}
 \includegraphics[width=0.124\textwidth,angle=-90]{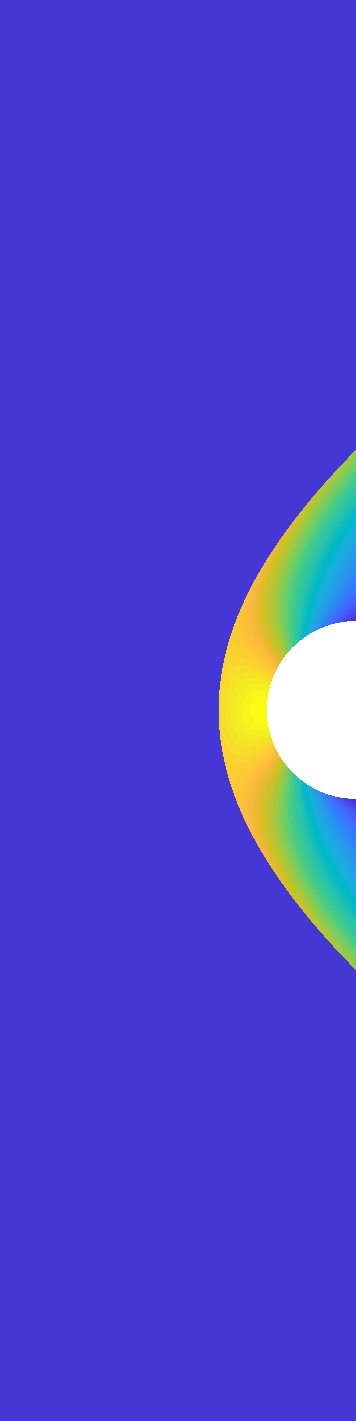}

\colorbarMatlabParula{1}{2.5}{4}{5.5}{7}
\caption{Solution (density) of inviscid flow over cylinder for
 $M_\infty=2$ (\textit{top left}), $M_\infty=2.51$ (\textit{top right}),
 $M_\infty=3$ (\textit{bottom left}), and $M_\infty=4$ (\textit{bottom right}).}
\label{fig:cyl0_hdm}
\end{figure}

Finally, we consider inviscid, supersonic flow of an ideal gas over a
cylinder. Consider a conservation law of the form (\ref{eqn:claw-phys})
defined over the domain $\Omega\subset\Rbb^d$ ($d=2$ in this problem)
with solution $\func{u}{\Omega\times\Dcal}{\Rbb^{d+2}}$ with components
\begin{equation}\label{eqn:euler-stvc}
u = \begin{bmatrix} \rho \\ \rho v \\ \rho E \end{bmatrix},
\end{equation}
where
$\func{\rho}{\Omega\times\Dcal}{\Rbb_{>0}}$,
$\func{v}{\Omega\times\Dcal}{\Rbb^d}$,
$\func{E}{\Omega\times\Dcal}{\Rbb_{>0}}$
are the density, velocity, and energy of the fluid, respectively.
The flux function and source term are defined as
\begin{equation} \label{eqn:euler}
 f(u) \coloneqq \begin{bmatrix} \rho v^T \\ \rho v\otimes v+P(u)I_d \\ (\rho E+P(u))v^T \end{bmatrix},
 \qquad
 s(u) \coloneqq 0,
\end{equation}
where $u\mapsto P(u) \coloneqq (\gamma-1)(\rho E - \rho \norm{v}^2/2)$ is the
pressure of the fluid and $\gamma\in\Rbb_{>0}$ is the ratio of specific
heats. In this problem, we take the domain $\Omega$ to be a cylinder
with freestream density $\rho_\infty=1$, velocity
$v_\infty = (0,-M_\infty c_\infty)$, pressure $P_\infty=1$,
and sound speed $c_\infty=\sqrt{\gamma}$ (Figure~\ref{fig:cyl0_geom}).
The freestream Mach number $M_\infty\in[2,4]$ is used to
parametrize the problem. The solution has a bow shock
that strengthens and moves closer to the cylinder surface as the freestream
Mach number increases (Figure~\ref{fig:cyl0_hdm}). We order all training
sets $\Dcal_\mathrm{tr}\subset\Dcal$ such that the centroid of the parameter
space ($M_\infty=3$) is the first parameter.

The system of conservation laws is discretized using the HOIST method
using $400$ elements of polynomial degree $p=2$ for a total of $N=9600$
degrees of freedom.
The Roe flux with Harten-Hyman entropy fix is used for the inviscid
numerical flux. We also use a high-order domain mapping $q=2$ to accurately
represent the curvature of the bow shock, which defines $\Gbb_{h,q}^\text{b}$
as the collection of piecewise quadratic functions that preserve the
boundaries of the domain in Figure~\ref{fig:cyl0_geom}.
We choose the reference domain to be identical to the
physical domain $\rdom=\pdom$ and the nominal mapping to be identity
$\bar\Gcal=\mathrm{Id}$. The HOIST method is used to define a mesh of
the domain $\rdom$ that conforms to the bow shock at $M_\infty=3$
(the first parameter of any training set).

\begin{figure}[t]
 \centering
 \includegraphics[width=0.124\textwidth,angle=-90]{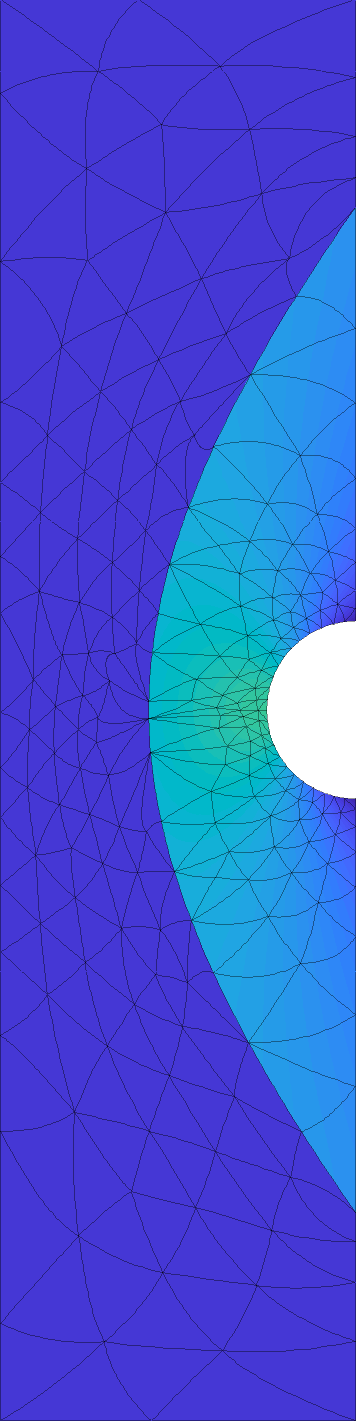}
 \includegraphics[width=0.124\textwidth,angle=-90]{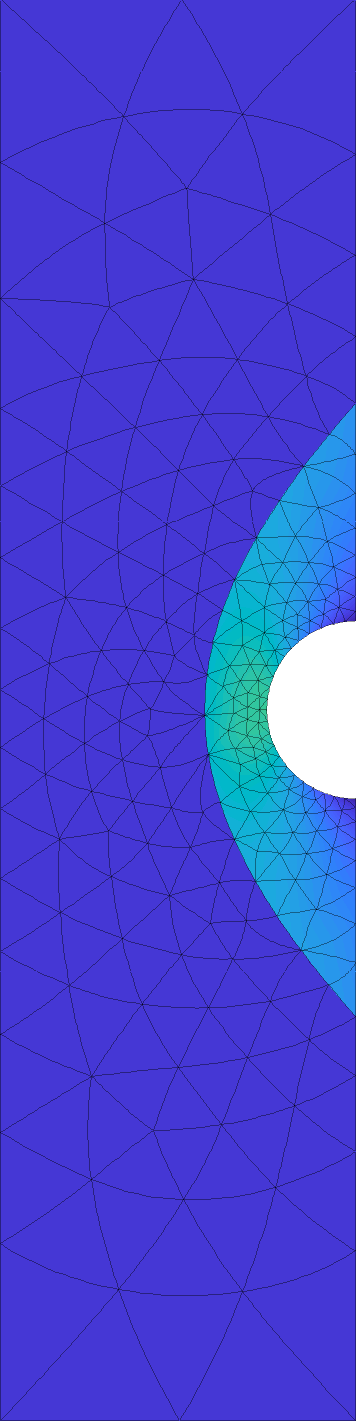} \\
 \includegraphics[width=0.124\textwidth,angle=-90]{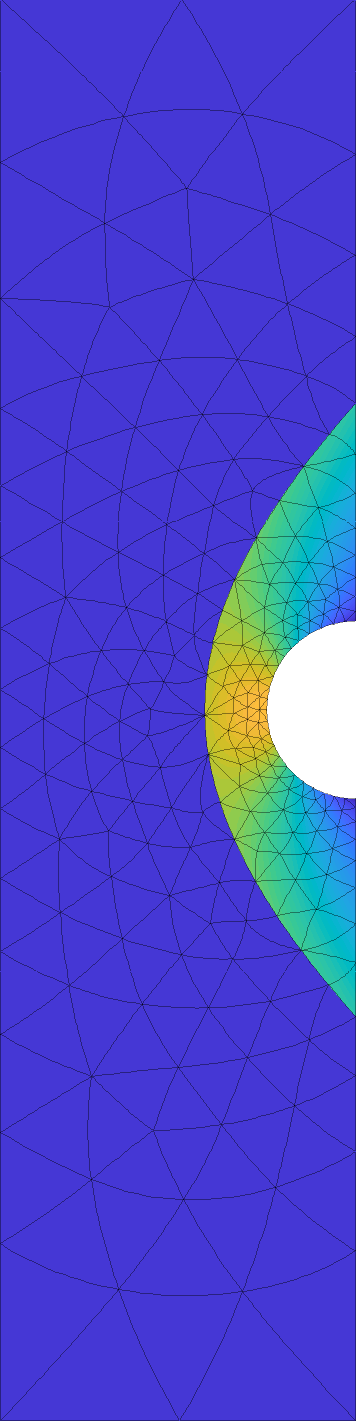}
 \includegraphics[width=0.124\textwidth,angle=-90]{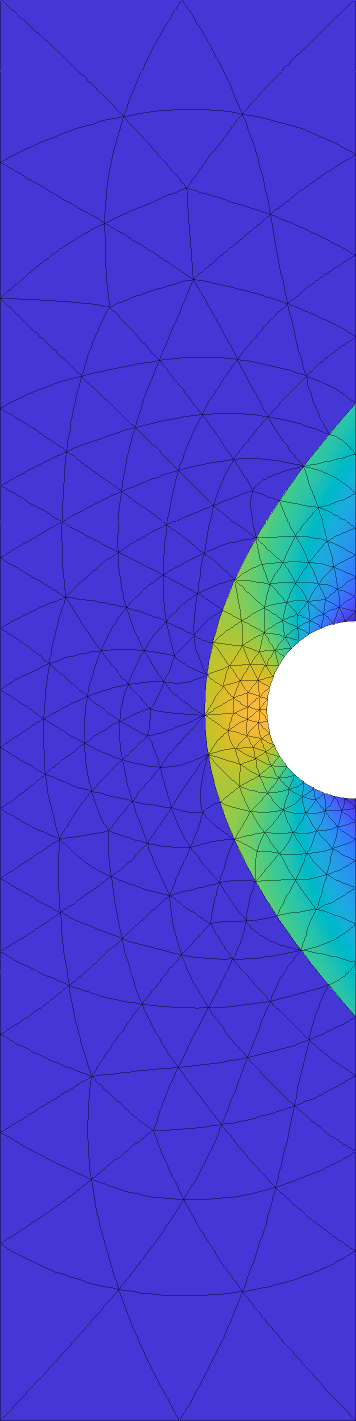} \\
 \includegraphics[width=0.124\textwidth,angle=-90]{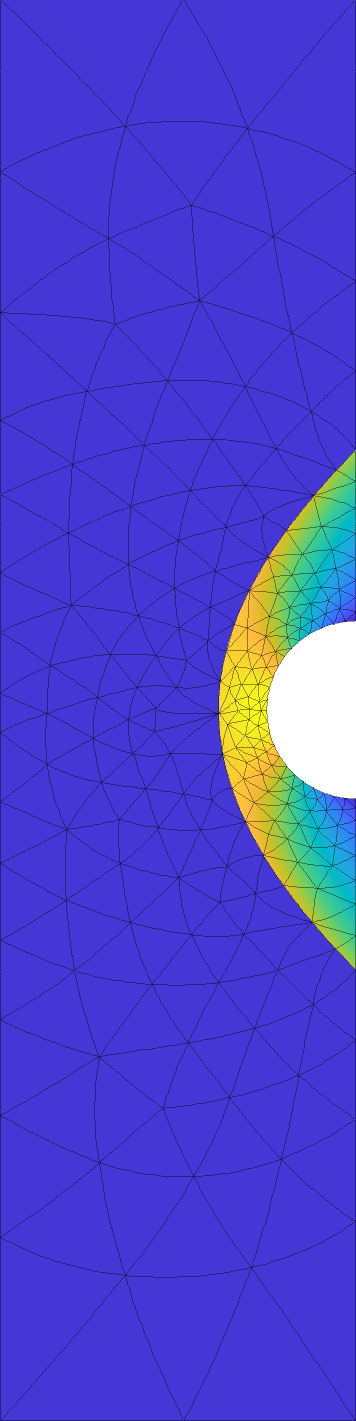}
 \includegraphics[width=0.124\textwidth,angle=-90]{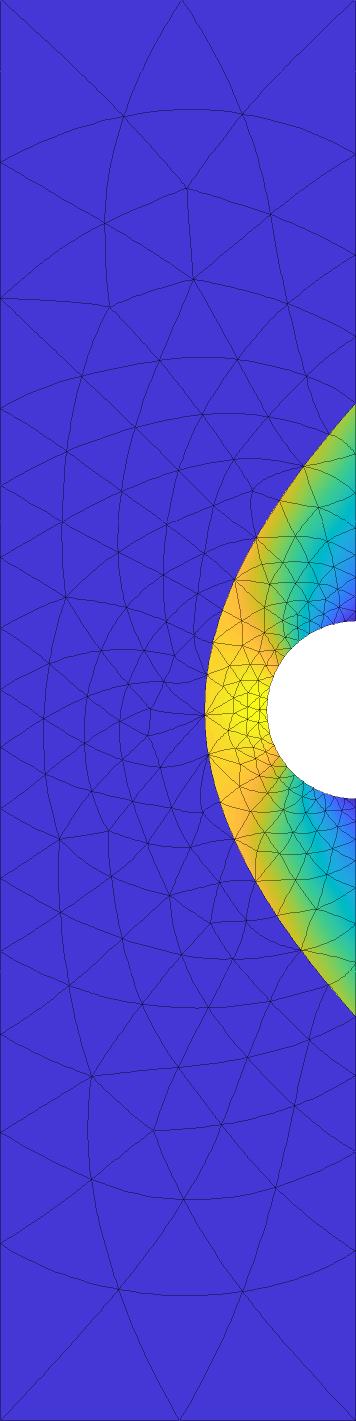}
\caption{Snapshots generated by the ROM-IFT method corresponding to the
 parameter set $\Dcal_3$. The left column shows the snapshot and corresponding
 domain deformation (mesh edges) in the physical domain, while the right
 column shows the corresponding snapshot in the fixed reference domain.
 Colorbar in Figure~\ref{fig:cyl0_hdm}.}
 \label{fig:cyl0_ift_snaps}
\end{figure}

\begin{figure}[t]
\centering
\input{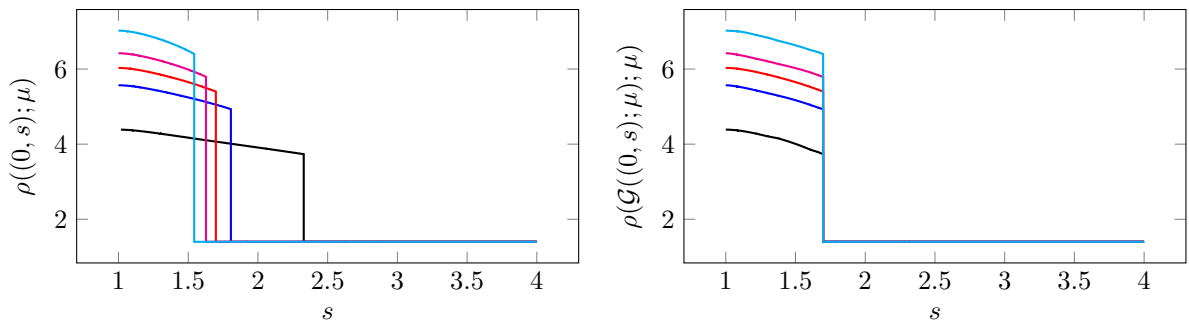}
\caption{Five snapshots of the flow over the cylinder (density) along the
 curve $\Gamma=\{(0,s)\mid s\in[1,4]\}$ associated with the training set
 $\Dcal_5$ in the physical domain (\textit{left}) and reference domain
 (\textit{right}).}
 \label{fig:cyl0_ift_snaps_slice}
\end{figure}

\begin{figure}[t]
 \centering
 \includegraphics[width=0.124\textwidth,angle=-90]{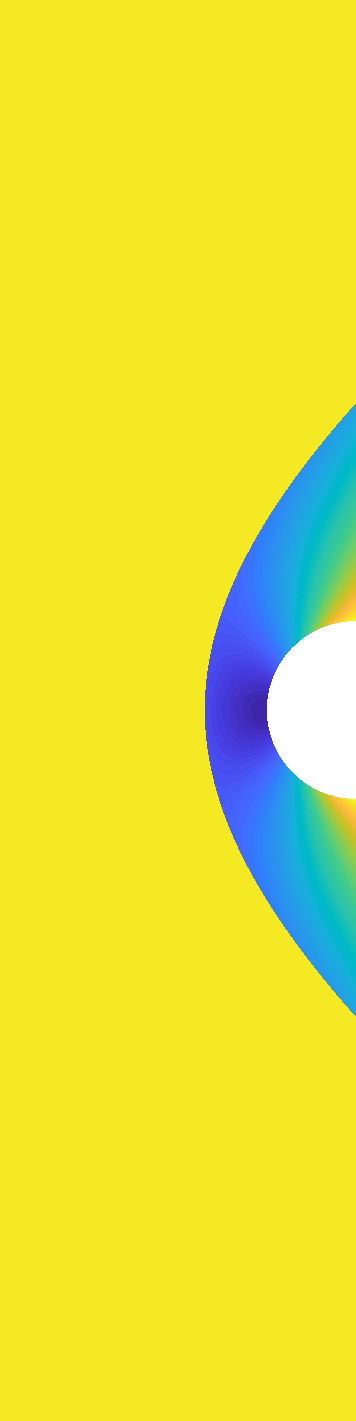}
 \includegraphics[width=0.124\textwidth,angle=-90]{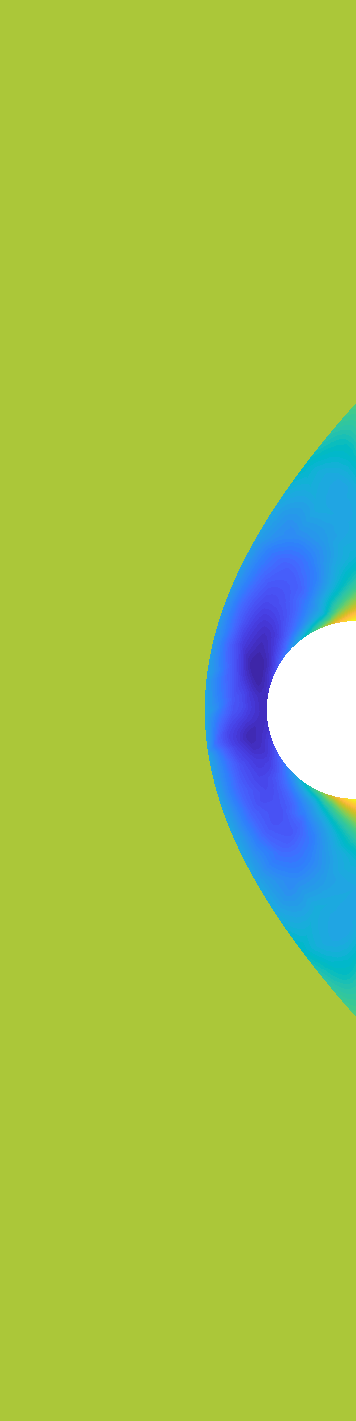} \\
 \includegraphics[width=0.124\textwidth,angle=-90]{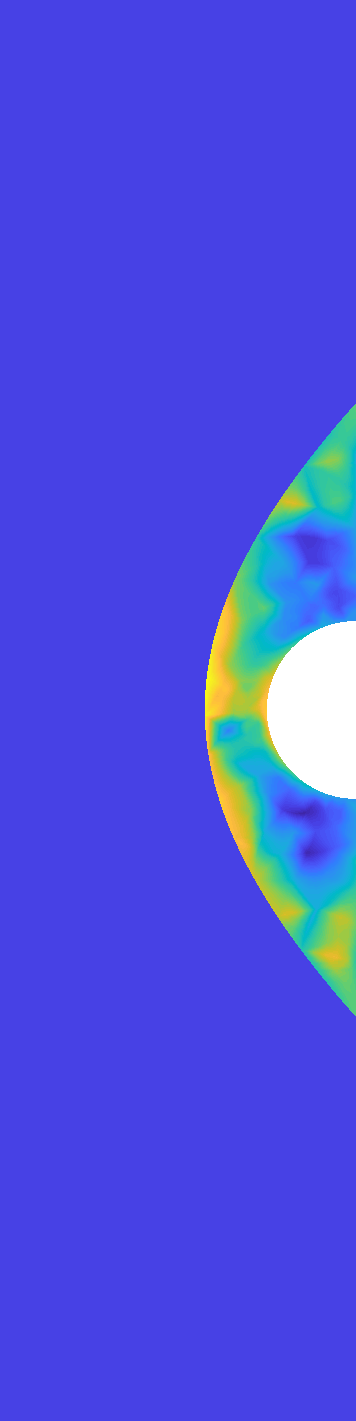}
 \includegraphics[width=0.124\textwidth,angle=-90]{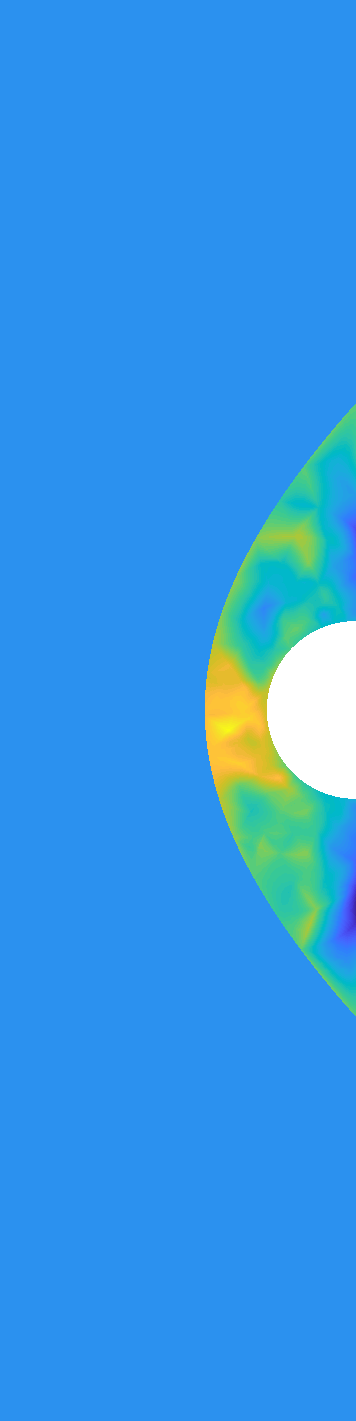}
\caption{First four POD modes (\textit{left-to-right}, \textit{top-to-bottom})
 corresponding to the aligned snapshots generated from the training set
 $\Dcal_5$. Each figure uses a different colorbar, scaled to its range
 to highlight features in the corresponding basis function.}
 \label{fig:cyl0_ift_pod}
\end{figure}

For this problem, we consider a targeted study to verify the
high accuracy of the ROM-IFT under limited training. In particular,
we consider two training sets, one with three parameters
$\Dcal_3 = \{2,3,4\}$ and one with five parameters
$\Dcal_5 = \{2,2.5,3,3.5,4\}$, that will be used to train the
ROM-IFT method; we will test the predictive accuracy on a set
containing $101$ uniformly spaced Mach numbers in the range
$[2,4]$ ($\Dcal_{101}$). Due to the limited amount of training,
we do not truncate the POD basis, which leads to a basis of
dimension $3$ and $5$, respectively. The proposed offline
procedure successfully determines a domain deformation for
each snapshot that causes the bow shock to align in the
reference domain (across all snapshots), despite being
far from aligned in the physical domain. This can be
seen from the full two-dimensional density field
(Figure~\ref{fig:cyl0_ift_snaps}) and its projection along the centerline
of the domain, $\Gamma = \{(0, s) \mid s\in[1, 4]\}$
(Figure~\ref{fig:cyl0_ift_snaps_slice}). Due to the alignment of the
snapshots, the location of the bow shock is independent of
Mach number in the reference domain, which makes compression
via POD highly effective. The first POD mode captures the
dominate behavior of the solution and subsequent modes provide
additional resolution in the smooth shock layer
(Figure~\ref{fig:cyl0_ift_pod}).

\begin{figure}[t]
\centering
\input{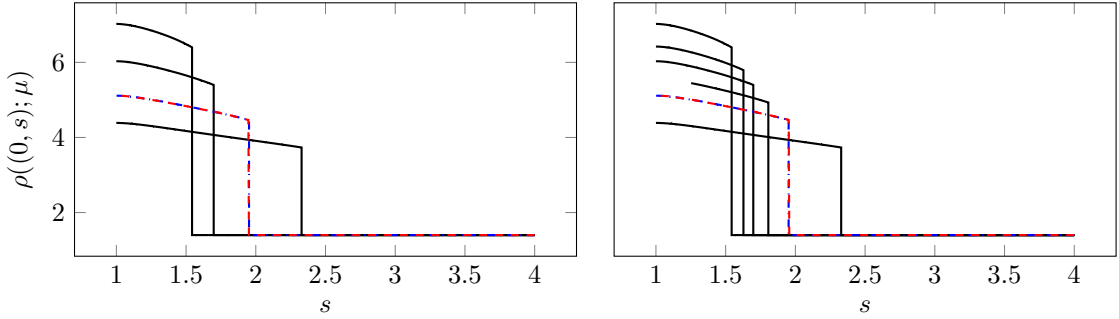}
\caption{Comparison of the HDM and ROM-IFT on the cylinder problem (density)
     along the line $\Gamma=\{(0,s)\mid s\in[1,4]\}$ with three (\textit{left})
     and five (\textit{right}) training parameters. The test set includes $101$
     Mach numbers ($\Dcal_{101}$), although only the Mach number that leads to
     a shock farthest from those in the training sets is shown
     ($M_\infty=2.36$). Legend: training snapshots
     (\ref{line:eulerii2d_cyl0_states_snaps}), HDM solution at test parameters
     (\ref{line:eulerii2d_cyl0_states_hdm}), and ROM-IFT solution at test
     parameter (\ref{line:eulerii2d_cyl0_states_ift}).}
 \label{fig:cyl0_ift_online}
\end{figure}

\begin{figure}[t]
\centering
\input{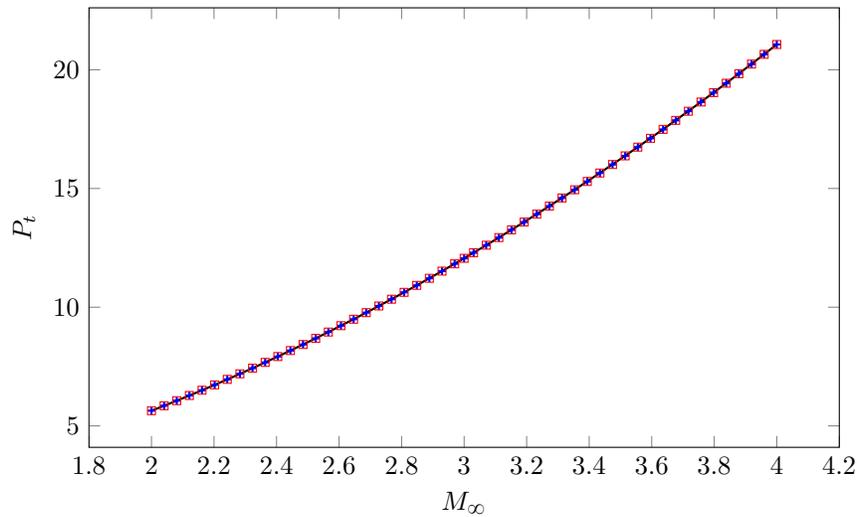}
\caption{Comparison of the analytical stagnation pressure
  (\ref{line:stagpres_exact}) as a function of freestream Mach number
  to the stagnation pressure predicted by the ROM-IFT method with training
  set $\Dcal_3$ (\ref{line:stagpres_ift3}) and $\Dcal_5$
  (\ref{line:stagpres_ift5}).}
\label{fig:cyl0_ift_stagpres}
\end{figure}

We assess the online performance of the ROM by considering the density
field projected along the centerline ($\Gamma$)
(Figure~\ref{fig:cyl0_ift_online}). In addition, we compare the stagnation
pressure $P_t$, i.e., the pressure at the stagnation point $x=(0,1)$,
predicted by the ROM-IFT method to the known value determined from
compressible flow theory (Figure~\ref{fig:cyl0_ift_stagpres})
\begin{equation}
 P_t = P_\infty \left(\frac{1-\gamma+2\gamma M_\infty^2}{\gamma+1}\right)\left(\frac{(\gamma+1)^2M_\infty^2}{4\gamma M_\infty^2-2(\gamma-1)}\right)^{\gamma/(\gamma-1)}.
\end{equation}
The prediction of the density field using the ROM-IFT method is highly
accurate using either training set, even at a Mach number that leads to
a bow shock far from those seen in the training
(Figure~\ref{fig:cyl0_ift_online}).
Similarly, we see the stagnation pressure is predicted very well using
the ROM-IFT method using either training set. The maximum relative
stagnation pressure error over the training set $\Dcal_3$ ($\Dcal_5$)
is $0.031\%$ ($0.020\%$), suggesting there is no reason to use more than
three training parameters and the ROM-IFT method is accurate with little
training.

\section{Conclusion}
\label{sec:concl}
This work introduced a new approach to reduce the cost of solving PDEs with
convection-dominated solutions: \textit{model reduction with implicit
feature tracking}. Traditional model reduction techniques that use
an affine subspace to reduce the dimensionality of the solution manifold
yield limited reduction and require extensive training due to the
slowly decaying Kolmogorov $n$-width of convection-dominated problems.
The ROM-IFT approach circumvents the slowly decaying $n$-width limitation
by using a nonlinear approximation manifold systematically defined by
composing a low-dimensional affine space with a space of bijections
of the underlying domain, similar to the approaches in
\cite{welper_interpolation_2017,cagniart_model_2019,nair_transported_2019,taddei_registration_2020}.
The novelty of the ROM-IFT method lies in 1) the simultaneous determination
of the reduced coordinates of the affine subspace and domain mapping via
residual minimization (Section~\ref{sec:rom:ift},~\ref{sec:rom:on}),
2) the general family of boundary-preserving domain mappings constructed
via high-order finite elements (Section~\ref{sec:govern:dommap})
\cite{huang2021robust} and its compression to a low-dimensional manifold
for efficient online computations (Section~\ref{sec:rom:off}),
3) an offline training phase that bootstraps the online ROM-IFT
formulation (Section~\ref{sec:rom:off}), and 4) a detailed description
of a Levenberg-Marquardt solver for the ROM-IFT residual minimization
problem (Section~\ref{sec:solver}). Several numerical
experiments demonstrate the offline ROM-IFT method effectively
aligns features in snapshots (Figure~\ref{fig:advec-react-snap},~\ref{fig:eulernozii1d_stdy1_snaps_basis},~\ref{fig:cyl0_ift_snaps_slice}),
which leads to highly effective compression of the snapshots into a reduced
basis (Figure~\ref{fig:advec-react-study0-cnvg},~\ref{fig:eulernozii1d_stdy0_cnvg}).
As a result, the ROM-IFT method leads to accurate predictions with limited
training (Figure~\ref{fig:eulernozii1d_stdy1_states},~\ref{fig:cyl0_ift_online}).
For the supersonic cylinder problem
that has a parametrized bow shock, only three training parameters are
needed to predict the stagnation pressure to within $0.03\%$ of the
exact solution over the entire parameter domain of interest
(Section~\ref{sec:numexp:euler-cyl0}).

Additional research is required to incorporate hyperreduction techniques
into the ROM-IFT framework to realize computational efficiency. Furthermore,
development of an efficient greedy sampling procedure would eliminate the
need to uniformly sample the parameter space \textit{a priori}, which can
be expensive, especially when many parameters are considered.


\section*{Acknowledgments}
This material is based upon work supported by the Air Force Office of
Scientific Research (AFOSR) under award number FA9550-20-1-0236. The
content of this publication does not necessarily reflect the position
or policy of any of these supporters, and no official endorsement
should be inferred.

\bibliographystyle{plain}
\bibliography{biblio}

\begin{thebibliography}{10}

\bibitem{agarwal_trust-region_2013}
Anshul Agarwal and Lorenz~T. Biegler.
\newblock A trust-region framework for constrained optimization using reduced
  order modeling.
\newblock {\em Optimization and Engineering}, 14(1):3--35, March 2013.

\bibitem{amsallem_nonlinear_2012}
David Amsallem, Matthew~J. Zahr, and Charbel Farhat.
\newblock Nonlinear model order reduction based on local reduced-order bases.
\newblock {\em International Journal for Numerical Methods in Engineering},
  92(10):891--916, December 2012.

\bibitem{arian_trust-region_2000}
E.~Arian, M.~Fahl, and E.~W. Sachs.
\newblock Trust-{Region} {Proper} {Orthogonal} {Decomposition} for {Flow}
  {Control}.
\newblock Technical Report ICASE-2000-25, Institue for Computer Applications in
  Science and Engineering, May 2000.

\bibitem{bansal2021model}
Harshit Bansal, Stephan Rave, Laura Iapichino, W~Schilders, and N~Wouw.
\newblock Model order reduction framework for problems with moving
  discontinuities.
\newblock In {\em Numerical Mathematics and Advanced Applications ENUMATH
  2019}, pages 83--91. Springer, 2021.

\bibitem{benner2015survey}
Peter Benner, Serkan Gugercin, and Karen Willcox.
\newblock A survey of projection-based model reduction methods for parametric
  dynamical systems.
\newblock {\em SIAM Review}, 57(4):483--531, 2015.

\bibitem{black_projection-based_2020}
Felix Black, Philipp Schulze, and Benjamin Unger.
\newblock Projection-based model reduction with dynamically transformed modes.
\newblock {\em ESAIM: Mathematical Modelling and Numerical Analysis},
  54(6):2011--2043, 2020.

\bibitem{blonigan2021model}
Patrick Blonigan, Francesco Rizzi, Micah Howard, Jeffrey Fike, and Kevin
  Carlberg.
\newblock Model reduction for steady hypersonic aerodynamics via conservative
  manifold least-squares {P}etrov--{G}alerkin projection.
\newblock {\em AIAA Journal}, 59(4):1296--1312, 2021.

\bibitem{bonfiglioli2014convergence}
Aldo Bonfiglioli and Renato Paciorri.
\newblock Convergence analysis of shock-capturing and shock-fitting solutions
  on unstructured grids.
\newblock {\em AIAA Journal}, 52(7):1404--1416, 2014.

\bibitem{brooks_karhunenloeve_2004}
Gregory~P. Brooks and Joseph~M. Powers.
\newblock A {Karhunen}--{Lo{\`e}ve} least-squares technique for optimization of
  geometry of a blunt body in supersonic flow.
\newblock {\em Journal of Computational Physics}, 195(1):387--412, March 2004.

\bibitem{budd2009adaptivity}
Chris Budd, Weizhang Huang, and Robert Russell.
\newblock Adaptivity with moving grids.
\newblock {\em Acta Numerica}, 18:111--241, 2009.

\bibitem{buffa_priori_2012}
Annalisa Buffa, Yvon Maday, Anthony~T. Patera, Christophe Prud'homme, and
  Gabriel Turinici.
\newblock \textit{{A} priori} convergence of the {Greedy} algorithm for the
  parametrized reduced basis method.
\newblock {\em ESAIM: Mathematical Modelling and Numerical Analysis},
  46(3):595--603, May 2012.

\bibitem{cagniart_model_2019}
Nicolas Cagniart, Yvon Maday, and Benjamin Stamm.
\newblock Model order reduction for problems with large convection effects.
\newblock In B.~N. Chetverushkin, W.~Fitzgibbon, Y.A. Kuznetsov,
  P.~Neittaanm{\"a}ki, J.~Periaux, and O.~Pironneau, editors, {\em
  Contributions to {Partial} {Differential} {Equations} and {Applications}},
  Computational {Methods} in {Applied} {Sciences}, pages 131--150. Springer
  International Publishing, Cham, 2019.

\bibitem{carlberg_adaptive_2015}
Kevin Carlberg.
\newblock Adaptive h-refinement for reduced-order models.
\newblock {\em International Journal for Numerical Methods in Engineering},
  102(5):1192--1210, 2015.

\bibitem{carlberg_efficient_2011}
Kevin Carlberg, Charbel Bou-Mosleh, and Charbel Farhat.
\newblock Efficient non-linear model reduction via a least-squares
  {Petrov}-{Galerkin} projection and compressive tensor approximations.
\newblock {\em International Journal for Numerical Methods in Engineering},
  86(2):155--181, April 2011.

\bibitem{constantine_reduced_2012}
Paul Constantine and Gianluca Iaccarino.
\newblock Reduced order models for parameterized hyperbolic conservations laws
  with shock reconstruction.
\newblock Technical report, Center for Turbulence Research Annual Brief, 2012.

\bibitem{corrigan_moving_2019}
Andrew Corrigan, Andrew~D. Kercher, and David~A. Kessler.
\newblock A moving discontinuous {Galerkin} finite element method for flows
  with interfaces: {A} moving discontinuous {Galerkin} finite element method
  for flows with interfaces.
\newblock {\em International Journal for Numerical Methods in Fluids},
  89(9):362--406, March 2019.

\bibitem{dihlmann_model_2011}
Markus Dihlmann, Martin Drohmann, and Bernard Haasdonk.
\newblock Model reduction of parametrized evolution problems using the reduced
  basis method with adaptive time-partitioning.
\newblock In {\em {Proceedings} of {ADMOS}}, page~64, 2011.

\bibitem{eftang_hp_2010}
Jens~L. Eftang, Anthony~T. Patera, and Einar~M. R{\o}nquist.
\newblock An $hp$ certified reduced basis method for parametrized elliptic
  partial differential equations.
\newblock {\em SIAM Journal on Scientific Computing}, 32(6):3170--3200, January
  2010.

\bibitem{ferrero2021registration}
Andrea Ferrero, Tommaso Taddei, and Lei Zhang.
\newblock Registration-based model reduction of parameterized two-dimensional
  conservation laws.
\newblock {\em arXiv preprint arXiv:2105.02024}, 2021.

\bibitem{gerbeau_approximated_2014}
Jean-Fr{\'e}d{\'e}ric Gerbeau and Damiano Lombardi.
\newblock Approximated {Lax} pairs for the reduced order integration of
  nonlinear evolution equations.
\newblock {\em Journal of Computational Physics}, 265:246--269, May 2014.

\bibitem{greif_decay_2019}
Constantin Greif and Karsten Urban.
\newblock Decay of the {Kolmogorov} {N}-width for wave problems.
\newblock {\em Applied Mathematics Letters}, 96:216--222, October 2019.

\bibitem{grepl2005posteriori}
Martin~A Grepl and Anthony~T Patera.
\newblock A posteriori error bounds for reduced-basis approximations of
  parametrized parabolic partial differential equations.
\newblock {\em ESAIM: Mathematical Modelling and Numerical Analysis},
  39(1):157--181, 2005.

\bibitem{haasdonk_training_2011}
Bernard Haasdonk, Markus Dihlmann, and Mario Ohlberger.
\newblock A training set and multiple bases generation approach for
  parameterized model reduction based on adaptive grids in parameter space.
\newblock {\em Mathematical and Computer Modelling of Dynamical Systems},
  17(4):423--442, August 2011.

\bibitem{harten1983self}
Ami Harten and James Hyman.
\newblock Self adjusting grid methods for one-dimensional hyperbolic
  conservation laws.
\newblock {\em Journal of Computational Physics}, 50(2):235--269, 1983.

\bibitem{hartman_deep_2017}
David Hartman and Lalit~K. Mestha.
\newblock A deep learning framework for model reduction of dynamical systems.
\newblock In {\em 2017 {IEEE} {Conference} on {Control} {Technology} and
  {Applications} ({CCTA})}, pages 1917--1922, August 2017.

\bibitem{holmes1996turbulence}
Philip Holmes, John Lumley, and Gahl Berkooz.
\newblock {\em Turbulence, coherent structures, dynamical systems and
  symmetry}.
\newblock Cambridge University Press, Cambridge, UK, 1996.

\bibitem{huang2021robust}
Tianci Huang and Matthew~J Zahr.
\newblock A robust, high-order implicit shock tracking method for simulation of
  complex, high-speed flows.
\newblock {\em arXiv preprint arXiv:2105.00139}, 2021.

\bibitem{kashima_nonlinear_2016}
Kenji Kashima.
\newblock Nonlinear model reduction by deep autoencoder of noise response data.
\newblock In {\em 2016 {IEEE} 55th {Conference} on {Decision} and {Control}
  ({CDC})}, pages 5750--5755, Las Vegas, NV, USA, December 2016. IEEE.

\bibitem{kim2020efficient}
Youngkyu Kim, Youngsoo Choi, David Widemann, and Tarek Zohdi.
\newblock Efficient nonlinear manifold reduced order model.
\newblock {\em arXiv preprint arXiv:2011.07727}, 2020.

\bibitem{knupp2001algebraic}
Patrick Knupp.
\newblock Algebraic mesh quality metrics.
\newblock {\em SIAM Journal on Scientific Computing}, 23(1):193--218, 2001.

\bibitem{lee_model_2020}
Kookjin Lee and Kevin~T. Carlberg.
\newblock Model reduction of dynamical systems on nonlinear manifolds using
  deep convolutional autoencoders.
\newblock {\em Journal of Computational Physics}, 404:108973, March 2020.

\bibitem{lee1999spurious}
Theodore~Kai Lee and Xiaolin Zhong.
\newblock Spurious numerical oscillations in simulation of supersonic flows
  using shock-capturing schemes.
\newblock {\em AIAA Journal}, 37(3):313--319, 1999.

\bibitem{legresley2006application}
Patrick~Allen LeGresley.
\newblock {\em Application of Proper Orthogonal Decomposition ({POD}) to Design
  Decomposition Methods}.
\newblock Stanford University, 2006.

\bibitem{lucia_reduced_2003}
David~J. Lucia, Paul~I. King, and Philip~S. Beran.
\newblock Reduced order modeling of a two-dimensional flow with moving shocks.
\newblock {\em Computers \& Fluids}, 32(7):917--938, August 2003.

\bibitem{maday_blackbox_2002}
Y.~Maday, A.~T. Patera, and D.~V. Rovas.
\newblock A blackbox reduced-basis output bound method for noncoercive linear
  problems.
\newblock In {\em Nonlinear Partial Differential Equations and their
  Applications - Coll{\`{e}}ge de France Seminar Volume {XIV}}, pages 533--569.
  Elsevier, 2002.

\bibitem{maulik_reduced-order_2020}
Romit Maulik, Bethany Lusch, and Prasanna Balaprakash.
\newblock Reduced-order modeling of advection-dominated systems with recurrent
  neural networks and convolutional autoencoders.
\newblock {\em Physics of Fluids}, 33(3):037106, 2021.

\bibitem{mojgani_lagrangian_2017}
Rambod Mojgani and Maciej Balajewicz.
\newblock Lagrangian basis method for dimensionality reduction of convection
  dominated nonlinear flows.
\newblock {\em arXiv:1701.04343 [physics]}, January 2017.
\newblock arXiv: 1701.04343.

\bibitem{mojgani_physics-aware_2020}
Rambod Mojgani and Maciej Balajewicz.
\newblock Physics-aware registration based auto-encoder for convection
  dominated {PDEs}.
\newblock {\em arXiv:2006.15655 [cs, math]}, June 2020.
\newblock arXiv: 2006.15655.

\bibitem{moretti1987computation}
Gino Moretti.
\newblock Computation of flows with shocks.
\newblock {\em Annual Review of Fluid Mechanics}, 19(1):313--337, 1987.

\bibitem{moretti2002thirty}
Gino Moretti.
\newblock Thirty-six years of shock fitting.
\newblock {\em Computers \& Fluids}, 31(4-7):719--723, 2002.

\bibitem{nair_transported_2019}
Nirmal~J. Nair and Maciej Balajewicz.
\newblock Transported snapshot model order reduction approach for parametric,
  steady-state fluid flows containing parameter-dependent shocks: {Model} order
  reduction for fluid flows containing shocks.
\newblock {\em International Journal for Numerical Methods in Engineering},
  117(12):1234--1262, March 2019.

\bibitem{nocedal_numerical_2006}
Jorge Nocedal and Stephen~J. Wright.
\newblock {\em Numerical Optimization}.
\newblock Springer Series in Operations Research. Springer, New York, 2006.

\bibitem{ohlberger_nonlinear_2013}
Mario Ohlberger and Stephan Rave.
\newblock Nonlinear reduced basis approximation of parameterized evolution
  equations via the method of freezing.
\newblock {\em Comptes Rendus Mathematique}, 351(23-24):901--906, December
  2013.

\bibitem{ohlberger_reduced_2016}
Mario Ohlberger and Stephan Rave.
\newblock Reduced {Basis} {Methods}: {Success}, {Limitations} and {Future}
  {Challenges}.
\newblock In {\em {Proceedings} of the {Conference} {Algoritmy}}, 2016.

\bibitem{omata_novel_2019}
Noriyasu Omata and Susumu Shirayama.
\newblock A novel method of low-dimensional representation for temporal
  behavior of flow fields using deep autoencoder.
\newblock {\em AIP Advances}, 9(1):015006, January 2019.

\bibitem{peherstorfer_model_2020}
Benjamin Peherstorfer.
\newblock Model reduction for transport-dominated problems via online adaptive
  bases and adaptive sampling.
\newblock {\em SIAM Journal on Scientific Computing}, 42(5):A2803--A2836,
  January 2020.

\bibitem{peraire2008compact}
Jaime Peraire and Per-Olof Persson.
\newblock The compact discontinuous {G}alerkin ({CDG}) method for elliptic
  problems.
\newblock {\em SIAM Journal on Scientific Computing}, 30(4):1806--1824, 2008.

\bibitem{persson2009discontinuous}
Per-Olof Persson, Javier Bonet, and Jaime Peraire.
\newblock {D}iscontinuous {G}alerkin solution of the {N}avier--{S}tokes
  equations on deformable domains.
\newblock {\em Computer Methods in Applied Mechanics and Engineering},
  198(17-20):1585--1595, 2009.

\bibitem{persson2006sub}
Per-Olof Persson and Jaime Peraire.
\newblock Sub-cell shock capturing for discontinuous {G}alerkin methods.
\newblock In {\em 44th AIAA Aerospace Sciences Meeting and Exhibit}, 2006.

\bibitem{powers2006exact}
Joseph~M. Powers and Tariq~D. Aslam.
\newblock Exact solution for multidimensional compressible reactive flow for
  verifying numerical algorithms.
\newblock {\em AIAA Journal}, 44(2):337--344, 2006.

\bibitem{prud2002reliable}
Christophe Prud'Homme, Dimitrios Rovas, Karen Veroy, Luc Machiels, Yvon Maday,
  Anthony Patera, and Gabriel Turinici.
\newblock Reliable real-time solution of parametrized partial differential
  equations: {R}educed-basis output bound methods.
\newblock {\em Journal of Fluids Engineering}, 124(1):70--80, 2002.

\bibitem{reiss_shifted_2018}
J.~Reiss, P.~Schulze, J.~Sesterhenn, and V.~Mehrmann.
\newblock The {Shifted} {Proper} {Orthogonal} {Decomposition}: A mode
  decomposition for multiple transport phenomena.
\newblock {\em SIAM Journal on Scientific Computing}, 40(3):A1322--A1344,
  January 2018.

\bibitem{rim_displacement_2018}
Donsub Rim and Kyle~T. Mandli.
\newblock Displacement interpolation using monotone rearrangement.
\newblock {\em SIAM/ASA Journal on Uncertainty Quantification},
  6(4):1503--1531, January 2018.

\bibitem{rim_transport_2018}
Donsub Rim, Scott Moe, and Randall~J. LeVeque.
\newblock Transport reversal for model reduction of hyperbolic partial
  differential equations.
\newblock {\em SIAM/ASA Journal on Uncertainty Quantification}, 6(1):118--150,
  January 2018.

\bibitem{rim_manifold_2020}
Donsub Rim, Benjamin Peherstorfer, and Kyle~T. Mandli.
\newblock Manifold approximations via transported subspaces: Model reduction
  for transport-dominated problems.
\newblock {\em arXiv:1912.13024 [cs, math]}, February 2020.
\newblock arXiv: 1912.13024.

\bibitem{roe1981approximate}
Philip~L. Roe.
\newblock Approximate {R}iemann solvers, parameter vectors, and difference
  schemes.
\newblock {\em Journal of Computational Physics}, 43(2):357--372, 1981.

\bibitem{sirovich1987turbulence}
Lawrence Sirovich.
\newblock Turbulence and the dynamics of coherent structures. {I.} {C}oherent
  structures.
\newblock {\em Quarterly of Applied Mathematics}, 45(3):561--571, 1987.

\bibitem{taddei_registration_2020}
Tommaso Taddei.
\newblock A registration method for model order reduction: Data compression and
  geometry reduction.
\newblock {\em SIAM Journal on Scientific Computing}, 42(2):A997--A1027,
  January 2020.

\bibitem{taddei_reduced_2015}
Tommaso Taddei, Simona Perotto, and Alfio Quarteroni.
\newblock Reduced basis techniques for nonlinear conservation laws.
\newblock {\em ESAIM: Mathematical Modelling and Numerical Analysis},
  49(3):787--814, May 2015.

\bibitem{taddei_space-time_2020}
Tommaso Taddei and Lei Zhang.
\newblock Space-time registration-based model reduction of parameterized
  one-dimensional hyperbolic {PDEs}.
\newblock {\em arXiv:2004.06693 [cs, math]}, April 2020.
\newblock arXiv: 2004.06693.

\bibitem{taddei2021registration}
Tommaso Taddei and Lei Zhang.
\newblock Registration-based model reduction in complex two-dimensional
  geometries.
\newblock {\em arXiv preprint arXiv:2101.10259}, 2021.

\bibitem{torlo_model_2020}
Davide Torlo.
\newblock Model reduction for advection dominated hyperbolic problems in an
  {ALE} framework: Offline and online phases.
\newblock {\em arXiv:2003.13735 [cs, math]}, March 2020.
\newblock arXiv: 2003.13735.

\bibitem{toro2013riemann}
Eleuterio~F Toro.
\newblock {\em Riemann Solvers and Numerical Methods for Fluid Dynamics: A
  Practical Introduction}.
\newblock Springer Science \& Business Media, 2013.

\bibitem{veroy2003posteriori}
Karen Veroy, Christophe Prud'Homme, Dimitrios Rovas, and Anthony Patera.
\newblock A posteriori error bounds for reduced-basis approximation of
  parametrized noncoercive and nonlinear elliptic partial differential
  equations.
\newblock In {\em 16th AIAA Computational Fluid Dynamics Conference}, page
  3847, 2003.

\bibitem{washabaugh_fast_2016}
Kyle~M. Washabaugh.
\newblock {\em Fast Fidelity for Better Design: A Scalable Model Order
  Reduction Framework for Steady Aerodynamic Design Applications}.
\newblock PhD thesis, Stanford University, 2016.

\bibitem{washabaugh_use_2016}
Kyle~M. Washabaugh, Matthew~J. Zahr, and Charbel Farhat.
\newblock On the use of discrete nonlinear reduced-order models for the
  prediction of steady-state flows past parametrically deformed complex
  geometries.
\newblock In {\em 54th {AIAA} {Aerospace} {Sciences} {Meeting}}, San Diego,
  California, USA, January 2016. American Institute of Aeronautics and
  Astronautics.

\bibitem{welper_interpolation_2017}
Garrit Welper.
\newblock Interpolation of functions with parameter dependent jumps by
  transformed snapshots.
\newblock {\em SIAM Journal on Scientific Computing}, 39(4):A1225--A1250,
  January 2017.

\bibitem{welper_transformed_2020}
Garrit Welper.
\newblock Transformed snapshot interpolation with high resolution transforms.
\newblock {\em SIAM Journal on Scientific Computing}, 42(4):A2037--A2061,
  January 2020.

\bibitem{xu_multi-level_2020}
Jiayang Xu and Karthik Duraisamy.
\newblock Multi-level convolutional autoencoder networks for parametric
  prediction of spatio-temporal dynamics.
\newblock {\em Computer Methods in Applied Mechanics and Engineering},
  372:113379, December 2020.

\bibitem{yano_globally_2020}
Masayuki Yano, Tianci Huang, and Matthew~J. Zahr.
\newblock A globally convergent method to accelerate topology optimization
  using on-the-fly model reduction.
\newblock {\em Computer Methods in Applied Mechanics and Engineering},
  375:113635, 2021.

\bibitem{yano_lp_2019}
Masayuki Yano and Anthony~T. Patera.
\newblock An {LP} empirical quadrature procedure for reduced basis treatment of
  parametrized nonlinear {PDEs}.
\newblock {\em Computer Methods in Applied Mechanics and Engineering},
  344:1104--1123, February 2019.

\bibitem{yue_accelerating_2013}
Yao Yue and Karl Meerbergen.
\newblock Accelerating optimization of parametric linear systems by model order
  reduction.
\newblock {\em SIAM Journal on Optimization}, 23(2):1344--1370, January 2013.

\bibitem{zahr_efficient_2019}
Matthew~J. Zahr, Kevin~T. Carlberg, and Drew~P. Kouri.
\newblock An efficient, globally convergent method for optimization under
  uncertainty using adaptive model reduction and sparse grids.
\newblock {\em SIAM/ASA Journal on Uncertainty Quantification}, 7(3):877--912,
  January 2019.

\bibitem{zahr_progressive_2015}
Matthew~J. Zahr and Charbel Farhat.
\newblock Progressive construction of a parametric reduced-order model for
  {PDE}-constrained optimization.
\newblock {\em International Journal for Numerical Methods in Engineering},
  102(5):1111--1135, May 2015.

\bibitem{zahr_optimization-based_2018}
Matthew~J. Zahr and Per-Olof Persson.
\newblock An optimization-based approach for high-order accurate discretization
  of conservation laws with discontinuous solutions.
\newblock {\em Journal of Computational Physics}, 365:105--134, July 2018.

\bibitem{zahr2020r}
Matthew~J. Zahr and Per-Olof Persson.
\newblock An $r$-adaptive, high-order discontinuous {G}alerkin method for flows
  with attached shocks.
\newblock In {\em AIAA Scitech 2020 Forum}, page 0537, 2020.

\bibitem{zahr_implicit_2020}
Matthew~J. Zahr, Andrew Shi, and Per-Olof Persson.
\newblock Implicit shock tracking using an optimization-based high-order
  discontinuous {Galerkin} method.
\newblock {\em Journal of Computational Physics}, 410:109385, June 2020.

\end{thebibliography}

\end{document}